\documentclass[A4j,11pt]{article}
\usepackage{amsfonts,amsmath,amscd}
\usepackage{graphics}
\usepackage{amssymb}
\usepackage[all]{xy}
\begin{document}

\title{{\bf The mean curvature flow\\
for equifocal submanifolds}}
\author{{\bf Naoyuki Koike}}
\date{}
\maketitle

\begin{abstract}
In this paper, we investigate the mean 
curvature flows having an equifocal submanifold in a symmetric space of 
compact type and its focal submanifolds as initial data.  The investigation 
is performed by investigating the lifts of the submanifolds and the flows to 
an (infinite dimensional separable) Hilbert space through a Riemannian 
submersion of the Hilbert space onto the symmetric space.  
\end{abstract}

\vspace{0.5truecm}






\section{Introduction}
The mean curvature flow of a (Riemannian) submanifold $f_0:M\hookrightarrow N$ 
is a map $f:M\times[0,\infty)\to N$ such that, for each $t\in[0,T)$, 
$f_t:M\to N\,(\displaystyle{\mathop{\Leftrightarrow}_{\rm def}}f_t(x)=f(x,t)
\,\,(x\in M))$ is an immersion and 
$f_{\ast}((\frac{\partial}{\partial t})_{(x,t)})$ is the mean curvature vector 
of $f_t:M\hookrightarrow N$, where $T$ is a positive constant or $T=\infty$ 
and $(t)$ is the natural coordinate of $[0,T)$.  Liu-Terng [LT] investigated 
the mean curvature flow having isoparametric submanifolds 
(or their focal submanifolds) in a Euclidean space 
as initial data and obtained the following facts.  

\vspace{0.5truecm}

\noindent
{\bf Fact 1([LT]).} {\sl Let $M$ be a compact isoparametric submanifold in a 
Euclidean space and $C$ be the Weyl domain of $M$ at $x_0\,(\in M)$.  
Then the following statements ${\rm (i)}$ and ${\rm (ii)}$ hold:

${\rm (i)}$ The mean curvature flow $M_t$ having $M$ as initial data converges 
to a focal submanifold of $M$ in finite time.  If a focal map of $M$ onto $F$ 
is spherical, then the mean curvature flow $M_t$ has type I singularity, 
that is, $\lim\limits_{t\to T-0}
{\rm max}_{v\in S^{\perp}M_t}\vert\vert A^t_v\vert\vert^2_{\infty}(T-t)\,<\,
\infty$, where $A^t_v$ is the shape operator of $M_t$ for $v$, 
$\vert\vert A^t_v\vert\vert_{\infty}$ is the sup norm of $A^t_v$ and 
$S^{\perp}M_t$ is the unit normal bundle of $M_t$.

{\rm(ii)} For any focal submanifold $F$ of $M$, there exists a parallel 
submanifold $M'$ of $M$ such that the mean curvature flow having $M'$ 
as initial data converges to $F$ in finite time.}

\vspace{0.5truecm}

\noindent
{\bf Fact 2([LT]).} {\sl Let $M$ and $C$ be as in Fact 1 and $\sigma$ be a 
stratum of dimension greater than zero of $\partial C$.  
Then the following statements ${\rm (i)}$ and ${\rm (ii)}$ hold:

{\rm(i)} For any focal submanifold $F$ (of $M$) through $\sigma$, the maen 
curvature flow $F_t$ having $F$ as initial data converges to a focal 
submanifold $F'$ (of $M$) through $\partial{\sigma}$ in finite time.  
If the fibration of $F$ onto $F'$ is spherical, then the mean curvature flow 
$F_t$ has type I singularity.

(ii) For any focal submanifold $F$ (of $M$) through $\partial\sigma$, there 
exists a focal submanifold $F'$ (of $M$) through $\sigma$ such that the mean 
curvature flow $F'_t$ having $F'$ as initial data converges to $F$ 
in finite time.}

\vspace{0.5truecm}

As a generalized notion of compact isoparametric 
hypersurfaces in a sphere and a hyperbolic space, and a compact isoparametric 
submanifolds in a Euclidean space, Terng-Thorbergsson [TT] defined 
the notion of an equifocal submanifold in a symmetric space as a compact 
submanifold $M$ satisfying the following three conditions:

\vspace{0.2truecm}

(i) the normal holonomy group of $M$ is trivial,

(ii) $M$ has a flat section, that is, for each $x\in M$, 
$\Sigma_x:=\exp^{\perp}(T^{\perp}_xM)$ is totally geodesic and the induced 
metric on $\Sigma_x$ is flat, where $T^{\perp}_xM$ is the normal space of $M$ 
at $x$ and $\exp^{\perp}$ is the normal exponential map of $M$.  

(iii) for each parallel normal vector field $v$ of $M$, the focal radii of $M$ 
along the normal geodesic $\gamma_{v_x}$ (with $\gamma'_{v_x}(0)=v_x$) are 
independent of the choice of $x\in M$, where $\gamma'_{v_x}(0)$ is the 
velocity vector of $\gamma_{v_x}$ at $0$.  

\vspace{0.2truecm}

\noindent
On the other hand, Heintze-Liu-Olmos [HLO] defined the notion of an 
isoparametric submanifold with flat section in a general Riemannian manifold 
as a submanifold $M$ satisfying the above condition (i) and the 
following conditions (ii$'$) and (iii$'$):

\vspace{0.2truecm}

(ii$'$) for each $x\in M$, there exists a neighborhood $U_x$ of 
the zero vector (of $T^{\perp}_xM$) in $T^{\perp}_xM$ such that 
$\Sigma_x:=\exp^{\perp}(U_x)$ is totally geodesic and the induced metric on 
$\Sigma_x$ is flat, 

\vspace{0.2truecm}

(iii$'$) sufficiently close parallel submanifolds of $M$ are CMC with 
respect to the radial direction.  

\vspace{0.2truecm}

\noindent
In the case where the ambient space is a symmetric space $G/K$ of compact 
type, they showed that the notion of an isoparametric submanifold with flat 
section coincides with that of an equifocal submanifold.  
The proof was performed by investigating its lift to $H^0([0,1],\mathfrak g)$ 
through a Riemannian submersion $\pi\circ\phi$, where $\pi$ is the natural 
projection of $G$ onto $G/K$ and $\phi$ is the parallel transport map for 
$G$ (which is a Riemannian submersion of $H^0([0,1],\mathfrak g)$ onto $G$ 
($\mathfrak g:$the Lie algebra of $G$)).  
Let $M$ be an equifocal submanifold in $G/K$ and $v$ be a parallel normal 
vector field of $M$.  The end-point map $\eta_v(:M\mapsto G/K)$ for $v$ 
is defined by $\eta_v(x)=\exp^{\perp}(v_x)$ ($x\in M$).  
Set $M_v:=\eta_v(M)$.  We call $M_v$ 
a parallel submanifold of $M$ when ${\rm dim}\,M_v={\rm dim}\,M$ and a focal 
submanifold of $M$ when ${\rm dim}\,M_v<{\rm dim}\,M$.  The parallel 
submanifolds of $M$ are equifocal.  Let $f:M\times[0,T)\to G/K$ be 
the mean curvature flow having $M$ as initial data.  
Then, it is shown that, for each 
$t\in[0,T)$, $f_t:M\hookrightarrow G/K$ is a parallel submanifold of $M$ 
and hence it is equifocal (see Lemma 3.1).  Fix $x_0\in M$.  Let 
$\widetilde C\,(\subset T^{\perp}_{x_0}M)$ be the fundamental domain 
containing the zero vector (of $T^{\perp}_{x_0}M$) of the Coxeter group 
(which acts on $T^{\perp}_{x_0}M$) of $M$ at $x_0$ 
and set $C:=\exp^{\perp}(\widetilde C)$, where we note that 
$\exp^{\perp}\vert_{\widetilde C}$ is a diffeomorphism onto $C$.  Without loss 
of generality, we may assume that $G$ is simply connected.  
Set $\widetilde M:=(\pi\circ\phi)^{-1}(M)$, which is an isoparametric 
submanifold in $H^0([0,1],\mathfrak g)$.  
Fix $u_0\in(\pi\circ\phi)^{-1}(x_0)$.  
The normal space $T^{\perp}_{x_0}M$ is identified with the normal space 
$T^{\perp}_{u_0}\widetilde M$ of $\widetilde M$ at $u_0$ through 
$(\pi\circ\phi)_{\ast u_0}$.  Each parallel submanifold of $M$ intersects 
with $C$ at the only point and each focal submanifold of $M$ intersects with 
$\partial C$ at the only point, where $\partial C$ is the boundary of $C$.  
Hence, for the mean curvature flow $f:M\times [0,T)\to G/K$ having $M$ as 
initial data, each $M_t(:=f_t(M))$ intersects with $C$ at the only point.  
Denote by $x(t)$ this intersection point and define $u:[0,T)\to\widetilde C\,
(\subset T^{\perp}_{x_0}M=T^{\perp}_{u_0}\widetilde M)$ by 
$\exp^{\perp}(u(t))=x(t)$ ($t\in[0,T)$).  Set 
$\widetilde M_t:=(\pi\circ\phi)^{-1}(M_t)$ ($t\in[0,T)$).  It is shown that 
$\widetilde M_t$ ($t\in[0,T)$) is the mean curvature flow having 
$\widetilde M$ as initial data because the mean curvature vector of 
$\widetilde M_t$ is the horizontal lift of that of $M_t$ through 
$\pi\circ\phi$.  By investigating $u:[0,T)\to T^{\perp}_{u_0}\widetilde M$, we 
obtain the following fact corresponding to the second-half part of (i) and 
(ii) of Fact 1.  

\vspace{0.5truecm}

\noindent
{\bf Theorem A.} {\sl Let $M$ be an equifocal submanifold in a symmetric space 
$G/K$ of compact type and $C$ be the image of the fundamental domain of 
the Coxeter group of $M$ at $x_0\,(\in M)$ by the normal exponential map.  
Then the following statements ${\rm (i)}$ and ${\rm (ii)}$ hold:

${\rm (i)}$ For any focal submanifold $F$ of $M$, there exists a parallel 
submanifold $M'$ of $M$ such that the mean curvature flow having $M'$ 
as initial data converges to $F$ in finite time.  

${\rm (ii)}$ Assume that $M$ is irreducible and the codimension of $M$ is 
greater than one.  
If the mean curvature flow $M'_t$ having  a parallel submanifold 
$M'$ of $M$ as initial data converges to a focal submanifold $F$ of $M$ 
and the fibration of $M'$ onto $F$ is spherical fibration, then the mean 
curvature flow $M'_t$ has type I singularity.}

\vspace{0.5truecm}

Also, we obtain the following fact corresponding to the second-half part of 
(i) and (ii) of Fact 2 
for the mean curvature flow having a focal submanifold of an equifocal 
submanifold as initial data.  

\vspace{0.5truecm}

\noindent
{\bf Theorem B.} {\sl Let $M$ and $C$ be as in the statement of Theorem A and 
$\sigma$ be a stratum of dimension greater than zero of $\partial C$ 
(which is a stratified space).  
Then the following statements ${\rm (i)}$ and ${\rm (ii)}$ hold:

${\rm (i)}$ For any focal submanifold $F$ of $M$ through $\partial\sigma$, 
there exists a focal submanifold $F'$ of $M$ through $\sigma$ such that 
the mean curvature flow having $F'$ as initial data converges to $F$ in finite 
time.  

${\rm (ii)}$ Assume that $M$ is irreducible and the codimension of $M$ is 
greater than one.  
If the mean curvature flow $F'_t$ having a focal submanifold 
$F'$ of $M$ through $\sigma$ as initial data converges to a focal submanifold 
$F$ of $M$ through $\partial\sigma$ and the fibration of $F'$ onto $F$ 
is spherical fibration, then the mean curvature flow $F'_t$ has type I 
singularity.}

\vspace{0.5truecm}

For an equifocal submanifold in a symmetric space of compact type, 
the first-half part of the statement (i) of Facts 1 and 2 do not hold.  
In fact, there exist 
Hermann actions admitting a minimal principal orbit (which is equifocal) 
and a minimal singular orbit (which is a focal submanifold of 
a principal orbit).  So the following question arises naturally.  

\vspace{0.5truecm}

\noindent
{\bf Question.} {\sl Let $M$ and $\sigma$ be as in Theorem B and $F$ be 
a focal submanifold of $M$ through $\sigma$.  

{\rm(i)} Does the mean curvature flow having $M$ as initial data converge to 
a focal submanifold of $M$ in finite time in the case where $M$ 
is not minimal ? 

{\rm(ii)} Does the mean curvature flow having $F$ as initial data converge 
to a focal submanifold of $F$ through $\partial\sigma$ in finite time in 
the case where $F$ is not minimal ?}

\vspace{0.5truecm}

According to the homogeneity theorem for an equifocal submanifold by 
Christ [Ch], all irreducible equifocal submanifolds of codimension greater 
than one in symmetric spaces of compact type are homogeneous.  
Hence, according to the result by Heintze-Palais-Terng-Thorbergsson [HPTT], 
they are principal orbits of hyperpolar actions.  
Furthermore, according to the classification by Kollross [Kol] of hyperpolar 
actions on irreducible symmetric spaces of compact type, all hyperpolar 
actions of cohomogeneity greater than one on the symmetric spaces are Hermann 
actions.  
Therefore, all equifocal submanifolds of codimension greater than one in 
irreducible symmetric spaces of compact type are principal orbits of 
Hermann actions.  In the last section, we describe explicitly the mean 
curvature flows having orbits of some Hermann actions as initial data.  

\section{Preliminaries}
In this section, we briefly review the quantities associated with an 
isoparametric submanifold in an (infinite dimensional separable) Hilbert 
space, which was introduced by Terng [T2].  
Let $M$ be an isoparametric submanifold in a Hilbert space $V$.  

\vspace{0.5truecm}

\noindent
{\bf 2.1. Principal curvatures, curvature normals and curvature distributions}
Let $E_0$ and $E_i$ ($i\in I$) be all the curvature distributions of $M$, 
where $E_0$ is defined by 
$(E_0)_x=\displaystyle{\mathop{\cap}_{v\in T^{\perp}_xM}
{\rm Ker}\,A_v}\,(x\in M)$.  For each $x\in M$, we have 
$T_xM=\overline{(E_0)_x\oplus
\displaystyle{\left(\mathop{\oplus}_{i\in I}(E_i)_x\right)}}$, 
which is the common eigenspace decomposition of $A_v$'s 
($v\in T^{\perp}_xM$).  
Also, let $\lambda_i$ ($i\in I$) be the principal curvatures of $M$, that is, 
$\lambda_i$ is the section of the dual bundle $(T^{\perp}M)^{\ast}$ of 
$T^{\perp}M$ such that $A_v\vert_{(E_i)_x}=(\lambda_i)_x(v){\rm id}$ 
holds for any $x\in M$ and any $v\in T^{\perp}_xM$, and ${\bf n}_i$ be 
the curvature normal corresponding to $\lambda_i$, that is, 
$\lambda_i(\cdot)=\langle{\bf n}_i,\cdot\rangle$.  

\vspace{0.5truecm}

\noindent
{\bf 2.2. The Coxeter group associated with an isoparametric submanifold}
Denote by ${\it l}^x_i$ the affine hyperplane $(\lambda_i)_x^{-1}(1)$ in 
$T^{\perp}_xM$.  The focal set of $M$ at $x$ is equal to the sum 
$\displaystyle{\mathop{\cup}_{i\in I}(x+{\it l}^x_i)}$ of the affine 
hyperplanes $x+{\it l}_i^x$'s ($i\in I$) in the affine subspace 
$x+T^{\perp}_xM$ of $V$.  Each affine hyperplane ${\it l}_i^x$ is called a 
focal hyperplane of $M$ at $x$.  Let $W$ be the group generated by the 
reflection $R_i^x$'s ($i\in I$) with respect to ${\it l}_i^x$.  
This group is independent of the choice $x$ of $M$ up to group isomorphism.  
This group is called the Coxeter group associated with $M$.  
The fundamental domain of the Coxeter group containing the zero vector of 
$T^{\perp}_xM$ is given by 
$\{v\in T^{\perp}_xM\,\vert\,\lambda_i(v)<1\,(i\in I)\}$.  

\vspace{0.5truecm}

\noindent
{\bf 2.3. Principal curvatures of parallel submanifolds}
Let $M_v$ be the parallel submanifold of $M$ for a (non-focal) parallel normal 
vector field $v$, that is $M_v=\eta_v(M)$, where $\eta_v$ is the end-point map 
for $v$.  Denote by $A^v$ the shape tensor of $M_v$.  
This submanifold $M_v$ also is isoparametric and 
$A^v_w\vert_{\eta_{v\ast}(E_i)_x}
=\frac{(\lambda_i)_x(w)}{1-(\lambda_i)_x(v_x)}{\rm id}
\,\,(i\in I)$ for any $w\in T^{\perp}_{\eta_v(x)}M_v$, that is, 
$\frac{\lambda_i}{1-\lambda_i(v)}$'s ($i\in I$) are the principal curvatures 
of $M_v$ and hence $\frac{{\bf n}_i}{1-\lambda_i(v)}$'s ($i\in I$) are 
the curvature normals of $M_v$, where we identify $T^{\perp}_{\eta_v(x)}M_v$ 
with $T^{\perp}_xM$.  

\vspace{0.5truecm}

\noindent
{\bf 2.4. The mean curvature vector of a regularizable submanifold}
Assume that $M$ is regularizable in sense of [HLO], that is, 
for each normal vector $v$ of $M$,  the regularizable trace ${\rm Tr}_r\,A_v$ 
and ${\rm Tr}\,A_v^2$ exist, where ${\rm Tr}_r\,A_v$ is defined by 
${\rm Tr}_r\,A_v:=\sum\limits_{i=1}^{\infty}(\mu^+_i+\mu^-_i)$ as 
${\rm Spec}\,A_v\setminus\{0\}
=\{\mu^+_1,\,\mu^-_1,\,\mu^+_2,\,\mu^-_2,\,\cdots,\}$ 
$(\mu^-_1<\mu^-_2<\cdots<0<\cdots<\mu^+_2<\mu^+_1)$, where ${\rm Spec}\,A_v$ 
is the spectrum of $A_v$.  
Then the mean curvature vector $H$ of $M$ is defined by 
$\langle H,v\rangle={\rm Tr}_r\,A_v\,\,(\forall\,v\in T^{\perp}M)$.  

\vspace{0.5truecm}

Let $M$ be an equifocal submanifold in a symmetric space $G/K$ of compact type 
and set $\widetilde M:=(\pi\circ\phi)^{-1}(M)$, where $\pi$ is the natural 
projection of $G$ onto $G/K$ and $\phi:H^0([0,1],\mathfrak g)\to G$ is the 
parallel transport map for $G$.  

\vspace{0.5truecm}

{\bf 2.5. The mean curvature vector of the lifted submanifold} 
Denote by $\widetilde H$ (resp. $H$) the mean 
curvature vector of $\widetilde M$ (resp. $M$).  Then $\widetilde M$ is 
a regularizable isoparametric submanifold and $\widetilde H$ 
is equal to the horizontal lift of $nH^L$ of $nH$ ($n:={\rm dim}\,M$) 
(see Lemma 5.2 of [HLO]).  

\section{Proofs of Theorems A and B}
In this section, we prove Theorem A.  Let $M$ be an equifocal submanifold in 
a symmetric space $G/K$ of compact type, 
$\pi:G\to G/K$ be the natural 
projection and $\phi$ be the parallel transport map for $G$.  Set 
$\widetilde M:=(\pi\circ\phi)^{-1}(M)$.  Let $\widetilde C(\subset 
T^{\perp}_{x_0}M)$ be the fundamental domain 
of the Coxeter group of $M$ at $x_0(\in M)$ containing 
the zero vector ${\bf 0}$ of $T^{\perp}_{x_0}M$ and set 
$C:=\exp^{\perp}(\widetilde C)$, where $\exp^{\perp}$ is 
the normal exponential map of $M$.  Fix $u_0\in (\pi\circ\phi)^{-1}(x_0)$.  
Under the identification of $T^{\perp}_{x_0}M$ and 
$T^{\perp}_{u_0}\widetilde M$, the Coxeter group of $M$ is regarded as that 
of the isoparametric submanifold $\widetilde M$.  Denote by $H$ (resp. 
$\widetilde H$) the mean curvature vector of $M$ (resp. $\widetilde M$).  
The mean curvature vector $H$ and $\widetilde H$ are a parallel normal 
vector field of $M$ and $\widetilde M$, respectively.  
Let $v$ be a parallel normal vector field of $M$ and $v^L$ be 
the horizontal lift of $v$ to $H^0([0,1],\mathfrak g)$, which is a parallel 
normal vector field of $\widetilde M$.  Denote by $M_v$ (resp. 
$\widetilde M_{v^L}$) the parallel submanifold $\eta_v(M)$ (resp. 
$\eta_{v^L}(\widetilde M)$) of $M$ (resp. $\widetilde M$), where $\eta_v$ 
(resp. $\eta_{v^L}$) is the end-point map for $v$ (resp. $v^L$).  
Then we have $\widetilde M_{v^L}=(\pi\circ\phi)^{-1}(M_v)$.  
Denote by $H^v$ (resp. $\widetilde H^{v^L}$) the mean curvature vector of 
$M_v$ (resp. $\widetilde M_{v^L}$).  Define a vector field $X$ on 
$\widetilde C\,(\subset T^{\perp}_{u_0}\widetilde M=T^{\perp}_{x_0}M)$ by 
$X_w:=(\widetilde H^{\widetilde w})_{u_0+w}$ ($w\in \widetilde C$), 
where $\widetilde w$ is the parallel normal 
vector field of $\widetilde M$ with $\widetilde w_{u_0}=w$.  Let 
$\xi\,:\,(-S,T)\to\widetilde C$ be the maximal integral curve of $X$ with 
$\xi(0)={\bf 0}$.  Note that $S$ and $T$ are possible be equal to $\infty$.  
Let $\widetilde{\xi(t)}$ be the parallel normal vector field of $M$ with 
$\widetilde{\xi(t)}_{x_0}=\xi(t)$.  Let $M_t$ (resp. $\widetilde M_t$) be the 
mean curvature flow having $M$ (resp. $\widetilde M$) as initial data.  

\vspace{0.5truecm}

\noindent
{\bf Lemma 3.1.} {\sl For all $t\in[0,T)$, we have 
$M_t=M_{\widetilde{\xi(t)}}$ and $\widetilde M_t=
\widetilde M_{\widetilde{\xi(t)}^L}$.}

\vspace{0.5truecm}

\noindent
{\it Proof.} Fix $t_0\in[0,T)$.  Define a flow $f:\widetilde M\times[0,T)\to 
H^0([0,1],\mathfrak g)$ by $f(u,t):=\eta_{\widetilde{\xi(t)}^L}(u)$ 
($(u,t)\in\widetilde M\times[0,T)$), where we note that 
$f_t(\widetilde M)=\widetilde M_{\widetilde{\xi(t)}^L}$.  
For simplicity, denote by ${\widetilde H}^{t_0}$ the mean 
curvature vector of ${\widetilde M}_{\widetilde{\xi(t_0)}^L}$.  It is easy to 
show that $f_{\ast}((\frac{\partial}{\partial t})_{(\cdot,t_0)})$ is a 
parallel normal vector field of ${\widetilde M}_{\widetilde{\xi(t)}^L}$ and 
that $f_{\ast}((\frac{\partial}{\partial t})_{(u_0,t_0)})
=({\widetilde H}^{t_0})_{f_{t_0}(u_0)}$.  On the other hand, since 
${\widetilde M}_{\widetilde{\xi(t_0)}^L}$ is isoparametric, 
${\widetilde H}^{t_0}$ 
is also a parallel normal vector field of 
${\widetilde M}_{\widetilde{\xi(t_0)}^L}$.  Hence we have 
$f_{\ast}((\frac{\partial}{\partial t})_{(\cdot,t_0)})=\widetilde H^{t_0}$.  
Therefore, it follows from the arbitrariness of $t_0$ that $f$ is the mean 
curvature flow having $\widetilde M$ as initial data, that is, 
${\widetilde M}_{\widetilde{\xi(t)}^L}={\widetilde M}_t$ ($t\in[0,T)$) holds.  
Define a flow $\overline f:M\times [0,T)\to G/K$ by $\overline f(x,t):=
\eta_{\widetilde{\xi(t)}}(x)$ ($(x,t)\in M\times[0,T)$), where we note that 
$\overline f_t(M)=M_{\widetilde{\xi(t)}}$ ($\overline f_t(\cdot):=\overline f
(\cdot,t)$).  For simplicity, denote by $H^t$ the mean curvature vector of 
$M_{\widetilde{\xi(t)}}$.  Fix $t_0\in[0,T)$.  Since 
$\widetilde M_{\widetilde{\xi(t_0)}^L}=(\pi\circ\phi)^{-1}
(M_{\widetilde{\xi(t_0)}})$, we have 
$(H^{t_0})^L=\widetilde H^{t_0}$.  On the other hand, we have 
$\overline f_{\ast}((\frac{\partial}{\partial t})_{(\cdot,t_0)})^L
=f_{\ast}((\frac{\partial}{\partial t})_{(\cdot,t_0)})(=\widetilde H^{t_0})$.  
Hence we have $\overline f_{\ast}((\frac{\partial}{\partial t})_{(\cdot,t_0)})
=H^{t_0}$.  Therefore, it follows from the arbitrariness of $t_0$ that 
$\overline f$ is the mean curvature flow having $M$ as initial data, that is, 
$M_{\widetilde{\xi(t)}}=M_t$ ($t\in[0,T)$).  \hspace{0.1truecm}q.e.d.

\vspace{0.5truecm}

Clearly we suffice to show the statement of Theorem A in the case where $M$ is 
full.  Hence, in the sequel, we assume that $M$ is full.  
Denote by $\Lambda$ the set of all principal curvatures of $\widetilde M$.  
Set $r:={\rm codim}\,M$.  It is shown that the set of all focal hyperplanes of 
$\widetilde M$ is given as the sum of finite pieces of infinite parallel 
families consisting of hyperplanes in $T^{\perp}_{u_0}\widetilde M$ which 
arrange at equal intervals.  Let 
$\{{\it l}^{u_0}_{aj}\,\vert\,j\in{\Bbb Z}\}$ ($1\leq a\leq \bar r$) be 
the finite pieces of infinite parallel families consisting of hyperplanes 
in $T^{\perp}_{u_0}\widetilde M$.  
Since ${\it l}^{u_0}_{aj}$'s ($j\in{\Bbb Z}$) arrange at equal intervals, 
we can express as 
$$\Lambda=\mathop{\cup}_{a=1}^{\bar r}\{\frac{\lambda_a}{1+b_aj}\,\vert\,
j\in{\bf Z}\},$$
where $\lambda_a$'s and $b_a$'s are parallel sections of 
$(T^{\perp}\widetilde M)^{\ast}$ and positive constants greater than one, 
respectively, 
which are defined by $((\lambda_a)_{u_0})^{-1}(1+b_aj)={\it l}^{u_0}_{aj}$.  
For simplicity, we set $\lambda_{aj}:=\frac{\lambda_a}{1+b_aj}$.  
Denote by ${\bf n}_{aj}$ and $E_{aj}$ the curvature normal and the curvature 
distribution corresponding to $\lambda_{aj}$, respectively.  
It is shown that, for each $a$, $\lambda_{a,2j}$'s ($j\in{\Bbb Z}$) have the 
same multiplicity and so are also $\lambda_{a,2j+1}$'s ($j\in{\Bbb Z}$).  
Denote by $m_a^e$ and $m_a^o$ the multiplicities of $\lambda_{a,2j}$ and 
$\lambda_{a,2j+1}$, respectively.  Take a parallel normal vector field 
$v$ of $\widetilde M$ with $v_{u_0}\in\widetilde C$.  Denote by 
${\widetilde A}^v$ (resp. ${\widetilde H}^v$) the shape tensor (resp. 
the mean curvature vector) of the parallel submanifold ${\widetilde M}_v$.  
Since ${\widetilde A}^v_w\vert_{\eta_{v\ast}(E_{aj})_u}
=\frac{(\lambda_{aj})_u(w)}{1-(\lambda_{aj})_u(v_u)}\,{\rm id}$ 
$(w\in T^{\perp}_u\widetilde M$), we have 
$$\begin{array}{l}
\displaystyle{{\rm Tr}_r\widetilde A^v_w=\sum_{a=1}^{\bar r}\left(
\sum_{j\in{\bf Z}}\frac{m_a^e\lambda_{a,2j}(w)}{1-\lambda_{a,2j}(v_u)}+
\sum_{j\in{\bf Z}}\frac{m_a^o\lambda_{a,2j+1}(w)}{1-\lambda_{a,2j+1}(v_u)}
\right)}\\
\hspace{1.1truecm}\displaystyle{=\sum_{a=1}^{\bar r}\left(
(m_a^e+m_a^o)\cot(\frac{\pi}{b_a}(1-\lambda_a(v_u)))\right.}\\
\hspace{2.4truecm}
\displaystyle{\left.
+(m_a^e-m_a^o){\rm cosec}(\frac{\pi}{b_a}(1-\lambda_a(v_u)))
\right)\frac{\pi}{2b_a}\lambda_a(w),}
\end{array}$$
where we use the relation 
$\frac{\cos\theta+1}{\sin\theta}=\sum\limits_{j\in{\bf Z}}
\frac{2}{\theta+2j\pi}$.  Therefore we have 
$$\begin{array}{l}
\displaystyle{\widetilde H^v=\sum_{a=1}^{\bar r}\left(
(m_a^e+m_a^o)\cot(\frac{\pi}{b_a}(1-\lambda_a(v)))\right.}\\
\hspace{1.2truecm}\displaystyle{\left.+(m_a^e-m_a^o){\rm cosec}
(\frac{\pi}{b_a}(1-\lambda_a(v)))\right)\frac{\pi}{2b_a}{\bf n}_a,}
\end{array}\leqno{(3.1)}$$
where ${\bf n}_a:={\bf n}_{a0}$.  

\vspace{0.5truecm}

Now we prove Theorem A.  

\vspace{0.5truecm}

\noindent
{\it Proof of Theorem A.} In this proof, we use the above notations.  
Let $X$ be the above vector field on $\widetilde C$.  From $(3.1)$, we have 
$$\begin{array}{l}
\displaystyle{X_w=\sum_{a=1}^{\bar r}\left(
(m_a^e+m_a^o)\cot(\frac{\pi}{b_a}(1-(\lambda_a)_{u_0}(w)))\right.}\\
\hspace{1.6truecm}\displaystyle{\left.+(m_a^e-m_a^o)
{\rm cosec}(\frac{\pi}{b_a}(1-(\lambda_a)_{u_0}(w)))\right)
\frac{\pi}{2b_a}({\bf n}_a)_{u_0}}\\
\hspace{0.6truecm}\displaystyle{=
\sum_{a=1}^{\bar r}\left(
\frac{\cos(\frac{\pi}{b_a}(1-(\lambda_a)_{u_0}(w)))+1}
{\sin(\frac{\pi}{b_a}(1-(\lambda_a)_{u_0}(w)))}\times m_a^e\right.}\\
\hspace{1.6truecm}\displaystyle{
\left.+\frac{\cos(\frac{\pi}{b_a}(1-(\lambda_a)_{u_0}(w)))-1}
{\sin(\frac{\pi}{b_a}(1-(\lambda_a)_{u_0}(w)))}\times m_a^o\right)
\frac{\pi}{2b_a}({\bf n}_a)_{u_0}.}
\end{array}\leqno{(3.2)}$$
Denote by $\bar r'$ the cardinal number of the set 
$\{a\,\vert\,(\lambda_a)_{u_0}^{-1}(1)\cap\partial\widetilde C\,:\,{\rm open}
\,\,{\rm in}\,\,\partial\widetilde C\}$.  Clearly we have 
$r+1\leq\bar r'\leq\bar r$.  By reordering $\{1,\cdots,\bar r\}$, we may 
assume that this set is equal to $\{1,\cdots,\bar r'\}$.  
Fix $w_0\in\widetilde C$ and $a_0\in\{1,\cdots,\bar r'\}$.  
For simplicity, set $\widetilde{\sigma}_{a_0}
:=(\lambda_{a_0})_{u_0}^{-1}(1)\cap\partial\widetilde C$.  Let $w'_0$ be 
the point of $\widetilde{\sigma}_{a_0}$ such that $w_0-w'_0$ is normal to 
$\widetilde{\sigma}_{a_0}$ and set $w_0^{\varepsilon}:=\varepsilon w_0
+(1-\varepsilon)w'_0$ for $\varepsilon\in(0,1)$.  
Then we have $\lim\limits_{\varepsilon\to+0}(\lambda_{a_0})_{u_0}
(w^{\varepsilon}_0)=1$ and $\displaystyle{\mathop{\sup}_{0<\varepsilon<1}
(\lambda_a)_{u_0}(w^{\varepsilon}_0)<1}$ for each 
$a\in\{1,\cdots,\bar r\}\setminus\{a_0\}$.  
Hence we have 
$$\lim_{\varepsilon\to+0}\frac{\cos(\frac{\pi}{b_{a_0}}
(1-(\lambda_{a_0})_{u_0}(w_0^{\varepsilon})))+1}
{\sin(\frac{\pi}{b_{a_0}}(1-(\lambda_{a_0})_{u_0}(w_0^{\varepsilon})))}=\infty,
\,\,
\lim_{\varepsilon\to+0}\frac{\cos(\frac{\pi}{b_{a_0}}
(1-(\lambda_{a_0})_{u_0}(w_0^{\varepsilon})))-1}
{\sin(\frac{\pi}{b_{a_0}}(1-(\lambda_{a_0})_{u_0}(w_0^{\varepsilon})))}=0$$
and
$$\mathop{\sup}_{0<\varepsilon<1}\frac{\cos(\frac{\pi}{b_a}
(1-(\lambda_a)_{u_0}(w_0^{\varepsilon})))\pm1}
{\sin(\frac{\pi}{b_a}(1-(\lambda_a)_{u_0}(w_0^{\varepsilon})))}\,<\,\infty.$$
Therefore, we have $\lim\limits_{\varepsilon\to+0}
\frac{X_{w_0^{\varepsilon}}}{\vert\vert X_{w_0^{\varepsilon}}\vert\vert}
=\frac{({\bf n}_{a_0})_{u_0}}{\vert\vert({\bf n}_{a_0})_{u_0}\vert\vert}$ and 
$\lim\limits_{\varepsilon\to+0}\vert\vert X_{w_0^{\varepsilon}}\vert\vert
=\infty$.  This implies that $X$ is as in Fig. 1 on a sufficiently small 
collar neighborhood of $\widetilde{\sigma}_{a_0}$.  
Therefore, the flow of $X$ is as in Fig. 2 on a sufficiently small collar 
neighborhood of $\partial\widetilde C$.  
Take an arbitrary focal submanifold $F$ of $M$ and set $\widetilde F:=
(\pi\circ\phi)^{-1}(F)$, which is a focal 

\vspace{0.5truecm}

\centerline{
\unitlength 0.1in
\begin{picture}( 32.1400, 16.4100)( 16.5000,-24.0100)
%
\special{pn 8}%
\special{pa 3600 1200}%
\special{pa 2800 2400}%
\special{pa 4400 2400}%
\special{pa 4400 2400}%
\special{pa 3600 1200}%
\special{fp}%
%
\special{pn 8}%
\special{pa 3790 1730}%
\special{pa 4314 1408}%
\special{fp}%
\special{sh 1}%
\special{pa 4314 1408}%
\special{pa 4248 1426}%
\special{pa 4270 1436}%
\special{pa 4268 1460}%
\special{pa 4314 1408}%
\special{fp}%
%
\special{pn 8}%
\special{pa 3730 1600}%
\special{pa 4256 1278}%
\special{fp}%
\special{sh 1}%
\special{pa 4256 1278}%
\special{pa 4188 1296}%
\special{pa 4210 1306}%
\special{pa 4210 1330}%
\special{pa 4256 1278}%
\special{fp}%
%
\special{pn 8}%
\special{pa 3680 1460}%
\special{pa 4206 1138}%
\special{fp}%
\special{sh 1}%
\special{pa 4206 1138}%
\special{pa 4138 1156}%
\special{pa 4160 1166}%
\special{pa 4160 1190}%
\special{pa 4206 1138}%
\special{fp}%
%
\special{pn 8}%
\special{pa 3640 1340}%
\special{pa 4164 1016}%
\special{fp}%
\special{sh 1}%
\special{pa 4164 1016}%
\special{pa 4098 1034}%
\special{pa 4120 1044}%
\special{pa 4118 1068}%
\special{pa 4164 1016}%
\special{fp}%
%
\special{pn 8}%
\special{pa 4230 2260}%
\special{pa 4754 1938}%
\special{fp}%
\special{sh 1}%
\special{pa 4754 1938}%
\special{pa 4688 1956}%
\special{pa 4710 1966}%
\special{pa 4708 1990}%
\special{pa 4754 1938}%
\special{fp}%
%
\special{pn 8}%
\special{pa 4020 2070}%
\special{pa 4546 1748}%
\special{fp}%
\special{sh 1}%
\special{pa 4546 1748}%
\special{pa 4478 1766}%
\special{pa 4500 1776}%
\special{pa 4500 1800}%
\special{pa 4546 1748}%
\special{fp}%
%
\special{pn 8}%
\special{pa 3930 1950}%
\special{pa 4456 1628}%
\special{fp}%
\special{sh 1}%
\special{pa 4456 1628}%
\special{pa 4388 1646}%
\special{pa 4410 1656}%
\special{pa 4410 1680}%
\special{pa 4456 1628}%
\special{fp}%
%
\special{pn 8}%
\special{pa 3860 1840}%
\special{pa 4384 1516}%
\special{fp}%
\special{sh 1}%
\special{pa 4384 1516}%
\special{pa 4318 1534}%
\special{pa 4340 1544}%
\special{pa 4338 1568}%
\special{pa 4384 1516}%
\special{fp}%
\put(36.9000,-18.3000){\makebox(0,0)[rt]{$\widetilde C$}}%
%
\special{pn 8}%
\special{ar 5500 930 1920 1800  2.1850120 2.1914636}%
\special{ar 5500 930 1920 1800  2.2108185 2.2172701}%
\special{ar 5500 930 1920 1800  2.2366249 2.2430765}%
\special{ar 5500 930 1920 1800  2.2624314 2.2688830}%
\special{ar 5500 930 1920 1800  2.2882378 2.2946894}%
\special{ar 5500 930 1920 1800  2.3140443 2.3204959}%
\special{ar 5500 930 1920 1800  2.3398507 2.3463023}%
\special{ar 5500 930 1920 1800  2.3656572 2.3721088}%
\special{ar 5500 930 1920 1800  2.3914636 2.3979152}%
\special{ar 5500 930 1920 1800  2.4172701 2.4237217}%
\special{ar 5500 930 1920 1800  2.4430765 2.4495281}%
\special{ar 5500 930 1920 1800  2.4688830 2.4753346}%
\special{ar 5500 930 1920 1800  2.4946894 2.5011410}%
\special{ar 5500 930 1920 1800  2.5204959 2.5269475}%
\special{ar 5500 930 1920 1800  2.5463023 2.5527539}%
\special{ar 5500 930 1920 1800  2.5721088 2.5785604}%
\special{ar 5500 930 1920 1800  2.5979152 2.6043669}%
\special{ar 5500 930 1920 1800  2.6237217 2.6301733}%
\special{ar 5500 930 1920 1800  2.6495281 2.6559798}%
\special{ar 5500 930 1920 1800  2.6753346 2.6817862}%
\special{ar 5500 930 1920 1800  2.7011410 2.7075927}%
\special{ar 5500 930 1920 1800  2.7269475 2.7333991}%
\special{ar 5500 930 1920 1800  2.7527539 2.7592056}%
\special{ar 5500 930 1920 1800  2.7785604 2.7850120}%
\special{ar 5500 930 1920 1800  2.8043669 2.8108185}%
\special{ar 5500 930 1920 1800  2.8301733 2.8366249}%
\special{ar 5500 930 1920 1800  2.8559798 2.8624314}%
\special{ar 5500 930 1920 1800  2.8817862 2.8882378}%
\special{ar 5500 930 1920 1800  2.9075927 2.9140443}%
\special{ar 5500 930 1920 1800  2.9333991 2.9398507}%
\special{ar 5500 930 1920 1800  2.9592056 2.9656572}%
\special{ar 5500 930 1920 1800  2.9850120 2.9876480}%
%
\special{pn 8}%
\special{pa 4120 2170}%
\special{pa 4644 1848}%
\special{fp}%
\special{sh 1}%
\special{pa 4644 1848}%
\special{pa 4578 1866}%
\special{pa 4600 1876}%
\special{pa 4598 1900}%
\special{pa 4644 1848}%
\special{fp}%
%
\special{pn 8}%
\special{pa 4340 2360}%
\special{pa 4864 2038}%
\special{fp}%
\special{sh 1}%
\special{pa 4864 2038}%
\special{pa 4798 2056}%
\special{pa 4820 2066}%
\special{pa 4818 2090}%
\special{pa 4864 2038}%
\special{fp}%
%
\special{pn 8}%
\special{pa 3860 940}%
\special{pa 3720 1370}%
\special{dt 0.045}%
\special{sh 1}%
\special{pa 3720 1370}%
\special{pa 3760 1314}%
\special{pa 3738 1320}%
\special{pa 3722 1300}%
\special{pa 3720 1370}%
\special{fp}%
\put(39.9000,-7.6000){\makebox(0,0)[rt]{$\widetilde{\sigma}_{a_0}$}}%
\end{picture}%
\hspace{3truecm}}

\vspace{1truecm}

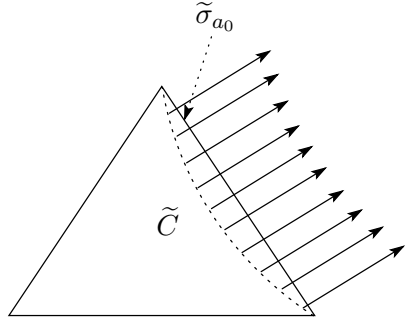
\centerline{{\bf Fig. 1.}}

\vspace{0.5truecm}

\centerline{
\unitlength 0.1in
\begin{picture}( 17.5000, 14.5000)( 24.0000,-25.8000)
%
\special{pn 8}%
\special{pa 2400 2574}%
\special{pa 4148 2574}%
\special{fp}%
%
\special{pn 8}%
\special{pa 2634 2434}%
\special{pa 3916 2434}%
\special{dt 0.045}%
%
\special{pn 8}%
\special{pa 2626 2434}%
\special{pa 3272 1376}%
\special{dt 0.045}%
%
\special{pn 8}%
\special{pa 2400 2574}%
\special{pa 2626 2434}%
\special{fp}%
%
\special{pn 8}%
\special{pa 4148 2574}%
\special{pa 3922 2434}%
\special{fp}%
%
\special{pn 8}%
\special{pa 3922 2434}%
\special{pa 3274 1368}%
\special{dt 0.045}%
%
\special{pn 8}%
\special{pa 2400 2574}%
\special{pa 3284 1130}%
\special{fp}%
%
\special{pn 8}%
\special{pa 3266 1138}%
\special{pa 4150 2580}%
\special{fp}%
%
\special{pn 8}%
\special{pa 3274 1382}%
\special{pa 3274 1152}%
\special{fp}%
%
\special{pn 8}%
\special{pa 3274 2574}%
\special{pa 3274 2434}%
\special{fp}%
%
\special{pn 8}%
\special{pa 3508 2574}%
\special{pa 3508 2434}%
\special{fp}%
%
\special{pn 8}%
\special{pa 3062 2574}%
\special{pa 3062 2434}%
\special{fp}%
%
\special{pn 8}%
\special{pa 2946 1904}%
\special{pa 2852 1848}%
\special{fp}%
%
\special{pn 8}%
\special{pa 2830 2100}%
\special{pa 2728 2036}%
\special{fp}%
%
\special{pn 8}%
\special{pa 3062 1730}%
\special{pa 2954 1668}%
\special{fp}%
%
\special{pn 8}%
\special{pa 3602 1898}%
\special{pa 3704 1848}%
\special{fp}%
%
\special{pn 8}%
\special{pa 3718 2100}%
\special{pa 3828 2058}%
\special{fp}%
%
\special{pn 8}%
\special{pa 3500 1730}%
\special{pa 3602 1682}%
\special{fp}%
%
\special{pn 8}%
\special{ar 3070 2574 240 230  3.1415927 3.7963951}%
%
\special{pn 8}%
\special{ar 3508 2574 220 252  5.7019104 6.2831853}%
%
\special{pn 8}%
\special{ar 3434 2574 518 216  5.5667723 6.2831853}%
%
\special{pn 8}%
\special{ar 3348 2574 722 288  3.1415927 3.6459400}%
%
\special{pn 8}%
\special{ar 3092 2184 218 866  4.9665570 5.3454873}%
%
\special{pn 8}%
\special{ar 3420 2044 162 710  4.1233793 4.5721859}%
%
\special{pn 8}%
\special{ar 3552 1808 212 328  3.9095599 4.3862332}%
%
\special{pn 8}%
\special{ar 2998 1884 220 412  5.0044263 5.4867582}%
%
\special{pn 8}%
\special{ar 3952 2008 240 378  1.2728089 1.8914058}%
%
\special{pn 8}%
\special{ar 3858 2016 212 230  1.1987620 1.8421182}%
%
\special{pn 8}%
\special{ar 2692 1974 176 294  1.2576866 2.0634369}%
%
\special{pn 8}%
\special{ar 2590 2148 220 224  1.1144746 1.8530537}%
\end{picture}%
\hspace{1truecm}}

\vspace{1truecm}

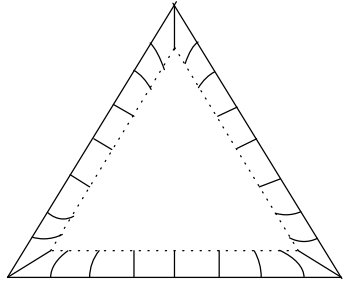
\centerline{{\bf Fig. 2.}}

\vspace{0.5truecm}

\noindent
submanifold of $\widetilde M$ and 
passes through $\partial\widetilde C$.  From the above fact for $X$, we see 
that there exists a parallel submanifold ${\widetilde M}'$ of $\widetilde M$ 
such that the mean curvature flow ${\widetilde M}'_t$ having ${\widetilde M}'$ 
as initial data converges to $\widetilde F$ in finite time.  Then 
$M':=(\pi\circ\phi)({\widetilde M}')$ is 
a parallel submanifold of $M$.  According to Lemma 3.1, $(\pi\circ\phi)
({\widetilde M}'_t)$ is the mean curvature flow having $M'$ as initial data 
and it converges to $F$ in finite time.  
Thus the statement (i) is shown.  Next we shall show the statement (ii).  
Assume that $M$ is irreducible and the codimension of $M$ is greater than 
one.  Let $M'$ and $F$ be as in the statement (ii).  Assume that the mean 
curvature flow having $M'$ as initial data converges to $F$ in finite time 
$T$.  Set 
${\widetilde M}':=(\pi\circ\phi)^{-1}(M')$ and $\widetilde F
:=(\pi\circ\phi)^{-1}(F)$.  Since the fibration of $M$ onto $F$ is a spherical 
fibration, $\widetilde F$ passes through a highest dimensional stratum 
$\widetilde{\sigma}$ of $\partial\widetilde C$.  Let $a_0$ be the element of 
$\{1,\cdots,\bar r'\}$ 
with $\widetilde{\sigma}\subset(\lambda_{a_0})_{u_0}^{-1}
(1)$.  Set ${\widetilde M}'_t:=(\pi\circ\phi)^{-1}(M'_t)$ ($t\in[0,T)$), which 
is the mean curvature flow having ${\widetilde M}'$ as initial data.  
Denote by $A^t$ (resp. $\widetilde A^t$) the shape tensor of $M'_t$ (resp. 
${\widetilde M}'_t$).  Let ${\widetilde M}'\cap(u_0+\widetilde C)=\{u_1\},\,
\xi:[0,T)\to\widetilde C$ be the integral curve 
of $X$ with $\xi(0)=u_1-u_0$ and $\widetilde{\xi(t)}$ be the parallel normal 
vector field of $M'$ with $\widetilde{\xi(t)}_{x_0}=\xi(t)$.  Then, since 
${\widetilde M}'_t$ is the parallel submanifold of ${\widetilde M}'$ for 
$\widetilde{\xi(t)}^L$, we have 
$${\rm Spec}\,\widetilde A^t_v\setminus\{0\}=\{\frac{(\lambda_{aj})_{u_0}(v)}
{1-(\lambda_{aj})_{u_0}(\xi(t))}\,\vert\,a=1,\cdots,\bar r',
\,\,j\in{\Bbb Z}\}$$
for each $v\in T^{\perp}_{u_0+\xi(t)}{\widetilde M}'_t=T^{\perp}_{u_0}
\widetilde M$.  
Since $\lim\limits_{t\to T-0}\xi(t)\in(\lambda_{a_0})_{u_0}^{-1}(1)$, we have 

\noindent
$\lim\limits_{t\to T-0}(\lambda_{a_0})_{u_0}(\xi(t))=1$ and 
$\lim\limits_{t\to T-0}(\lambda_a)_{u_0}(\xi(t))<1$ ($a\not=a_0$).  Hence we 
have 
$$\begin{array}{l}
\hspace{0.6truecm}
\displaystyle{\lim_{t\to T-0}\vert\vert\widetilde A^t_v\vert\vert_{\infty}^2
(T-t)}\\
\displaystyle{=\lim_{t\to T-0}\frac{(\lambda_{a_0})_{u_0}(v)^2}
{(1-(\lambda_{a_0})_{u_0}(\xi(t)))^2}(T-t)}\\
\displaystyle{=\frac12(\lambda_{a_0})_{u_0}(v)^2
\lim_{t\to T-0}\frac{1}{(1-(\lambda_{a_0})_{u_0}(\xi(t)))(\lambda_{a_0})_{u_0}
(\xi'(t))}.}
\end{array}\leqno{(3.3)}
$$
Since $\xi'(t)=(\widetilde H^{\widetilde{\xi(t)}^L})_{u_0+\xi(t)}$, we have 
$$\begin{array}{l}
\hspace{0.6truecm}\displaystyle{\lim_{t\to T-0}
(1-(\lambda_{a_0})_{u_0}(\xi(t)))(\lambda_{a_0})_{u_0}(\xi'(t))}\\
\displaystyle{=\lim_{t\to T-0}\sum_{a=1}^{\bar r}\left(
\frac{(m_a^e+m_a^o)(1-(\lambda_{a_0})_{u_0}(\xi(t)))}
{\tan(\frac{\pi}{b_a}(1-(\lambda_a)_{u_0}(\xi(t))))}\right.}\\
\hspace{2.5truecm}\displaystyle{
\left.+\frac{(m_a^e-m_a^o)(1-(\lambda_{a_0})_{u_0}(\xi(t)))}
{\sin(\frac{\pi}{b_a}(1-(\lambda_a)_{u_0}(\xi(t))))}\right)
\frac{\pi}{2b_a}\langle({\bf n}_a)_{u_0},({\bf n}_{a_0})_{u_0}\rangle}\\
\displaystyle{=\frac12\lim_{t\to T-0}\left((m_{a_0}^e+m_{a_0}^o)\cos^2
(\frac{\pi}{b_{a_0}}(1-(\lambda_{a_0})_{u_0}
(\xi(t))))\right.}\\
\hspace{2.5truecm}\displaystyle{
\left.+(m_{a_0}^e-m_{a_0}^o)\frac{1}{\cos(\frac{\pi}{b_a}
(1-(\lambda_{a_0})_{u_0}(\xi(t))))}\right)\vert\vert({\bf n}_{a_0})_{u_0}
\vert\vert^2}\\
\displaystyle{
=m_{a_0}^e\vert\vert({\bf n}_{a_0})_{u_0}\vert\vert^2,}
\end{array}$$
which together with $(3.3)$ deduces
$$\lim_{t\to T-0}\vert\vert\widetilde A^t_v\vert\vert_{\infty}^2(T-t)
=\frac{(\lambda_{a_0})_{u_0}(v)^2}{2m_{a_0}^e\vert\vert({\bf n}_{a_0})_{u_0}
\vert\vert^2}
\leqno{(3.4)}$$
and hence 
$$\lim_{t\to T-0}\mathop{\max}_{v\in S^{\perp}_{u_0+\xi(t)}{\widetilde M}'_t}
\vert\vert\widetilde A^t_v\vert\vert_{\infty}^2(T-t)
=\frac{1}{2m_{a_0}^e}.$$
Thus the mean curvature flow ${\widetilde M}'_t$ has type I singularity.  Set 
$\bar v_t:=(\pi\circ\phi)_{\ast u_0+\xi(t)}(v)$ and let 
$\{\lambda^t_1,\cdots,\lambda^t_n\}\,\,(\lambda^t_1\leq\cdots\leq\lambda^t_n)$ 
(resp. $\{\mu^t_1,\cdots,\mu^t_n\}\,\,(0\leq\mu^t_1\leq\cdots\leq\mu^t_n)$) 
be all the eigenvalues of $A^t_{\bar v_t}$ 
(resp. $R(\cdot,\bar v_t)\bar v_t$), where $n:={\rm dim}\,M$.  
Since $M$ is an irreducible equifocal 
submanifold of codimension greater than one by the assumption, it is 
homogeneous by the homogeneity theorem of Christ (see [Ch]) and hence it is 
a principal orbit of a Hermann action by the result of 
Heintze-Palais-Terng-Thorbergsson (see [HPTT]) and the classification of 
hyperpolar actions by Kollross (see [Kol]).  Furthermore, $M$ and its parallel 
submanifolds are curvature-adapted by the result of Goertsches-Thorbergsson 
(see [GT]).  
Therefore, $A^t_{\bar v_t}$ and $R(\cdot,\bar v_t)\bar v_t$ commute and hence 
we have 
$$\sum_{i=1}^n\sum_{j=1}^n\left(
{\rm Ker}(A^t_{\bar v_t}-\lambda^t_i\,{\rm id})\cap{\rm Ker}(R(\cdot,\bar v_t)
\bar v_t-\mu_j^t\,{\rm id})\right)=T_{(\pi\circ\phi)(u_0+\xi(t))}M'_t.$$
Set $\bar E_{ij}^t:={\rm Ker}(A^t_{\bar v_t}-\lambda_i^t\,{\rm id})\cap 
{\rm Ker}(R(\cdot,\bar v_t)\bar v_t-\mu_j^t\,{\rm id})$ 
($i,j\in\{1,\cdots,n\}$) and 
$I_t:=\{(i,j)\in\{1,\cdots,n\}^2\,\vert\,\bar E_{ij}^t\not=\{0\}\}$.  
For each $(i,j)\in I_t$, we have 
$${\rm Spec}(\widetilde A^t_v\vert_{(\pi\circ\phi)_{\ast}^{-1}(\bar E_{ij}^t)})
=
\left\{
\begin{array}{ll}
\displaystyle{\left\{\frac{\sqrt{\mu_j^t}}{\arctan\frac{\sqrt{\mu_j^t}}
{\lambda_i^t}+k\pi}\,\vert\,k\in{\Bbb Z}\right\}} & 
\displaystyle{(\mu_j^t\not=0)}\\
\displaystyle{\{\lambda_i^t\}} & \displaystyle{(\mu_j^t=0)}
\end{array}\right.$$
in terms of Proposition 3.2 of [Koi1] and hence 
$$\vert\vert\widetilde A^t_v\vert\vert_{\infty}=\max\left(\left\{
\frac{\sqrt{\mu_j^t}}{\arctan\frac{\sqrt{\mu_j^t}}{\vert\lambda_i^t\vert}}\,
\vert\,(i,j)\in I_t\,\,{\rm s.t.}\,\,\mu_j^t\not=0\right\}\cup
\{\vert\lambda_i^t\vert\,\vert\,(i,j)\in I_t\,\,{\rm s.t.}\,\,\mu_j^t=0\}
\right).$$
It is clear that 
$\displaystyle{\mathop{{\rm sup}}_{0\leq t<T}\mu_n^t\,<\,\infty}$.  
If $\lim\limits_{t\to T-0}\vert\lambda_i^t\vert=\infty$, then we have 
$\displaystyle{\lim_{t\to T-0}\left(\frac{\sqrt{\mu_j^t}}
{\arctan\frac{\sqrt{\mu_j^t}}{\lambda_i^t}}\right)/\lambda_i^t}$\newline
$=1$.  Hence we have 
$$\begin{array}{l}
\hspace{0.6truecm}\displaystyle{\lim\limits_{t\to T-0}
\vert\vert\widetilde A^t_v\vert\vert^2_{\infty}(T-t)}\\
\displaystyle{=\max\left\{\lim\limits_{t\to T-0}
(\lambda_i^t)^2(T-t)\,\vert\,i=1,\cdots,n\right\}}\\
\displaystyle{=\lim\limits_{t\to T-0}\max\{(\lambda_i^t)^2(T-t)\,\vert\,
i=1,\cdots,n\}}\\
\displaystyle{=\lim\limits_{t\to T-0}
\vert\vert A^t_{\bar v_t}\vert\vert^2_{\infty}(T-t),}
\end{array}$$
which together with $(3.4)$ deduces
$$\lim\limits_{t\to T-0}\vert\vert A^t_{\bar v_t}\vert\vert^2_{\infty}(T-t)=
\frac{(\lambda_{a_0})_{u_0}(v)^2}{2m_{a_0}^e\vert\vert
({\bf n}_{a_0})_{u_0}\vert\vert^2}.$$
Therefore we obtain 
$$\lim\limits_{t\to T-0}\mathop{\max}_{v\in S^{\perp}_{\exp^{\perp}(\xi(t))}
M'_t}\vert\vert A^t_v\vert\vert^2_{\infty}(T-t)=\frac{1}{2m_{a_0}^e}<\infty.$$
Thus the mean curvature flow $M'_t$ has type I singularity.  
\hspace{2.5truecm}q.e.d.

\vspace{0.5truecm}

Next we prove Theorem B.  

\vspace{0.5truecm}

\noindent
{\it Proof of Theorem B.} For simplicity, set $I:=\{1,\cdots,\bar r\}$, 
where $\bar r$ is as above.  
Let $\widetilde{\sigma}$ be a stratum of dimension greater than zero of 
$\partial\widetilde C$ and 
$I_{\widetilde{\sigma}}:=\{a\in I\,\vert\,\widetilde{\sigma}\subset 
(\lambda_a)_{u_0}^{-1}(1)\}$.  Let $w\in\widetilde{\sigma}$ and 
$\widetilde w$ be the parallel normal vector field of $M$ with 
$\widetilde w_{x_0}=w$.  The horizontal lift ${\widetilde w}^L$ of 
$\widetilde w$ to $H^0([0,1],\mathfrak g)$ is a parallel normal vector field 
of $\widetilde M$.  Denote by $F$ (resp. $\widetilde F$) the focal submanifold 
of $M$ (resp. $\widetilde M$) for $\widetilde w$ (resp. 
${\widetilde w}^L$) and by $\widetilde A^{\widetilde F}$ (resp. 
$\widetilde H^{\widetilde F}$) the shape tensor (resp. the mean curvature 
vector) of $\widetilde F$.  Then we have 
$$T_{u_0+w}\widetilde F
=\left(\mathop{\oplus}_{a\in I\setminus I_{\widetilde{\sigma}}}
\mathop{\oplus}_{j\in{\bf Z}}\eta_{\widetilde w\ast}(E_{aj})\right)
\oplus\left(\mathop{\oplus}_{a\in I_{\widetilde{\sigma}}}
\mathop{\oplus}_{j\in{\bf Z}\setminus\{0\}}\eta_{\widetilde w\ast}(E_{aj})
\right).$$
Also we have 
$$T^{\perp}_{u_0+w}\widetilde F
=\left(\mathop{\oplus}_{a\in I_{\widetilde{\sigma}}}
(E_{a0})_{u_0}\right)\oplus T^{\perp}_{u_0}\widetilde M,$$
where we identify $T_{u_0+w}H^0([0,1],\mathfrak g)$ with 
$T_{u_0}H^0([0,1],\mathfrak g)$.  
For $v\in T^{\perp}_{u_0}\widetilde M\,(\subset 
T^{\perp}_{u_0+w}\widetilde F)$, 
we have 
$$\widetilde A^{\widetilde F}_v\vert_{\eta_{\widetilde w\ast}((E_{aj})_{u_0})}
=\frac{(\lambda_{aj})_{u_0}(v)}{1-(\lambda_{aj})_{u_0}(w)}{\rm id}\,\,\,\,
((a,j)\in((I\setminus I_{\widetilde{\sigma}})\times{\Bbb Z})\cup
(I_{\widetilde{\sigma}}\times({\Bbb Z}\setminus\{0\}))).$$
Hence we have 
$$\begin{array}{l}
\hspace{0.6truecm}\displaystyle{{\rm Tr}_r\widetilde A^{\widetilde F}_v}\\
\displaystyle{=\sum_{a\in I\setminus I_{\widetilde{\sigma}}}\left(
\sum_{j\in{\bf Z}}\frac{m_a^e(\lambda_{a,2j})_{u_0}(v)}
{1-(\lambda_{a,2j})_{u_0}(w)}+\sum_{j\in{\bf Z}}
\frac{m_a^o(\lambda_{a,2j+1})_{u_0}(v)}
{1-(\lambda_{a,2j+1})_{u_0}(w)}\right)}\\
\hspace{0.6truecm}\displaystyle{+\sum_{a\in I_{\widetilde{\sigma}}}
\left(\sum_{j\in{\bf Z}\setminus\{0\}}\frac{m_a^e(\lambda_{a,2j})_{u_0}(v)}
{1-(\lambda_{a,2j})_{u_0}(w)}+\sum_{j\in{\bf Z}}
\frac{m_a^o(\lambda_{a,2j+1})_{u_0}(v)}
{1-(\lambda_{a,2j+1})_{u_0}(w)}\right)}\\
\displaystyle{=\sum_{a\in I\setminus I_{\widetilde{\sigma}}}\left(
(m_a^e+m_a^o)\cot(\frac{\pi(1-(\lambda_a)_{u_0}(w))}{b_a})\right.}\\
\hspace{2truecm}\displaystyle{
\left.+(m_a^e-m_a^o){\rm cosec}(\frac{\pi(1-(\lambda_a)_{u_0}(w))}{b_a})\right)
\frac{\pi}{2b_a}(\lambda_a)_{u_0}(v)}\\
\hspace{0.6truecm}\displaystyle{+\sum_{a\in I_{\widetilde{\sigma}}}
\left(\frac{m_a^e(\lambda_a)_{u_0}(v)}{b_a}
\sum_{j\in{\bf Z}\setminus\{0\}}\frac{1}{2j}+
\frac{m_a^o(\lambda_a)_{u_0}(v)}{b_a}
\sum_{j\in{\bf Z}}\frac{1}{2j+1}\right)}\\
\displaystyle{=\sum_{a\in I\setminus I_{\widetilde{\sigma}}}\left(
(m_a^e+m_a^o)\cot(\frac{\pi(1-(\lambda_a)_{u_0}(w))}{b_a})\right.}\\
\hspace{2truecm}\displaystyle{
\left.+(m_a^e-m_a^o){\rm cosec}(\frac{\pi(1-(\lambda_a)_{u_0}(w))}{b_a})\right)
\frac{\pi}{2b_a}(\lambda_a)_{u_0}(v),}
\end{array}$$
that is, the $T_{u_0}^{\perp}\widetilde M$-component 
$((\widetilde H^{\widetilde F})_{u_0+w})_{T^{\perp}_{u_0}\widetilde M}$ of 
$(\widetilde H^{\widetilde F})_{u_0+w}$ is equal to 
$$\begin{array}{l}
\displaystyle{\sum_{a\in I\setminus I_{\widetilde{\sigma}}}\left(
(m_a^e+m_a^o)\cot(\frac{\pi(1-(\lambda_a)_{u_0}(w))}{b_a})\right.}\\
\hspace{1.2truecm}\displaystyle{
\left.+(m_a^e-m_a^o){\rm cosec}(\frac{\pi(1-(\lambda_a)_{u_0}(w))}{b_a})\right)
\frac{\pi}{2b_a}({\bf n}_a)_{u_0}.}
\end{array}$$
Denote by $\Phi_{u_0+w}$ the normal 
holonomy group of $\widetilde F$ at $u_0+w$ and $L_{u_0}$ be the focal leaf 
through $u_0$ for $\widetilde w$.  Since $L_{u_0}=\Phi_{u_0+w}\cdot u_0$, 
there exists $\mu\in\Phi_{u_0+w}$ such that 
$\mu(T^{\perp}_{u_0}\widetilde M)=T^{\perp}_{u_1}\widetilde M$ for any point 
$u_1$ of $L_{u_0}$.  From this fact, we have 
$\mu(((\widetilde H^{\widetilde F})_{u_0+w})_{T^{\perp}_{u_0}\widetilde M})
=((\widetilde H^{\widetilde F})_{u_0+w})_{T^{\perp}_{u_1}\widetilde M}$.  
Furthermore, from this fact, we have \newline
$\displaystyle{(\widetilde H^{\widetilde F})_{u_0+w}
\in\mathop{\cap}_{u\in L_{u_0}}T^{\perp}_u\widetilde M}$, which contains 
$\widetilde{\sigma}$ as an open subset.  Therefore, we obtain 
$$\begin{array}{l}
\displaystyle{(\widetilde H^{\widetilde F})_{u_0+w}=
\sum_{a\in I\setminus I_{\widetilde{\sigma}}}\left(
(m_a^e+m_a^o)\cot(\frac{\pi(1-(\lambda_a)_{u_0}(w))}{b_a})\right.}\\
\hspace{2.7truecm}\displaystyle{\left.+(m_a^e-m_a^o){\rm cosec}
(\frac{\pi(1-(\lambda_a)_{u_0}(w))}{b_a})\right)
\frac{\pi}{2b_a}({\bf n}_a)_{u_0}\,\,\,\,(\in T\widetilde{\sigma}).}
\end{array}\leqno{(3.5)}$$
Define a tangent vector field $X^{\widetilde{\sigma}}$ on 
$\widetilde{\sigma}$ by $X^{\widetilde{\sigma}}_w
:=(\widetilde H^{\widetilde F})_{u_0+w}$ ($w\in\widetilde{\sigma}$).  
Fix $w\in\widetilde{\sigma}$.  
Let $\xi:[0,T)\to\widetilde{\sigma}$ be the (maximal) integral curve of 
$X^{\widetilde{\sigma}}$ with $\xi(0)=w$, where $T$ is possible to be equal 
to $\infty$.  
Since $\xi(t)-w\in\widetilde{\sigma}\,
(\subset T^{\perp}_{u_0+w}\widetilde F)$, $\Phi_{u_0+w}$ preserves 
$\widetilde{\sigma}$ invariantly and it acts on $\widetilde{\sigma}$ 
trivially, $\xi(t)-w$ extends to a parallel normal vector field of 
$\widetilde F$.  Denote by $\widetilde{\xi(t)-w}$ this parallel normal vector 
field and by $\widehat{\xi(t)-w}:=(\pi\circ\phi)_{\ast}
(\widetilde{\xi(t)-w})$, which is a parallel normal vector field of $F$.  
Also, denote by $\widetilde F_{\widetilde{\xi(t)-w}}$ the parallel 
submanifold of $\widetilde F$ for $\widetilde{\xi(t)-w}$ and by 
$F_{\widehat{\xi(t)-w}}$ the parallel submanifold of $F$ for 
$\widehat{\xi(t)-w}$.  Let $F_t$ (resp. $\widetilde F_t$) be the mean 
curvature flow having $F$ (resp. $\widetilde F$) as initial data.  Then, by 
imitating the proof of Lemma 3.1, we can show $\widetilde F_t
=\widetilde F_{\widetilde{\xi(t)-w}}$ and 
$F_t=F_{\widehat{\xi(t)-w}}$ ($t\in[0,T)$).  By using 
these facts and $(3.5)$ and imitating the proof of the statement (i) of 
Theorem A, we can show the statement (i).  Also, by imitating the proof of 
the statement (ii) of Theorem A, we can show the statement (ii).  
\hspace{0.5truecm}q.e.d.

\vspace{1truecm}

\centerline{
\unitlength 0.1in
\begin{picture}( 80.9000, 28.8100)(-20.9000,-34.4100)
%
\special{pn 8}%
\special{pa 2242 2220}%
\special{pa 2248 2250}%
\special{pa 2244 2282}%
\special{pa 2234 2314}%
\special{pa 2218 2340}%
\special{pa 2200 2366}%
\special{pa 2178 2390}%
\special{pa 2156 2412}%
\special{pa 2132 2434}%
\special{pa 2104 2452}%
\special{pa 2078 2470}%
\special{pa 2050 2486}%
\special{pa 2022 2500}%
\special{pa 1994 2514}%
\special{pa 1964 2526}%
\special{pa 1934 2538}%
\special{pa 1904 2548}%
\special{pa 1872 2556}%
\special{pa 1842 2564}%
\special{pa 1810 2570}%
\special{pa 1778 2574}%
\special{pa 1746 2578}%
\special{pa 1714 2580}%
\special{pa 1682 2580}%
\special{pa 1652 2576}%
\special{pa 1620 2572}%
\special{pa 1588 2564}%
\special{pa 1558 2552}%
\special{pa 1530 2538}%
\special{pa 1508 2516}%
\special{pa 1490 2490}%
\special{pa 1486 2480}%
\special{sp}%
%
\special{pn 8}%
\special{pa 1486 2480}%
\special{pa 1480 2450}%
\special{pa 1484 2418}%
\special{pa 1494 2388}%
\special{pa 1510 2360}%
\special{pa 1528 2334}%
\special{pa 1550 2310}%
\special{pa 1574 2288}%
\special{pa 1598 2266}%
\special{pa 1624 2248}%
\special{pa 1650 2232}%
\special{pa 1678 2214}%
\special{pa 1706 2200}%
\special{pa 1736 2186}%
\special{pa 1766 2174}%
\special{pa 1794 2162}%
\special{pa 1826 2152}%
\special{pa 1856 2144}%
\special{pa 1888 2136}%
\special{pa 1918 2130}%
\special{pa 1950 2126}%
\special{pa 1982 2122}%
\special{pa 2014 2120}%
\special{pa 2046 2120}%
\special{pa 2078 2124}%
\special{pa 2110 2128}%
\special{pa 2140 2136}%
\special{pa 2170 2148}%
\special{pa 2198 2162}%
\special{pa 2222 2184}%
\special{pa 2240 2210}%
\special{pa 2242 2220}%
\special{sp -0.045}%
%
\special{pn 8}%
\special{ar 1600 1470 400 200  0.1721908 0.2001811}%
%
\special{pn 8}%
\special{pa 2398 3098}%
\special{pa 2396 3130}%
\special{pa 2384 3158}%
\special{pa 2366 3186}%
\special{pa 2344 3208}%
\special{pa 2318 3228}%
\special{pa 2292 3246}%
\special{pa 2262 3260}%
\special{pa 2234 3274}%
\special{pa 2204 3284}%
\special{pa 2174 3294}%
\special{pa 2142 3302}%
\special{pa 2110 3308}%
\special{pa 2080 3314}%
\special{pa 2048 3318}%
\special{pa 2016 3322}%
\special{pa 1984 3322}%
\special{pa 1952 3322}%
\special{pa 1920 3320}%
\special{pa 1888 3318}%
\special{pa 1856 3316}%
\special{pa 1824 3310}%
\special{pa 1794 3302}%
\special{pa 1762 3294}%
\special{pa 1732 3282}%
\special{pa 1704 3270}%
\special{pa 1676 3254}%
\special{pa 1650 3234}%
\special{pa 1628 3212}%
\special{pa 1610 3184}%
\special{pa 1602 3154}%
\special{pa 1600 3146}%
\special{sp}%
%
\special{pn 8}%
\special{pa 1600 3146}%
\special{pa 1604 3114}%
\special{pa 1616 3084}%
\special{pa 1634 3058}%
\special{pa 1656 3036}%
\special{pa 1682 3016}%
\special{pa 1708 2998}%
\special{pa 1736 2982}%
\special{pa 1766 2970}%
\special{pa 1796 2958}%
\special{pa 1826 2950}%
\special{pa 1856 2942}%
\special{pa 1888 2934}%
\special{pa 1920 2928}%
\special{pa 1952 2924}%
\special{pa 1984 2922}%
\special{pa 2016 2922}%
\special{pa 2048 2920}%
\special{pa 2080 2922}%
\special{pa 2112 2924}%
\special{pa 2144 2928}%
\special{pa 2174 2934}%
\special{pa 2206 2940}%
\special{pa 2236 2950}%
\special{pa 2266 2960}%
\special{pa 2296 2974}%
\special{pa 2324 2990}%
\special{pa 2348 3008}%
\special{pa 2372 3032}%
\special{pa 2388 3058}%
\special{pa 2398 3088}%
\special{pa 2398 3098}%
\special{sp -0.045}%
%
\special{pn 8}%
\special{ar 700 3660 930 2210  5.0709435 6.1828940}%
%
\special{pn 8}%
\special{ar 1130 3860 1320 3060  5.1600065 6.1237801}%
%
\special{pn 8}%
\special{ar 950 3450 1060 2290  5.0852484 6.2791495}%
%
\special{pn 20}%
\special{sh 1}%
\special{ar 1880 2360 10 10 0  6.28318530717959E+0000}%
\special{sh 1}%
\special{ar 1590 1630 10 10 0  6.28318530717959E+0000}%
\special{sh 1}%
\special{ar 2000 3130 10 10 0  6.28318530717959E+0000}%
\special{sh 1}%
\special{ar 2000 3130 10 10 0  6.28318530717959E+0000}%
%
\special{pn 8}%
\special{pa 910 2240}%
\special{pa 580 3060}%
\special{fp}%
%
\special{pn 8}%
\special{pa 580 3060}%
\special{pa 2920 2420}%
\special{fp}%
%
\special{pn 8}%
\special{pa 2920 2420}%
\special{pa 3230 1650}%
\special{fp}%
\special{pa 3230 1650}%
\special{pa 3230 1650}%
\special{fp}%
%
\special{pn 8}%
\special{pa 3230 1660}%
\special{pa 2190 1920}%
\special{fp}%
%
\special{pn 8}%
\special{pa 920 2240}%
\special{pa 1340 2140}%
\special{fp}%
%
\special{pn 8}%
\special{pa 1400 2120}%
\special{pa 2130 1940}%
\special{dt 0.045}%
%
\special{pn 8}%
\special{pa 1916 1440}%
\special{pa 1924 1472}%
\special{pa 1926 1502}%
\special{pa 1920 1534}%
\special{pa 1908 1564}%
\special{pa 1892 1592}%
\special{pa 1874 1618}%
\special{pa 1856 1644}%
\special{pa 1834 1668}%
\special{pa 1810 1690}%
\special{pa 1786 1710}%
\special{pa 1762 1730}%
\special{pa 1734 1748}%
\special{pa 1708 1766}%
\special{pa 1680 1782}%
\special{pa 1652 1798}%
\special{pa 1624 1812}%
\special{pa 1594 1826}%
\special{pa 1564 1836}%
\special{pa 1534 1848}%
\special{pa 1504 1856}%
\special{pa 1472 1864}%
\special{pa 1442 1870}%
\special{pa 1410 1874}%
\special{pa 1378 1876}%
\special{pa 1346 1876}%
\special{pa 1314 1870}%
\special{pa 1282 1864}%
\special{pa 1254 1852}%
\special{pa 1226 1834}%
\special{pa 1206 1810}%
\special{pa 1202 1802}%
\special{sp}%
%
\special{pn 8}%
\special{pa 1202 1802}%
\special{pa 1192 1772}%
\special{pa 1192 1740}%
\special{pa 1196 1708}%
\special{pa 1208 1680}%
\special{pa 1222 1650}%
\special{pa 1240 1624}%
\special{pa 1260 1600}%
\special{pa 1282 1576}%
\special{pa 1306 1554}%
\special{pa 1330 1532}%
\special{pa 1354 1512}%
\special{pa 1382 1494}%
\special{pa 1408 1476}%
\special{pa 1436 1462}%
\special{pa 1464 1446}%
\special{pa 1492 1432}%
\special{pa 1522 1418}%
\special{pa 1552 1408}%
\special{pa 1582 1396}%
\special{pa 1612 1386}%
\special{pa 1644 1380}%
\special{pa 1674 1374}%
\special{pa 1706 1370}%
\special{pa 1738 1368}%
\special{pa 1770 1368}%
\special{pa 1802 1372}%
\special{pa 1834 1380}%
\special{pa 1862 1390}%
\special{pa 1890 1408}%
\special{pa 1910 1432}%
\special{pa 1916 1440}%
\special{sp -0.045}%
%
\special{pn 8}%
\special{pa 1410 2830}%
\special{pa 2350 1880}%
\special{dt 0.045}%
%
\special{pn 8}%
\special{pa 1420 2830}%
\special{pa 1660 2580}%
\special{fp}%
%
\special{pn 20}%
\special{sh 1}%
\special{ar 1660 2570 10 10 0  6.28318530717959E+0000}%
\special{sh 1}%
\special{ar 1660 2570 10 10 0  6.28318530717959E+0000}%
%
\special{pn 8}%
\special{pa 2640 2490}%
\special{pa 1110 2200}%
\special{dt 0.045}%
%
\special{pn 8}%
\special{pa 2630 2490}%
\special{pa 2160 2400}%
\special{fp}%
%
\special{pn 20}%
\special{sh 1}%
\special{ar 2170 2410 10 10 0  6.28318530717959E+0000}%
\special{sh 1}%
\special{ar 2170 2410 10 10 0  6.28318530717959E+0000}%
\put(22.5000,-24.5000){\makebox(0,0)[rt]{$u_0$}}%
\put(17.9000,-26.2000){\makebox(0,0)[rt]{$u_1$}}%
%
\special{pn 8}%
\special{pa 2220 2020}%
\special{pa 2350 1880}%
\special{fp}%
%
\special{pn 8}%
\special{pa 1380 2250}%
\special{pa 1140 2200}%
\special{fp}%
%
\special{pn 8}%
\special{pa 2380 1200}%
\special{pa 1950 1760}%
\special{dt 0.045}%
\special{sh 1}%
\special{pa 1950 1760}%
\special{pa 2006 1720}%
\special{pa 1982 1718}%
\special{pa 1976 1696}%
\special{pa 1950 1760}%
\special{fp}%
%
\special{pn 8}%
\special{pa 2220 1700}%
\special{pa 1880 2350}%
\special{dt 0.045}%
\special{sh 1}%
\special{pa 1880 2350}%
\special{pa 1930 2300}%
\special{pa 1906 2304}%
\special{pa 1894 2282}%
\special{pa 1880 2350}%
\special{fp}%
%
\special{pn 8}%
\special{pa 580 2310}%
\special{pa 900 2730}%
\special{dt 0.045}%
\special{sh 1}%
\special{pa 900 2730}%
\special{pa 876 2666}%
\special{pa 868 2688}%
\special{pa 844 2690}%
\special{pa 900 2730}%
\special{fp}%
%
\special{pn 8}%
\special{pa 1600 740}%
\special{pa 1420 1410}%
\special{dt 0.045}%
\special{sh 1}%
\special{pa 1420 1410}%
\special{pa 1458 1352}%
\special{pa 1434 1358}%
\special{pa 1418 1340}%
\special{pa 1420 1410}%
\special{fp}%
%
\special{pn 8}%
\special{pa 3120 1420}%
\special{pa 2470 2450}%
\special{dt 0.045}%
\special{sh 1}%
\special{pa 2470 2450}%
\special{pa 2522 2404}%
\special{pa 2498 2406}%
\special{pa 2490 2384}%
\special{pa 2470 2450}%
\special{fp}%
%
\special{pn 8}%
\special{pa 770 1860}%
\special{pa 1510 2720}%
\special{dt 0.045}%
\special{sh 1}%
\special{pa 1510 2720}%
\special{pa 1482 2656}%
\special{pa 1476 2680}%
\special{pa 1452 2684}%
\special{pa 1510 2720}%
\special{fp}%
%
\special{pn 8}%
\special{pa 700 1320}%
\special{pa 1790 2580}%
\special{dt 0.045}%
\special{sh 1}%
\special{pa 1790 2580}%
\special{pa 1762 2516}%
\special{pa 1756 2540}%
\special{pa 1732 2544}%
\special{pa 1790 2580}%
\special{fp}%
\put(16.9000,-5.6000){\makebox(0,0)[rt]{$\widetilde F$}}%
\put(25.3000,-10.2000){\makebox(0,0)[rt]{$\widetilde M$}}%
\put(25.9000,-15.3000){\makebox(0,0)[rt]{$u_0+w$}}%
\put(32.9000,-12.1000){\makebox(0,0)[rt]{$T^{\perp}_{u_0}\widetilde M$}}%
\put(8.0000,-16.5000){\makebox(0,0)[rt]{$T^{\perp}_{u_1}\widetilde M$}}%
\put(7.8000,-11.5000){\makebox(0,0)[rt]{$L_{u_0}$}}%
\put(7.0000,-20.7000){\makebox(0,0)[rt]{$T^{\perp}_{u_0+w}\widetilde F$}}%
%
\special{pn 8}%
\special{pa 5000 2000}%
\special{pa 4000 2400}%
\special{pa 4000 3000}%
\special{pa 5000 2600}%
\special{pa 5000 2000}%
\special{pa 5000 2000}%
\special{pa 5000 2000}%
\special{fp}%
%
\special{pn 8}%
\special{pa 5000 2000}%
\special{pa 5000 2600}%
\special{pa 6000 3000}%
\special{pa 6000 3000}%
\special{pa 6000 2400}%
\special{pa 5000 2000}%
\special{pa 5000 2000}%
\special{pa 5000 2000}%
\special{fp}%
%
\special{pn 8}%
\special{sh 0.300}%
\special{pa 5000 2100}%
\special{pa 5000 2540}%
\special{pa 5410 2430}%
\special{pa 5410 2430}%
\special{pa 5000 2100}%
\special{fp}%
%
\special{pn 4}%
\special{pa 5000 2150}%
\special{pa 5018 2114}%
\special{dt 0.027}%
\special{pa 5000 2240}%
\special{pa 5050 2140}%
\special{dt 0.027}%
\special{pa 5000 2330}%
\special{pa 5082 2166}%
\special{dt 0.027}%
\special{pa 5000 2420}%
\special{pa 5114 2192}%
\special{dt 0.027}%
\special{pa 5000 2510}%
\special{pa 5146 2218}%
\special{dt 0.027}%
\special{pa 5036 2532}%
\special{pa 5178 2244}%
\special{dt 0.027}%
\special{pa 5088 2518}%
\special{pa 5210 2270}%
\special{dt 0.027}%
\special{pa 5140 2504}%
\special{pa 5242 2296}%
\special{dt 0.027}%
\special{pa 5192 2490}%
\special{pa 5276 2322}%
\special{dt 0.027}%
\special{pa 5244 2476}%
\special{pa 5308 2348}%
\special{dt 0.027}%
\special{pa 5296 2462}%
\special{pa 5340 2374}%
\special{dt 0.027}%
\special{pa 5346 2448}%
\special{pa 5372 2398}%
\special{dt 0.027}%
\special{pa 5398 2434}%
\special{pa 5404 2424}%
\special{dt 0.027}%
%
\special{pn 13}%
\special{sh 1}%
\special{ar 5190 2380 10 10 0  6.28318530717959E+0000}%
\special{sh 1}%
\special{ar 5190 2380 10 10 0  6.28318530717959E+0000}%
%
\special{pn 13}%
\special{pa 5000 2100}%
\special{pa 4990 2540}%
\special{fp}%
%
\special{pn 13}%
\special{pa 4980 2360}%
\special{pa 4990 2150}%
\special{fp}%
\special{sh 1}%
\special{pa 4990 2150}%
\special{pa 4968 2216}%
\special{pa 4988 2204}%
\special{pa 5008 2218}%
\special{pa 4990 2150}%
\special{fp}%
%
\special{pn 13}%
\special{sh 1}%
\special{ar 4990 2370 10 10 0  6.28318530717959E+0000}%
\special{sh 1}%
\special{ar 4990 2370 10 10 0  6.28318530717959E+0000}%
%
\special{pn 8}%
\special{pa 5280 1900}%
\special{pa 5000 2050}%
\special{dt 0.045}%
\special{sh 1}%
\special{pa 5000 2050}%
\special{pa 5068 2036}%
\special{pa 5048 2026}%
\special{pa 5050 2002}%
\special{pa 5000 2050}%
\special{fp}%
%
\special{pn 8}%
\special{pa 4760 2480}%
\special{pa 4990 2460}%
\special{dt 0.045}%
\special{sh 1}%
\special{pa 4990 2460}%
\special{pa 4922 2446}%
\special{pa 4938 2466}%
\special{pa 4926 2486}%
\special{pa 4990 2460}%
\special{fp}%
%
\special{pn 8}%
\special{pa 4670 2990}%
\special{pa 4340 2680}%
\special{dt 0.045}%
\special{sh 1}%
\special{pa 4340 2680}%
\special{pa 4376 2740}%
\special{pa 4380 2718}%
\special{pa 4402 2712}%
\special{pa 4340 2680}%
\special{fp}%
%
\special{pn 8}%
\special{pa 5500 3030}%
\special{pa 5740 2670}%
\special{dt 0.045}%
\special{sh 1}%
\special{pa 5740 2670}%
\special{pa 5686 2714}%
\special{pa 5710 2714}%
\special{pa 5720 2738}%
\special{pa 5740 2670}%
\special{fp}%
%
\special{pn 8}%
\special{pa 5350 2560}%
\special{pa 5200 2390}%
\special{dt 0.045}%
\special{sh 1}%
\special{pa 5200 2390}%
\special{pa 5230 2454}%
\special{pa 5236 2430}%
\special{pa 5260 2428}%
\special{pa 5200 2390}%
\special{fp}%
\put(60.0000,-16.7000){\makebox(0,0)[rt]{$\displaystyle{\mathop{\cap}_{u\in L_{u_0}}T^{\perp}_u\widetilde M}$}}%
\put(48.5000,-30.4000){\makebox(0,0)[rt]{$T^{\perp}_{u_1}\widetilde M$}}%
\put(56.4000,-30.8000){\makebox(0,0)[rt]{$T^{\perp}_{u_0}\widetilde M$}}%
\put(46.8000,-24.0000){\makebox(0,0)[rt]{$\widetilde{\sigma}$}}%
\put(54.5000,-26.1000){\makebox(0,0)[rt]{$u_0$}}%
\put(49.2000,-21.6000){\makebox(0,0)[rt]{$(\widetilde H^{\widetilde F})_{u_0+w}$}}%
%
\special{pn 8}%
\special{pa 5780 2090}%
\special{pa 5100 2260}%
\special{dt 0.045}%
\special{sh 1}%
\special{pa 5100 2260}%
\special{pa 5170 2264}%
\special{pa 5152 2248}%
\special{pa 5160 2224}%
\special{pa 5100 2260}%
\special{fp}%
\put(59.5000,-20.0000){\makebox(0,0)[rt]{$\widetilde C$}}%
\end{picture}%
\hspace{6truecm}}

\vspace{1truecm}

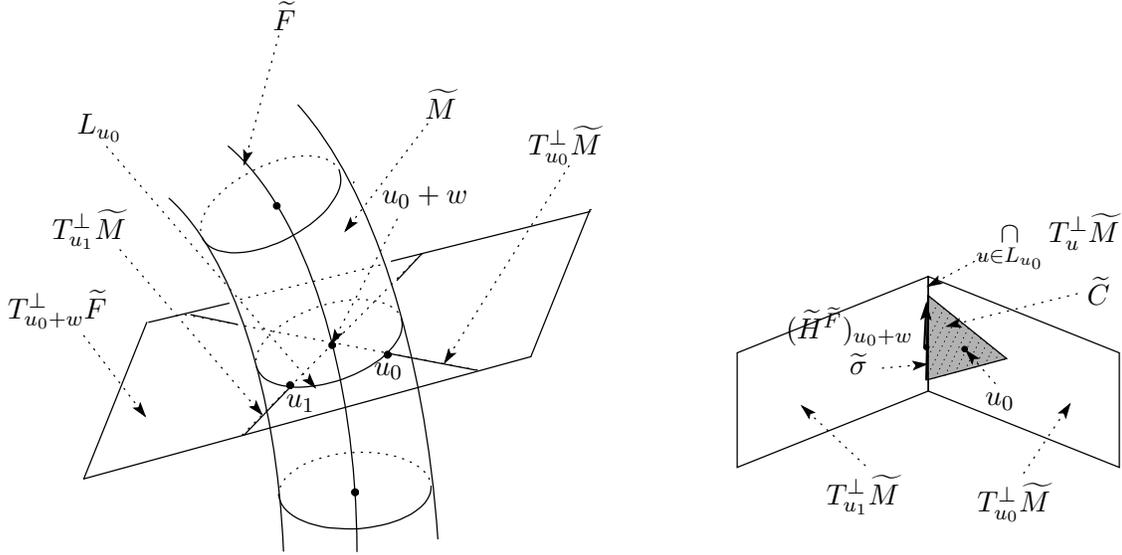
\centerline{{\bf Fig. 3.}}

\section{The mean curvature flows of orbits of Hermann actions} 
In this section, we describe explicitly the mean curvature flows of orbits of 
Hermann actions whose principal orbits are full and irreducible.  
Let $G/K$ be a symmetric space of compact type and $H$ be 
a symmetric subgroup of $G$.  Also, let $\theta$ be an involution of 
$G$ with $({\rm Fix}\,\theta)_0\subset K\subset{\rm Fix}\,\theta$ and $\tau$ 
be an invloution of $G$ with $({\rm Fix}\,\tau)_0\subset H\subset{\rm Fix}\,
\tau$, where ${\rm Fix}\,\theta$ (resp. ${\rm Fix}\,\tau$) is the fixed point 
group of $\theta$ (resp. $\tau$) and $({\rm Fix}\,\theta)_0$ (resp. 
$({\rm Fix}\,\tau)_0$) is the identity component of ${\rm Fix}\,\theta$ 
(resp. ${\rm Fix}\,\tau$).  We assume that 
$\tau\circ\theta=\theta\circ\tau$.  Set $L:={\rm Fix}(\theta\circ\tau)$.  
Denote by the same symbol $\theta$ (resp. $\tau$) the involution of the Lie 
algebra $\mathfrak g$ of $G$ induced from $\theta$ (resp. $\tau$).  Set 
$\mathfrak k:={\rm Ker}(\theta-{\rm id}),\,\mathfrak p:={\rm Ker}(\theta+
{\rm id}),\,\mathfrak h:={\rm Ker}(\tau-{\rm id})$ and $\mathfrak q:=
{\rm Ker}(\tau+{\rm id})$.  The space $\mathfrak p$ is identified with 
$T_{eK}(G/K)$.  From $\theta\circ\tau=\tau\circ\theta$, we have 
$\mathfrak p=\mathfrak p\cap\mathfrak h+\mathfrak p\cap\mathfrak q$.  Take 
a maximal abelian subspace $\mathfrak b$ of $\mathfrak p\cap\mathfrak q$ and 
let $\mathfrak p=\mathfrak z_{\mathfrak p}(\mathfrak b)
+\sum\limits_{\beta\in\triangle'_+}\mathfrak p_{\beta}$ be the root space 
decomposition with respect to $\mathfrak b$, where 
$\mathfrak z_{\mathfrak p}(\mathfrak b)$ is the centralizer of $\mathfrak b$ 
in $\mathfrak p$, $\triangle'_+$ is the positive root system of 
$\triangle':=\{\beta\in\mathfrak b^{\ast}\setminus\{0\}
\,\vert\,\exists\,X(\not=0)\in
\mathfrak p\,\,{\rm s.t.}\,\,{\rm ad}(b)^2(X)=-\beta(b)^2X\,\,(\forall\,b\in
\mathfrak b)\}$ under some lexicographic ordering of $\mathfrak b^{\ast}$ and 
$\mathfrak p_{\beta}:=\{X\in\mathfrak p\,\vert\,{\rm ad}(b)^2(X)=-\beta(b)^2X
\,\,(\forall\,b\in\mathfrak b)\}$ ($\beta\in\triangle'_+$). Also, let 
${\triangle'}^V_+:=\{\beta\in\triangle'_+\,\vert\,\mathfrak p_{\beta}\cap
\mathfrak q\not=\{0\}\}$ and ${\triangle'}^H_+:=\{\beta\in\triangle'_+\,\vert\,
\mathfrak p_{\beta}\cap\mathfrak h\not=\{0\}\}$.  
Then we have $\mathfrak q=\mathfrak b+\sum\limits_{\beta\in{\triangle'}^V_+}
(\mathfrak p_{\beta}\cap\mathfrak q)$ and $\mathfrak h=
\mathfrak z_{\mathfrak h}(\mathfrak b)+\sum\limits_{\beta\in{\triangle'}^H_+}
(\mathfrak p_{\beta}\cap\mathfrak h)$, where $\mathfrak z_{\mathfrak h}
(\mathfrak b)$ is the centralizer of $\mathfrak b$ in $\mathfrak h$.  The 
orbit $H(eK)$ is a reflective submanifold and it is isometric to the symmetric 
space $H/H\cap K$ (equipped with a metric scaled suitably).  Also, 
$\exp^{\perp}(T^{\perp}_{eK}(H(eK)))$ is also a reflective submanifold and it 
is isometric to the symmetric space $L/H\cap K$ (equipped with a metric scaled 
suitably), where $\exp^{\perp}$ is the normal exponential map of $H(eK)$.  
The system ${\triangle'}^V:={\triangle'}^V_+\cup(-{\triangle'}^V_+)$ is the 
root system of $L/H\cap K$.  Define a subset $\widetilde C$ of $\mathfrak b$ 
by 
$$\begin{array}{l}
\displaystyle{\widetilde C:=\{b\in\mathfrak b\,\vert\,0<\beta(b)<\pi\,
(\forall\,\beta\in{\triangle'}^V_+\setminus{\triangle'}^H_+),\,\,
-\frac{\pi}{2}<\beta(b)<\frac{\pi}{2}\,
(\forall\,\beta\in{\triangle'}^H_+\setminus{\triangle'}^V_+),}\\
\hspace{2.2truecm}\displaystyle{
0<\beta(b)<\frac{\pi}{2}\,
(\forall\,\beta\in{\triangle'}^V_+\cap{\triangle'}^H_+)
\}.}
\end{array}$$
Let $\Pi$ be the simple root system of $\triangle'_+$, and set $\Pi_V:=\Pi\cap
{\triangle'}^V_+$ and $\Pi_H:=\Pi\cap{\triangle'}^H_+$.  Also, let $\delta$ be 
the highest root of ${\triangle'}^V_+\cup2{\triangle'}^H_+$.  
Set $\widetilde{\Pi}:=\Pi\cup\{\delta\},\,\widetilde{\Pi}_V:=
\widetilde{\Pi}\cap{\triangle'_+}^V$ and 
$\widetilde{\Pi}_H:=
\widetilde{\Pi}\cap({\triangle'_+}^H\cup2{\triangle'_+}^H)$.  
Then we have 
$$\widetilde C=\{b\in\mathfrak b\,\vert\,\beta(b)>0\,(\forall\,\beta\in\Pi_V),
\,\,\beta(b)>-\frac{\pi}{2}\,(\forall\,\beta\in\Pi_H\setminus\Pi_V),\,\,
\delta(b)<\pi\},$$
Set $C:={\rm Exp}(\widetilde C)$, where ${\rm Exp}$ is the exponential map 
of $G/K$ at $eK$.  
Let $P(G,H\times K):=\{g\in H^1([0,1],G)\,\vert\,(g(0),g(1))\in H\times K\}$, 
where $H^1([0,1],G)$ is the Hilbert Lie group of all $H^1$-paths in $G$.  
This group acts on $H^0([0,1],\mathfrak g)$ as gauge action.  The orbits of 
the $P(G,H\times K)$-action are the inverse images of orbits of the $H$-action 
by $\pi\circ\phi$.  The set $\Sigma:={\rm Exp}(\mathfrak b)$ is a section of 
the $H$-action and $\mathfrak b$ is a section of the $P(G,H\times K)$-action 
on $H^0([0,1],\mathfrak g)$, where $\mathfrak b$ is identified with the 
horizontal lift of $\mathfrak b$ to the zero element $\hat 0$ of 
$H^0([0,1],\mathfrak g)$ ($\hat0\,:\,$the constant path at the zero element 
$0$ of $\mathfrak g$).  
The set $\widetilde C$ is the fundamental domain of the Coxeter group of 
a principal $P(G,H\times K)$-orbit (hence a principal $H$-orbit) and 
each prinicipal $H$-orbit meets $C$ at one point and each singular $H$-orbit 
meets $\partial C$ at one point.  
The focal set of the principal orbit $P(G,H\times K)\cdot Z_0$ 
($Z_0\in\widetilde C$) at $Z_0$ consists of the hyperplanes 
$\beta^{-1}(j\pi)$'s 
($\beta\in{\triangle'}^V_+\setminus{\triangle'}^H_+,\,\,j\in{\Bbb Z}$), 
$\beta^{-1}((j+\frac12)\pi)$'s 
($\beta\in{\triangle'}^H_+\setminus{\triangle'}^V_+,\,\,j\in{\Bbb Z}$), 
$\beta^{-1}(\frac{j\pi}{2})$'s 
($\beta\in{\triangle'}^V_+\cap{\triangle'}^H_+,\,\,j\in{\Bbb Z}$) in 
$\mathfrak b(=T^{\perp}_{Z_0}(P(G,H\times K)\cdot Z_0))$.  
Denote by $\exp^G$ the exponential map of $G$.  Note that 
$\pi\circ\exp^G\vert_{\mathfrak p}={\rm Exp}$.  
Let $Y_0\in\widetilde C$ and $M(Y_0):=H({\rm Exp}(Y_0))$.  
Then we have $T^{\perp}_{{\rm Exp}(Y_0)}M(Y_0)
=(\exp^G(Y_0))_{\ast}(\mathfrak b)$.  Denote by $A^{Y_0}$ the shape tensor of 
$M(Y_0)$.  Take $v\in T^{\perp}_{{\rm Exp}(Y_0)}M(Y_0)$ and set 
$\bar v:=(\exp^G(Y_0))_{\ast}^{-1}(v)$.  By scaling the metric of $G/K$ by 
a suitable positive constant, we have 
$$A^{Y_0}_v\vert_{\exp^G(Y_0)_{\ast}(\mathfrak p_{\beta}\cap\mathfrak q)}=
-\frac{\beta(\bar v)}{\tan\beta(Y_0)}{\rm id}\,\,\,\,
(\beta\in{\triangle'}^V_+)\leqno{(4.1)}$$
and 
$$A^{Y_0}_v\vert_{\exp^G(Y_0)_{\ast}(\mathfrak p_{\beta}\cap\mathfrak h)}=
\beta(\bar v)\tan\beta(Y_0){\rm id}\,\,\,\,(\beta\in{\triangle'}^H_+).
\leqno{(4.2)}$$
Set $m_{\beta}^V:={\rm dim}(\mathfrak p_{\beta}\cap\mathfrak q)$ 
($\beta\in{\triangle'}^V_+$) and 
$m_{\beta}^H:={\rm dim}(\mathfrak p_{\beta}\cap\mathfrak h)$ 
($\beta\in{\triangle'}^H_+$).  
Set 
$\widetilde M(Y_0):=(\pi\circ\phi)^{-1}(M(Y_0))(=P(G,H\times K)\cdot Y_0)$.  
We can show $(\pi\circ\phi)(Y_0)={\rm Exp}(Y_0)$.  

\vspace{1truecm}

\centerline{
\unitlength 0.1in
\begin{picture}( 52.8000, 31.7100)( 15.2000,-35.7100)
%
\special{pn 8}%
\special{pa 2400 1000}%
\special{pa 1800 1600}%
\special{pa 3600 1600}%
\special{pa 4200 1000}%
\special{pa 4200 1000}%
\special{pa 2400 1000}%
\special{fp}%
%
\special{pn 8}%
\special{pa 5000 1000}%
\special{pa 4400 1600}%
\special{pa 6200 1600}%
\special{pa 6800 1000}%
\special{pa 6800 1000}%
\special{pa 5000 1000}%
\special{fp}%
%
\special{pn 8}%
\special{pa 2390 1400}%
\special{pa 3290 1530}%
\special{pa 3230 1120}%
\special{pa 3230 1120}%
\special{pa 2390 1400}%
\special{fp}%
%
\special{pn 8}%
\special{pa 5030 1400}%
\special{pa 5930 1530}%
\special{pa 5870 1120}%
\special{pa 5870 1120}%
\special{pa 5030 1400}%
\special{fp}%
%
\special{pn 13}%
\special{sh 1}%
\special{ar 2390 1400 10 10 0  6.28318530717959E+0000}%
\special{sh 1}%
\special{ar 2390 1400 10 10 0  6.28318530717959E+0000}%
%
\special{pn 13}%
\special{sh 1}%
\special{ar 5750 1350 10 10 0  6.28318530717959E+0000}%
\special{sh 1}%
\special{ar 5750 1350 10 10 0  6.28318530717959E+0000}%
%
\special{pn 8}%
\special{pa 2390 1400}%
\special{pa 3110 1340}%
\special{fp}%
\special{sh 1}%
\special{pa 3110 1340}%
\special{pa 3042 1326}%
\special{pa 3058 1344}%
\special{pa 3046 1366}%
\special{pa 3110 1340}%
\special{fp}%
%
\special{pn 8}%
\special{ar 3610 3450 1230 1104  3.2423400 4.4011242}%
%
\special{pn 8}%
\special{ar 2352 3856 2024 532  4.7317414 5.7176258}%
%
\special{pn 8}%
\special{ar 5090 3590 1050 900  3.1625477 4.3331796}%
%
\special{pn 8}%
\special{ar 3020 3180 1996 778  4.8178040 5.7162869}%
%
\special{pn 8}%
\special{ar 2820 3790 1160 860  4.8554059 5.6578153}%
%
\special{pn 8}%
\special{ar 3570 3490 1070 660  4.1322560 5.1025259}%
%
\special{pn 8}%
\special{ar 4040 3450 302 584  3.4269447 4.4736244}%
%
\special{pn 13}%
\special{sh 1}%
\special{ar 2990 2940 10 10 0  6.28318530717959E+0000}%
\special{sh 1}%
\special{ar 2990 2940 10 10 0  6.28318530717959E+0000}%
%
\special{pn 8}%
\special{pa 2400 1410}%
\special{pa 2990 2930}%
\special{dt 0.045}%
%
\special{pn 8}%
\special{ar 3020 3240 1170 300  4.6877391 5.3265484}%
%
\special{pn 13}%
\special{sh 1}%
\special{ar 3690 3000 10 10 0  6.28318530717959E+0000}%
\special{sh 1}%
\special{ar 3690 3000 10 10 0  6.28318530717959E+0000}%
%
\special{pn 8}%
\special{pa 3690 3000}%
\special{pa 5740 1350}%
\special{dt 0.045}%
\put(23.3000,-13.5000){\makebox(0,0)[rt]{${\bf 0}$}}%
\put(58.0000,-12.1000){\makebox(0,0)[rt]{${\bf 0}$}}%
%
\special{pn 8}%
\special{pa 2610 840}%
\special{pa 2950 1350}%
\special{dt 0.045}%
\special{sh 1}%
\special{pa 2950 1350}%
\special{pa 2930 1284}%
\special{pa 2920 1306}%
\special{pa 2896 1306}%
\special{pa 2950 1350}%
\special{fp}%
%
\special{pn 8}%
\special{pa 3060 840}%
\special{pa 3160 1240}%
\special{dt 0.045}%
\special{sh 1}%
\special{pa 3160 1240}%
\special{pa 3164 1170}%
\special{pa 3148 1188}%
\special{pa 3124 1180}%
\special{pa 3160 1240}%
\special{fp}%
%
\special{pn 8}%
\special{pa 3470 2290}%
\special{pa 3690 2990}%
\special{dt 0.045}%
\special{sh 1}%
\special{pa 3690 2990}%
\special{pa 3690 2920}%
\special{pa 3674 2940}%
\special{pa 3652 2932}%
\special{pa 3690 2990}%
\special{fp}%
%
\special{pn 8}%
\special{pa 6320 810}%
\special{pa 5830 1430}%
\special{dt 0.045}%
\special{sh 1}%
\special{pa 5830 1430}%
\special{pa 5888 1390}%
\special{pa 5864 1388}%
\special{pa 5856 1366}%
\special{pa 5830 1430}%
\special{fp}%
%
\special{pn 8}%
\special{pa 3940 2180}%
\special{pa 3740 2910}%
\special{dt 0.045}%
\special{sh 1}%
\special{pa 3740 2910}%
\special{pa 3778 2852}%
\special{pa 3754 2860}%
\special{pa 3738 2840}%
\special{pa 3740 2910}%
\special{fp}%
\put(26.4000,-6.6000){\makebox(0,0)[rt]{$Y_0$}}%
\put(31.3000,-6.7000){\makebox(0,0)[rt]{$\widetilde C_0$}}%
\put(64.3000,-6.4000){\makebox(0,0)[rt]{$\widetilde C$}}%
\put(38.1000,-11.0000){\makebox(0,0)[rt]{$\mathfrak b$}}%
\put(70.1000,-4.0000){\makebox(0,0)[rt]{$T^{\perp}_{{\rm Exp}(Y_0)}M$}}%
%
\special{pn 8}%
\special{pa 6750 630}%
\special{pa 6430 1120}%
\special{dt 0.045}%
\special{sh 1}%
\special{pa 6430 1120}%
\special{pa 6484 1076}%
\special{pa 6460 1076}%
\special{pa 6450 1054}%
\special{pa 6430 1120}%
\special{fp}%
\put(36.2000,-21.3000){\makebox(0,0)[rt]{${\rm Exp}(Y_0)$}}%
\put(40.3000,-20.2000){\makebox(0,0)[rt]{$C$}}%
\put(29.2000,-28.9000){\makebox(0,0)[rt]{$eK$}}%
%
\special{pn 8}%
\special{pa 4720 2490}%
\special{pa 4460 2790}%
\special{dt 0.045}%
\special{sh 1}%
\special{pa 4460 2790}%
\special{pa 4520 2754}%
\special{pa 4496 2750}%
\special{pa 4490 2728}%
\special{pa 4460 2790}%
\special{fp}%
\put(50.1000,-23.1000){\makebox(0,0)[rt]{${\rm Exp}(\mathfrak b)$}}%
\end{picture}%
}

\vspace{1truecm}

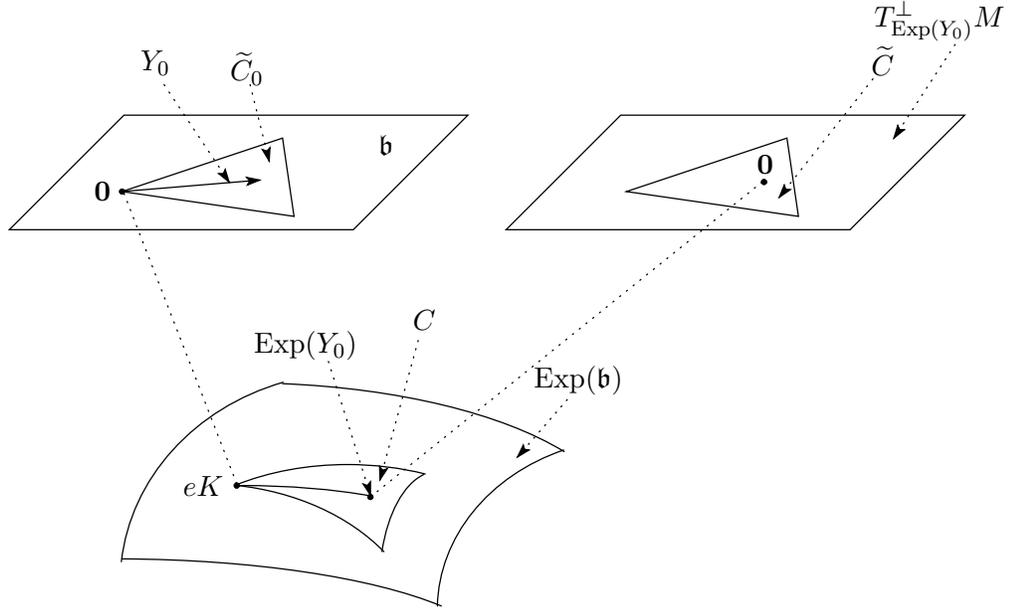
\centerline{{\bf Fig. 4.}}

\vspace{0.5truecm}

\noindent
Denote by ${\widetilde A}^{Y_0}$ the shape tensor of $\widetilde M(Y_0)$.  
According to Proposition 3.2 of [Koi1], we have 
$$\begin{array}{l}
\displaystyle{{\rm Spec}({\widetilde A}^{Y_0}_{\bar v}
\vert_{(\pi\circ\phi)^{-1}_{\ast Y_0}(\exp^G(Y_0)_{\ast}
(\mathfrak p_{\beta}\cap\mathfrak q))})\setminus\{0\}=
\{\frac{-\beta(\bar v)}{\beta(Y_0)+j\pi}\,\vert\,j\in{\Bbb Z}\}\,\,\,\,
(\beta\in{\triangle'}^V_+),}\\
\displaystyle{{\rm Spec}({\widetilde A}^{Y_0}_{\bar v}
\vert_{(\pi\circ\phi)^{-1}_{\ast Y_0}(\exp^G(Y_0)_{\ast}
(\mathfrak p_{\beta}\cap\mathfrak h))})\setminus\{0\}=
\{\frac{-\beta(\bar v)}{\beta(Y_0)+(j+\frac12)\pi}\,\vert\,j\in{\Bbb Z}\}\,\,\,
\,(\beta\in{\triangle'}^H_+),}
\end{array}$$
and 
$${\rm Spec}({\widetilde A}^{Y_0}_{\bar v}
\vert_{(\pi\circ\phi)^{-1}_{\ast Y_0}(\exp^G(Y_0)_{\ast}
(\mathfrak z_{\mathfrak h}(\mathfrak b)))})=\{0\}.$$
Hence the set ${\cal PC}_{\widetilde M(Y_0)}$ of all principal curvatures of 
$\widetilde M(Y_0)$ is given by 
$${\cal PC}_{\widetilde M(Y_0)}=\{\frac{-\widetilde{\beta}}{\beta(Y_0)+j\pi}\,
\vert\,\beta\in{\triangle'}^V_+,\,\,j\in\Bbb Z\}\cup
\{\frac{-\widetilde{\beta}}{\beta(Y_0)+(j+\frac12)\pi}\,\vert\,
\beta\in{\triangle'}^H_+,\,\,j\in\Bbb Z\},$$
where $\widetilde{\beta}$ is the parallel section of $(T^{\perp}\widetilde M
(Y_0))^{\ast}$ with $\widetilde{\beta}_{u_0}
=\beta\circ\exp^G(Y_0)_{\ast}^{-1}$.  Also, we can show that the multiplicity 
of $\frac{-\widetilde{\beta}}
{\beta(Y_0)+j\pi}$ ($\beta\in{\triangle'}^V_+$) is equal to $m_{\beta}^V$ 
and that of 
$\frac{-\widetilde{\beta}}{\beta(Y_0)+(j+\frac12)\pi}$ 
($\beta\in{\triangle'}^H_+$) is equal to $m_{\beta}^H$.  
Define $\lambda_{\beta}^{Y_0}$ and $b_{\beta}^{Y_0}$ 
($\beta\in{\triangle'}_+$) by 
$$(\lambda^{Y_0}_{\beta},b^{Y_0}_{\beta}):=\left\{
\begin{array}{ll}
\displaystyle{(\frac{-\widetilde{\beta}}{\beta(Y_0)},\,\frac{\pi}
{\beta(Y_0)})} & \displaystyle{(\beta\in{\triangle'}^V_+\setminus
{\triangle'}^H_+)}\\
\displaystyle{(\frac{-\widetilde{\beta}}{\beta(Y_0)+\frac{\pi}{2}},\,
\frac{\pi}{\beta(Y_0)+\frac{\pi}{2}})} & \displaystyle{
(\beta\in{\triangle'}^H_+\setminus{\triangle'}^V_+)}\\
\displaystyle{(\frac{-\widetilde{\beta}}{\beta(Y_0)},\,\frac{\pi}
{2\beta(Y_0)})} & 
\displaystyle{(\beta\in{\triangle'}^V_+\cap{\triangle'}^H_+).}\end{array}
\right.$$
Then we have $\frac{-\widetilde{\beta}}{\beta(Y_0)+j\pi}
=\frac{\lambda_{\beta}^{Y_0}}{1+jb_{\beta}^{Y_0}}$ when 
$\beta\in{\triangle'}^V_+\setminus{\triangle'}^H_+$, 
$\frac{-\widetilde{\beta}}{\beta(Y_0)+(j+\frac12)\pi}
=\frac{\lambda_{\beta}^{Y_0}}{1+jb_{\beta}^{Y_0}}$ when 
$\beta\in{\triangle'}^H_+\setminus{\triangle'}^V_+$ and 
$(\frac{-\widetilde{\beta}}{\beta(Y_0)+j\pi},\frac{-\widetilde{\beta}}
{\beta(Y_0)+(j+\frac12)\pi})
=(\frac{\lambda_{\beta}^{Y_0}}{1+2jb_{\beta}^{Y_0}},\,
\frac{\lambda_{\beta}^{Y_0}}{1+(2j+1)b_{\beta}^{Y_0}})$ when 
$\beta\in{\triangle'}^V_+\cap{\triangle'}^H_+$.  That is, we have 
$${\cal PC}_{\widetilde M(Y_0)}=\{\frac{\lambda_{\beta}^{Y_0}}
{1+jb_{\beta}^{Y_0}}\,\vert\,\beta\in\triangle'_+,\,j\in{\Bbb Z}\}.$$
Denote by $m_{\beta j}$ the multiplicity of 
$\frac{\lambda_{\beta}}{1+jb_{\beta}}$.  Then we have 
$$
m_{\beta,2j}=\left\{
\begin{array}{ll}
\displaystyle{m_{\beta}^V} & \displaystyle{(\beta\in{\triangle'}^V_+\setminus
{\triangle'}^H_+)}\\
\displaystyle{m_{\beta}^H} & \displaystyle{(\beta\in{\triangle'}^H_+
\setminus{\triangle'}^V_+)}\\
\displaystyle{m_{\beta}^V} & \displaystyle{(\beta\in{\triangle'}^V_+\cap
{\triangle'}^H_+),}
\end{array}\right.\qquad
m_{\beta,2j+1}=\left\{
\begin{array}{ll}
\displaystyle{m_{\beta}^V} & \displaystyle{(\beta\in{\triangle'}^V_+
\setminus{\triangle'}^H_+)}\\
\displaystyle{m_{\beta}^H} & \displaystyle{(\beta\in{\triangle'}^H_+\setminus
{\triangle'}^V_+)}\\
\displaystyle{m_{\beta}^H} & \displaystyle{(\beta\in{\triangle'}^V_+\cap
{\triangle'}^H_+),}
\end{array}\right.$$
where $j\in{\Bbb Z}$.  
Denote by $\widetilde H^{Y_0}$ the mean curvature vector of 
$\widetilde M(Y_0)$ and ${\bf n}^{Y_0}_{\beta}$ the curvature normal 
corresponding to $\lambda^{Y_0}_{\beta}$.  From $(3.1)$ (the case of $v=0$), 
we have 
$$\begin{array}{l}
\displaystyle{\widetilde H^{Y_0}=\sum_{\beta\in{\triangle'}^V_+\setminus
{\triangle'}^H_+}2m_{\beta}^V\cot\frac{\pi}{b^{Y_0}_{\beta}}\times
\frac{\pi}{2b_{\beta}^{Y_0}}{\bf n}_{\beta}^{Y_0}}\\
\hspace{1.2truecm}\displaystyle{
+\sum_{\beta\in{\triangle'}^H_+\setminus{\triangle'}^V_+}2m_{\beta}^H
\cot\frac{\pi}{b^{Y_0}_{\beta}}\times\frac{\pi}{2b_{\beta}^{Y_0}}
{\bf n}_{\beta}^{Y_0}}\\
\hspace{1.2truecm}\displaystyle{
+\sum_{\beta\in{\triangle'}^V_+\cap{\triangle'}^H_+}
\left((m_{\beta}^V+m_{\beta}^H)\cot\frac{\pi}{b^{Y_0}_{\beta}}+
(m_{\beta}^V-m_{\beta}^H){\rm cosec}\frac{\pi}{b^{Y_0}_{\beta}}\right)
\frac{\pi}{2b_{\beta}^{Y_0}}{\bf n}_{\beta}^{Y_0}}\\
\hspace{0.6truecm}\displaystyle{=
\sum_{\beta\in{\triangle'}^V_+\setminus{\triangle'}^H_+}m_{\beta}^V\beta(Y_0)
\cot\beta(Y_0){\bf n}_{\beta}^{Y_0}}\\
\hspace{1.2truecm}\displaystyle{
-\sum_{\beta\in{\triangle'}^H_+\setminus{\triangle'}^V_+}
m_{\beta}^H(\beta(Y_0)+\frac{\pi}{2})\tan\beta(Y_0){\bf n}_{\beta}^{Y_0}}\\
\hspace{1.2truecm}\displaystyle{
+\sum_{\beta\in{\triangle'}^V_+\cap{\triangle'}^H_+}\beta(Y_0)
\left((m_{\beta}^V+m_{\beta}^H)\cot2\beta(Y_0)\right.}\\
\hspace{3.5truecm}\displaystyle{\left.+
(m_{\beta}^V-m_{\beta}^H){\rm cosec}2\beta(Y_0)\right){\bf n}_{\beta}^{Y_0}.}
\end{array}
\leqno{(4.3)}$$
Define $\beta^{\sharp}\,(\in\mathfrak b)$ by $\beta(\cdot)=\langle
\beta^{\sharp},\cdot\rangle$ and let $\widetilde{\beta^{\sharp}}^{Y_0}$ be 
the parallel normal vector field of $\widetilde M(Y_0)$ with 
$(\widetilde{\beta^{\sharp}}^{Y_0})_{Y_0}=\beta^{\sharp}$, where we identify 
$\mathfrak b$ with $T^{\perp}_{Y_0}\widetilde M(Y_0)$.  From the definition of 
$\lambda^{Y_0}_{\beta}$, we have 
$${\bf n}_{\beta}^{Y_0}=\left\{
\begin{array}{ll}
\displaystyle{-\frac{1}{\beta(Y_0)}\widetilde{\beta^{\sharp}}^{Y_0}} & 
\displaystyle{(\beta\in{\triangle'}^V_+\setminus{\triangle'}^H_+\,\,{\rm or}
\,\,\beta\in{\triangle'}^V_+\cap{\triangle'}^H_+)}\\
\displaystyle{-\frac{1}{\beta(Y_0)+\frac{\pi}{2}}
\widetilde{\beta^{\sharp}}^{Y_0}} & 
\displaystyle{(\beta\in{\triangle'}^H_+\setminus{\triangle'}^V_+).}
\end{array}
\right.$$
Substituting this relation into $(4.3)$, we have 
$$\begin{array}{l}
\displaystyle{
\widetilde H^{Y_0}=-\sum_{\beta\in{\triangle'}^V_+\setminus{\triangle'}^H_+}
m_{\beta}^V\cot\beta(Y_0)\widetilde{\beta^{\sharp}}^{Y_0}
+\sum_{\beta\in{\triangle'}^H_+\setminus{\triangle'}^V_+}
m_{\beta}^H\tan\beta(Y_0)\widetilde{\beta^{\sharp}}^{Y_0}}\\
\hspace{1.2truecm}\displaystyle{
-\sum_{\beta\in{\triangle'}^V_+\cap{\triangle'}^H_+}
\left((m_{\beta}^V+m_{\beta}^H)\cot2\beta(Y_0)\right.}\\
\hspace{3.4truecm}\displaystyle{\left.+(m_{\beta}^V-m_{\beta}^H)
{\rm cosec}2\beta(Y_0)\right)
\widetilde{\beta^{\sharp}}^{Y_0}.}
\end{array}
\leqno{(4.4)}$$
Define a tangent vector field $X$ on $\widetilde C$ by $X_{Y_0}
:=({\widetilde H}^{Y_0})_{Y_0}\,(\in T^{\perp}_{Y_0}\widetilde M(Y_0)
=\mathfrak b(\subset V))$.  From $(4.4)$, we have 
$$\begin{array}{l}
\displaystyle{X_{Y_0}=-\sum_{\beta\in{\triangle'}^V_+\setminus{\triangle'_+}^H}
m_{\beta}^V\cot\beta(Y_0)\beta^{\sharp}
+\sum_{\beta\in{\triangle'}^H_+\setminus{\triangle'}_+^V}
m_{\beta}^H\tan\beta(Y_0)\beta^{\sharp}}\\
\hspace{1.2truecm}\displaystyle{
-\sum_{\beta\in{\triangle'}_+^V\cap{\triangle'}_+^H}
\left((m_{\beta}^V+m_{\beta}^H)\cot2\beta(Y_0)\right.}\\
\hspace{3truecm}\displaystyle{\left.+(m_{\beta}^V-m_{\beta}^H)
{\rm cosec}2\beta(Y_0)\right)\beta^{\sharp}.}
\end{array}\leqno{(4.5)}$$
In order to analyze the mean curvature flows having principal orbits of the 
$H$-actions as initial data, we have only to analyze this vector field $X$.  
By using $(4.5)$, we can explicitly describe $X$ for each Hermann action.  
Set 
$${\it l}_{\beta}:=\left\{
\begin{array}{ll}
\displaystyle{\beta^{-1}(0)} & 
\displaystyle{(\beta\in\Pi_V)}\\
\displaystyle{\beta^{-1}(-\frac{\pi}{2})} & 
\displaystyle{(\beta\in\Pi_H\setminus\Pi_V)}\\
\displaystyle{\delta^{-1}(\pi)} & 
\displaystyle{(\beta=\delta).}
\end{array}\right.$$
Let $a=1$ when $\delta\in{\triangle'}_+^V\setminus{\triangle'}_+^H$ 
and $a=\frac12$ when $\delta\in{\triangle'}_+^H$.  
Fix a stratum $\widetilde{\sigma}$ of $\partial\widetilde C$.  
Define ${\triangle'}^{\widetilde{\sigma}}_+$ by 
${\triangle'}^{\widetilde{\sigma}}_+:
=(\triangle'_+\setminus(\Pi\cup\{a\delta\}))\cup\{\beta\in\Pi\cup\{a\delta\}\,
\vert\,\widetilde{\sigma}\not\subset{\it l}_{\beta}\}$ and set 
$({\triangle'}^V_+)^{\widetilde{\sigma}}:={\triangle'}^V_+\cap
{\triangle'}^{\widetilde{\sigma}}_+$ and 
$({\triangle'}^H_+)^{\widetilde{\sigma}}:={\triangle'}^H_+\cap
{\triangle'}^{\widetilde{\sigma}}_+$.  Fix $Z_0\in\widetilde{\sigma}$.  Set 
$F(Z_0):=H({\rm Exp}Z_0)$ and $\widetilde F(Z_0):=P(G,H\times K)\cdot Z_0$.  
Denote by $\widetilde H^{Z_0}$ (resp. $H^{Z_0}$) the mean curvature vector of 
$\widetilde F(Z_0)$ (resp. $F(Z_0)$).  
Define a tangent vector field 
$X^{\widetilde{\sigma}}$ on $\widetilde{\sigma}$ by 
$(X^{\widetilde{\sigma}})_{Z_0}:=(\widetilde H^{Z_0})_{Z_0}$ 
($Z_0\in\widetilde{\sigma}$).  
From $(3.5)$ (the case of $w=Z_0-Y_0$), 
we have 
$$\begin{array}{l}
\displaystyle{
(X^{\widetilde{\sigma}})_{Z_0}
=-\sum_{\beta\in({\triangle'}^V_+)^{\widetilde{\sigma}}\setminus
({\triangle'}^H_+)^{\widetilde{\sigma}}}m_{\beta}^V
\cot\beta(Z_0)\beta^{\sharp}}\\
\hspace{1.2truecm}\displaystyle{
+\sum_{\beta\in({\triangle'}^H_+)^{\widetilde{\sigma}}\setminus
({\triangle'}^V_+)^{\widetilde{\sigma}}}
m_{\beta}^H\tan\beta(Z_0)\beta^{\sharp}}\\
\hspace{1.2truecm}\displaystyle{
-\sum_{\beta\in({\triangle'}^V_+)^{\widetilde{\sigma}}\cap
({\triangle'}^H_+)^{\widetilde{\sigma}}}
\left((m_{\beta}^V+m_{\beta}^H)\cot2\beta(Z_0)\right.}\\
\hspace{4.5truecm}\displaystyle{\left.+(m_{\beta}^V-m_{\beta}^H){\rm cosec}
2\beta(Z_0)\right)\beta^{\sharp}.}
\end{array}\leqno{(4.6)}$$
In order to analyze the mean curvature flows having singular orbits through 
$\widetilde{\sigma}$ of the $H$-action as initial data, we have only to 
analyze the vector field $X^{\widetilde{\sigma}}$.  
Now we shall describe explicitly the above vector fields $X$ and 
$X^{\widetilde{\sigma}}$ for some Hermann actions.  
We shall use the above notations.  

\vspace{0.5truecm}

\noindent
{\it Example 1.} We consider the isotropy action 
$\displaystyle{Sp(n) \curvearrowright SU(2n)/Sp(n)}$.  
Then, since $\triangle'$ is equal to the root system of 
$SU(2n)/Sp(n)$, it is of $({\rm A}_{n-1})$-type and each of its roots 
is of multiplicity $4$.  
We can describe 
$\mathfrak b,\triangle'_+,\Pi$ and $\delta$ as 
$$\begin{array}{l}
\displaystyle{\mathfrak b=\{\sum\limits_{i=1}^nx_ie_i\,\vert\,
\sum_{i=1}^nx_i=0\}(\subset{\bf R}^n),}\\
\displaystyle{\triangle'_+=\{(\beta_i-\beta_j)\vert_{\mathfrak b}\,\vert\,
1\leq i<j\leq n\},}\\
\displaystyle{\Pi=\{(\beta_i-\beta_{i+1})\vert_{\mathfrak b}\,\vert\,
1\leq i\leq n-1\},\,\,\,\,\delta=(\beta_1-\beta_n)\vert_{\mathfrak b},}
\end{array}\leqno{(4.7)}$$
where $(e_1,\cdots,e_n)$ is the orthonormal base of ${\bf R}^n$ and 
$(\beta_1,\cdots,\beta_n)$ is the dual base of $(e_1,\cdots,e_n)$.  
For simplicity, we set $\beta_{ij}:=(\beta_i-\beta_j)\vert_{\mathfrak b}\,
(1\leq i<j\leq n),\,I:=\{1,\cdots,n\}$ and 
$\widehat I:=\{(i,j)\in I^2\,\vert\,
1\leq i<j\leq n\}$.  
Since we consider the isotropy action, we have 
${\triangle'}^V_+={\triangle'}_+$ and ${\triangle'}^H_+=\emptyset$, that is, 
$$\begin{array}{l}
\displaystyle{\widetilde C=\{{\bf x}\in\mathfrak b\,\vert\,0<\beta_{ij}
({\bf x})<\pi\,\,((i,j)\in\widehat I)\}}\\
\hspace{0.4truecm}\displaystyle{=\{{\bf x}\in\mathfrak b\,\vert\,
\beta_{i,i+1}({\bf x})>0\,(i=1,\cdots,n-1),\,\beta_{1n}({\bf x})<\pi\}.}
\end{array}$$
From $(4.5)$, we can describe $X$ explicitly as 
$$\begin{array}{l}
\displaystyle{X_{\bf x}=-4\sum_{(i,j)\in\widehat I}\cot\beta_{ij}({\bf x})
\beta_{ij}^{\sharp}=-4\sum_{(i,j)\in\widehat I}\cot(x_i-x_j)(e_i-e_j)}\\
\hspace{0.7truecm}\displaystyle{
=-4\sum_{i\in I}(\sum_{j\in I\setminus\{i\}}\cot(x_i-x_j))e_i\qquad\,\,
({\bf x}(=\sum_{i=1}^nx_ie_i)\in\widetilde C).}
\end{array}\leqno{(4.8)}$$
Take a stratum $\widetilde{\sigma}$ of $\partial\widetilde C$.  
Set $\widehat I_{\widetilde{\Pi}}:=\{(i,i+1)\,\vert\,1\leq i\leq n-1\}\cup
\{(1,n)\}$, ${\it l}_{i,i+1}:=\beta_{i,i+1}^{-1}(0)$ ($1\leq i\leq n-1$) 
and ${\it l}_{1n}:=\beta_{1n}^{-1}(\pi)$.  Set 
$\widehat I_{\widetilde{\sigma}}:=\{(i,j)\in
\widehat I_{\widetilde{\Pi}}\,\vert\,
\widetilde{\sigma}\subset{\it l}_{ij}\}$.  
Also, set $I^{\widetilde{\sigma}}_i:=\{j\in I\,\vert\,(i,j)\in
\widehat I\setminus\widehat I_{\widetilde{\sigma}}\,\,{\rm or}\,\,
(j,i)\in\widehat I\setminus\widehat I_{\widetilde{\sigma}}\}$ 
($i=1,\cdots,n$).  
Then, from $(4.6)$, we can describe $X^{\widetilde{\sigma}}$ explicitly as 
$$\begin{array}{l}
\displaystyle{X^{\widetilde{\sigma}}_{\bf x}
=-4\sum_{(i,j)\in\widehat I\setminus\widehat I_{\widetilde{\sigma}}}
\cot(x_i-x_j)(e_i-e_j)}\\
\hspace{0.6truecm}\displaystyle{
=-4\sum_{i\in I}\left(\sum_{j\in I^{\widetilde{\sigma}}_i}\cot(x_i-x_j)\right)
e_i\qquad\,\,\,\,\,\,({\bf x}(=\sum_{i=1}^nx_ie_i)\in\widetilde{\sigma}).}
\end{array}\leqno{(4.9)}$$
According to $(4.8)$, a principal orbit $Sp(n)({\rm Exp}\,{\bf x})$ 
(${\bf x}\in\widetilde C$) is minimal if and only if the following relations 
hold:
$$\sum_{j\in I\setminus\{i\}}\cot(x_i-x_j)=0\,\,\,\,(i=1,\cdots,n). 
\leqno{(4.10)}$$
Also, according to $(4.9)$, a singular orbit $Sp(n)({\rm Exp}\,{\bf x})$ 
(${\bf x}\in\widetilde{\sigma}$) is minimal if and only if the following 
relations hold:
$$\sum_{j\in I^{\widetilde{\sigma}}_i}\cot(x_i-x_j)=0\,\,\,\,(i=1,\cdots,n). 
\leqno{(4.11)}$$
In the sequel, we consider the case of $n=3$.  By solving $(4.10)$ under 
$0<x_i-x_j<\pi$ ($(i,j)\in\widehat I$), we have 
$(x_1,x_2,x_3)=(\frac{\pi}{3},0,-\frac{\pi}{3})$.  
Therefore the orbit $Sp(3)({\rm Exp}(\frac{\pi}{3}e_1-\frac{\pi}{3}e_3))$ 
is the only minimal principal orbit of the $Sp(3)$-action.  
Denote by $\widetilde{\sigma}_{12},\,\widetilde{\sigma}_{23}$ and 
$\widetilde{\sigma}_{13}$ one dimensional stratums of 
$\partial\widetilde C$ which are contained in $\beta_{12}^{-1}(0),\,
\beta_{23}^{-1}(0)$ and $\beta_{13}^{-1}(\pi)$, respectively.  For 
$(k,{\it l})\in\widehat I$, by solving $(4.11)$ 
(the case of $\widetilde{\sigma}=\widetilde{\sigma}_{k{\it l}}$) under 
$0<x_i-x_j<\pi\,\,((i,j)\in\widehat I\setminus
\widehat I_{\widetilde{\sigma}_{k{\it l}}})$, 
we have 
$$(x_1,x_2,x_3)=\left\{
\begin{array}{ll}
\displaystyle{(\frac{\pi}{6},\frac{\pi}{6},-\frac{\pi}{3})} & 
\displaystyle{((k,{\it l})=(1,2))}\\
\displaystyle{(\frac{\pi}{3},-\frac{\pi}{6},-\frac{\pi}{6})} & 
\displaystyle{((k,{\it l})=(2,3))}\\
\displaystyle{(\frac{\pi}{2},0,-\frac{\pi}{2})} & 
\displaystyle{((k,{\it l})=(1,3)).}
\end{array}\right.$$
Therefore the orbits $Sp(3)({\rm Exp}(\frac{\pi}{6}(e_1+e_2-2e_3))),\,
Sp(3)({\rm Exp}(\frac{\pi}{6}(2e_1-e_2-e_3)))$ and 
$Sp(3)({\rm Exp}(\frac{\pi}{2}(e_1-e_3)))$ are the only minimal 
singular orbits through one dimensional stratums of $\partial\widetilde C$ 
of the $Sp(3)$-action.  

\vspace{1truecm}

\centerline{
\unitlength 0.1in
\begin{picture}( 43.2000, 21.5000)(  9.3000,-30.8000)
%
\special{pn 8}%
\special{pa 4000 1000}%
\special{pa 2840 2800}%
\special{fp}%
%
\special{pn 8}%
\special{pa 3600 1000}%
\special{pa 4770 2800}%
\special{fp}%
%
\special{pn 8}%
\special{pa 2590 2600}%
\special{pa 5000 2600}%
\special{fp}%
%
\special{pn 8}%
\special{pa 3360 1310}%
\special{pa 4300 1310}%
\special{fp}%
%
\special{pn 8}%
\special{pa 4820 2310}%
\special{pa 4480 2850}%
\special{fp}%
%
\special{pn 8}%
\special{pa 3800 1310}%
\special{pa 3800 2600}%
\special{dt 0.045}%
%
\special{pn 8}%
\special{pa 4640 2600}%
\special{pa 3390 1970}%
\special{dt 0.045}%
%
\special{pn 8}%
\special{pa 2980 2600}%
\special{pa 4220 1960}%
\special{dt 0.045}%
%
\special{pn 13}%
\special{sh 1}%
\special{ar 3800 2170 10 10 0  6.28318530717959E+0000}%
\special{sh 1}%
\special{ar 3800 2170 10 10 0  6.28318530717959E+0000}%
%
\special{pn 13}%
\special{sh 1}%
\special{ar 3800 2600 10 10 0  6.28318530717959E+0000}%
\special{sh 1}%
\special{ar 3800 2600 10 10 0  6.28318530717959E+0000}%
%
\special{pn 13}%
\special{sh 1}%
\special{ar 3380 1960 10 10 0  6.28318530717959E+0000}%
\special{sh 1}%
\special{ar 3380 1960 10 10 0  6.28318530717959E+0000}%
%
\special{pn 13}%
\special{sh 1}%
\special{ar 4220 1960 10 10 0  6.28318530717959E+0000}%
\special{sh 1}%
\special{ar 4220 1960 10 10 0  6.28318530717959E+0000}%
%
\special{pn 8}%
\special{pa 3690 2860}%
\special{pa 3800 2610}%
\special{dt 0.045}%
\special{sh 1}%
\special{pa 3800 2610}%
\special{pa 3756 2664}%
\special{pa 3780 2660}%
\special{pa 3792 2680}%
\special{pa 3800 2610}%
\special{fp}%
%
\special{pn 8}%
\special{pa 3160 3060}%
\special{pa 3290 2600}%
\special{dt 0.045}%
\special{sh 1}%
\special{pa 3290 2600}%
\special{pa 3254 2660}%
\special{pa 3276 2652}%
\special{pa 3292 2670}%
\special{pa 3290 2600}%
\special{fp}%
%
\special{pn 8}%
\special{pa 3230 1710}%
\special{pa 3370 1940}%
\special{dt 0.045}%
\special{sh 1}%
\special{pa 3370 1940}%
\special{pa 3352 1874}%
\special{pa 3342 1894}%
\special{pa 3318 1894}%
\special{pa 3370 1940}%
\special{fp}%
%
\special{pn 8}%
\special{pa 2890 2090}%
\special{pa 3180 2240}%
\special{dt 0.045}%
\special{sh 1}%
\special{pa 3180 2240}%
\special{pa 3130 2192}%
\special{pa 3134 2216}%
\special{pa 3112 2228}%
\special{pa 3180 2240}%
\special{fp}%
%
\special{pn 8}%
\special{pa 4500 1960}%
\special{pa 4230 1950}%
\special{dt 0.045}%
\special{sh 1}%
\special{pa 4230 1950}%
\special{pa 4296 1972}%
\special{pa 4284 1952}%
\special{pa 4298 1932}%
\special{pa 4230 1950}%
\special{fp}%
%
\special{pn 8}%
\special{pa 4260 1660}%
\special{pa 3810 2170}%
\special{dt 0.045}%
\special{sh 1}%
\special{pa 3810 2170}%
\special{pa 3870 2134}%
\special{pa 3846 2130}%
\special{pa 3840 2108}%
\special{pa 3810 2170}%
\special{fp}%
%
\special{pn 8}%
\special{pa 5120 2440}%
\special{pa 4770 2400}%
\special{dt 0.045}%
\special{sh 1}%
\special{pa 4770 2400}%
\special{pa 4834 2428}%
\special{pa 4824 2406}%
\special{pa 4840 2388}%
\special{pa 4770 2400}%
\special{fp}%
%
\special{pn 8}%
\special{pa 4290 1130}%
\special{pa 4150 1310}%
\special{dt 0.045}%
\special{sh 1}%
\special{pa 4150 1310}%
\special{pa 4208 1270}%
\special{pa 4184 1268}%
\special{pa 4176 1246}%
\special{pa 4150 1310}%
\special{fp}%
\put(44.6000,-29.0000){\makebox(0,0)[rt]{$\frac{\pi}{6}(e_1+e_2-2e_3)$}}%
\put(33.3000,-30.8000){\makebox(0,0)[rt]{$\beta_{12}^{-1}(0)$}}%
\put(51.4000,-18.8000){\makebox(0,0)[rt]{$\frac{\pi}{2}(e_1-e_3)$}}%
%
\put(46.6000,-20.3000){\makebox(0,0)[lb]{}}%
\put(48.4000,-14.5000){\makebox(0,0)[rt]{$\frac{\pi}{3}(e_1-e_3)$}}%
\put(35.4000,-14.9000){\makebox(0,0)[rt]{$\frac{\pi}{6}(2e_1-e_2-e_3)$}}%
\put(28.4000,-19.9000){\makebox(0,0)[rt]{$\beta_{23}^{-1}(0)$}}%
\put(56.1000,-23.8000){\makebox(0,0)[rt]{$\beta_{23}^{-1}(\pi)$}}%
\put(44.9000,-9.3000){\makebox(0,0)[rt]{$\beta_{12}^{-1}(\pi)$}}%
%
\special{pn 8}%
\special{pa 5250 2110}%
\special{pa 4410 2220}%
\special{dt 0.045}%
\special{sh 1}%
\special{pa 4410 2220}%
\special{pa 4480 2232}%
\special{pa 4464 2214}%
\special{pa 4474 2192}%
\special{pa 4410 2220}%
\special{fp}%
\put(57.6000,-20.0000){\makebox(0,0)[rt]{$\beta_{13}^{-1}(\pi)$}}%
%
\special{pn 13}%
\special{sh 1}%
\special{ar 2980 2600 10 10 0  6.28318530717959E+0000}%
\special{sh 1}%
\special{ar 2980 2600 10 10 0  6.28318530717959E+0000}%
\put(30.4000,-26.4000){\makebox(0,0)[rt]{${\bf 0}$}}%
\end{picture}%
\hspace{5truecm}}

\vspace{0.5truecm}

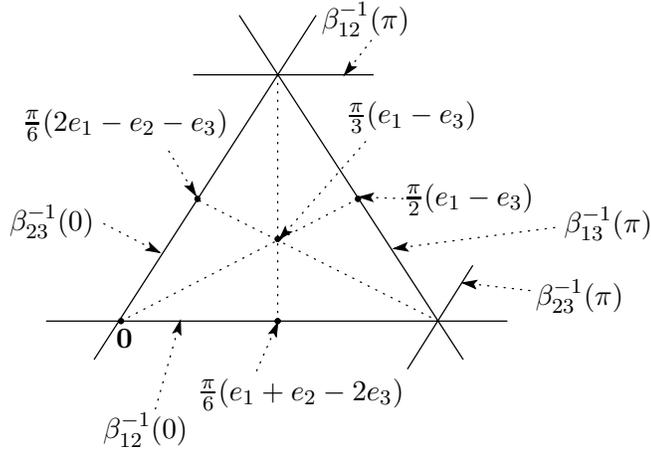
\centerline{{\bf Fig. 5}}

\vspace{0.5truecm}

\noindent
Next we shall investigate the divergences ${\rm div}\,X$ and 
${\rm div}\,X^{\widetilde{\sigma}_{k{\it l}}}$ of $X$ and 
$X^{\widetilde{\sigma}_{k{\it l}}}$ ($(k,{\it l})\in\widehat I$).  
Take an orthonormal base $\{v_1:=\frac{1}{\sqrt2}(e_1-e_2),\,
v_2:=\frac{1}{\sqrt6}(e_1+e_2-2e_3)\}$ of $\mathfrak b$.  Let ${\bf x}=\sum\limits_{i=1}^3x_ie_i=\sum\limits_{i=1}^2
y_iv_i$.  Then, from $(4.8)$, we have 
$$\begin{array}{l}
\displaystyle{X_{\bf x}=-4\left(
\sqrt2\cot\sqrt2y_1+(\sqrt2-\frac{1}{\sqrt2})
\cot(\frac{1}{\sqrt2}y_1+\frac{\sqrt6}{2}y_2)\right.}\\
\hspace{2truecm}\displaystyle{\left.-\frac{1}{\sqrt2}\cot
(-\frac{1}{\sqrt2}y_1+\frac{\sqrt6}{2}y_2)\right)v_1}\\
\hspace{1.4truecm}\displaystyle{-2\sqrt6\left(
\cot(\frac{1}{\sqrt2}y_1+\frac{\sqrt6}{2}y_2)+
\cot(-\frac{1}{\sqrt2}y_1+\frac{\sqrt6}{2}y_2)\right)v_2}
\end{array}$$
and hence 
$$({\rm div}\,X)_{\bf x}=8\left(\frac{1}{\sin^2\sqrt2y_1}+
\frac{1}{\sin^2(\frac{1}{\sqrt2}y_1+\frac{\sqrt6}{2}y_2)}+
\frac{1}{\sin^2(-\frac{1}{\sqrt2}y_1+\frac{\sqrt6}{2}y_2)}\right)\,>\,0.$$
Therefore, since $X$ has the only zero point, ${\rm div}\,X>0$ on 
$\widetilde C$ and $X$ is as in Fig.1 over a collar neighborhood of 
$\widetilde{\sigma}_{k{\it l}}\,((k,{\it l})\in\widehat I)$ 
(see the proof of Theorem A), all the integral curves 
of $X$ through points other than its zero point $\frac{\pi}{3}(e_1-e_3)$ 
converge to points of $\partial\widetilde C$ in finite time.  
Similarly we can show ${\rm div}\,X^{\widetilde{\sigma}_{k{\it l}}}\,>\,0$ on 
$\widetilde{\sigma}_{k{\it l}}$ ($(k,{\it l})\in\widehat I$).  Hence all 
integral curves of $X^{\widetilde{\sigma}_{k{\it l}}}$ through points other 
than its zero point converge to points of 
$\partial\widetilde{\sigma}_{k{\it l}}$ in finite time.  
Therefore we obtain the following fact.  

\vspace{0.5truecm}

\noindent
{\bf Proposition 4.1.} {\sl The mean curvature flow having any non-minimal 
principal orbit of the $Sp(3)$-action (on $SU(6)/Sp(3)$) as initial data 
converges to some singular orbit in finite time.  Also, the mean curvature 
flow having any non-minimal singular orbit through 
$\exp^{\perp}(\widetilde{\sigma}_{k{\it l}})$ of the 
$Sp(3)$-action as initial data converges to one of two singular orbits through 
$\exp^{\perp}(\partial\widetilde{\sigma}_{k{\it l}})$ in finite time.}

\vspace{0.5truecm}

\noindent
{\it Example 2.} 
We consider the Hermann action 
$\displaystyle{SO(2n)\curvearrowright SU(2n)/Sp(n)}$.  
Since $L/H\cap K$ is equal to $(SU(n)\times SU(n))/SU(n)$, the cohomogeneity 
of this action is equal to the rank $n-1$ of $(SU(n)\times SU(n))/SU(n)$.  
On the other hand, the rank of $SU(2n)/Sp(n)$ also is equal to $n-1$.  Thus 
the cohomogeneity of this action is equal to the rank of $SU(2n)/Sp(n)$.  
Hence a maximal abelian subspace $\mathfrak b$ of 
$\mathfrak p\cap\mathfrak q$ is also a maximal abelian subspace of 
$\mathfrak p$ and hence $\triangle'$ is the root system of $SU(2n)/Sp(n)$, 
it is of $(A_{n-1})$-type and each 
of its roots is of multiplicity $4$.  Hence the quantities 
$\mathfrak b,\,\triangle'_+,\,\Pi$ and $\delta$ are given as in $(4.7)$.  Let 
$\beta_i,\,\beta_{ij},\,I$ and $\widehat I$ be as in Example 1.  Also, since 
${\triangle'}^V$ is the root system of $(SU(n)\times SU(n))/SU(n)$, it is of 
$(A_{n-1})$-type and each of its roots is of multiplicity $2$.  Hence we have 
${\triangle'}^H_+=\triangle'_+$ and each of its roots is of multiplicity $2$.  
Therefore we have 
$$\widetilde C=\{{\bf x}\in\mathfrak b\,\vert\,
\beta_{i,i+1}({\bf x})>0\,\,(i=1,\cdots,n-1),\,\,\beta_{1n}({\bf x})
<\frac{\pi}{2}\}$$
and 
$$\begin{array}{l}
\displaystyle{X_{\bf x}=-4\sum_{(i,j)\in\widehat I}\cot2\beta_{ij}({\bf x})
\beta_{ij}^{\sharp}=-4\sum_{(i,j)\in\widehat I}\cot2(x_i-x_j)(e_i-e_j)}\\
\hspace{0.7truecm}\displaystyle{
=-4\sum_{i\in I}(\sum_{j\in I\setminus\{i\}}\cot2(x_i-x_j))e_i.}
\end{array}\leqno{(4.12)}$$
Take a stratum $\widetilde{\sigma}$ of $\partial\widetilde C$.  
Let $\widehat I_{\widetilde{\sigma}}$ and $I^{\widetilde{\sigma}}_i$ be 
as in Example 1, where ${\it l}_{1n}:=\beta_{1n}^{-1}(\frac{\pi}{2})$.  
Then, from $(4.6)$, we can describe $X^{\widetilde{\sigma}}$ explicitly as 
$$X^{\widetilde{\sigma}}_{\bf x}
=-4\sum_{i\in I}\left(\sum_{j\in I^{\widetilde{\sigma}}_i}\cot2(x_i-x_j)\right)
e_i\qquad\,\,\,\,\,\,({\bf x}(=\sum_{i=1}^nx_ie_i)\in\widetilde{\sigma}).
\leqno{(4.13)}$$
According to $(4.12)$, a principal orbit $SO(2n)({\rm Exp}\,{\bf x})$ 
(${\bf x}\in\widetilde C$) is minimal if and only if the following relations 
hold:
$$\sum_{j\in I\setminus\{i\}}\cot2(x_i-x_j)=0\,\,\,\,(i=1,\cdots,n). 
\leqno{(4.14)}$$
Also, according to $(4.13)$, a singular orbit $Sp(n)({\rm Exp}\,{\bf x})$ 
(${\bf x}\in\widetilde{\sigma}$) is minimal if and only if the following 
relations hold:
$$\sum_{j\in I^{\widetilde{\sigma}}_i}\cot2(x_i-x_j)=0\,\,\,\,(i=1,\cdots,n). 
\leqno{(4.15)}$$
In the sequel, we consider the case of $n=3$.  By solving $(4.14)$ under 
$0<x_i-x_j<\frac{\pi}{2}$ ($(i,j)\in\widehat I$), we have 
$(x_1,x_2,x_3)=(\frac{\pi}{6},0,-\frac{\pi}{6})$.  
Therefore the orbit 
$SO(6)({\rm Exp}(\frac{\pi}{6}e_1-\frac{\pi}{6}e_3))$ 
is the only minimal principal orbit of the $SO(6)$-action.  
Denote by $\widetilde{\sigma}_{12},\,\widetilde{\sigma}_{23}$ and 
$\widetilde{\sigma}_{13}$ one dimensional stratums of 
$\partial\widetilde C$ which are contained in $\beta_{12}^{-1}(0),\,
\beta_{23}^{-1}(0)$ and $\beta_{13}^{-1}(\frac{\pi}{2})$, respectively.  
For $(k,{\it l})\in\widehat I$, by solving $(4.15)$ 
(the case of $\widetilde{\sigma}=\widetilde{\sigma}_{k{\it l}}$) under 
$0<x_i-x_j<\frac{\pi}{2}\,\,
((i,j)\in\widehat I\setminus\widehat I_{\widetilde{\sigma}_{k{\it l}}})$, 
we have 
$$(x_1,x_2,x_3)=\left\{
\begin{array}{ll}
\displaystyle{(\frac{\pi}{12},\frac{\pi}{12},-\frac{\pi}{6})} & 
\displaystyle{((k,{\it l})=(1,2))}\\
\displaystyle{(\frac{\pi}{6},-\frac{\pi}{12},-\frac{\pi}{12})} & 
\displaystyle{((k,{\it l})=(2,3))}\\
\displaystyle{(\frac{\pi}{4},0,-\frac{\pi}{4})} & 
\displaystyle{((k,{\it l})=(1,3)).}
\end{array}\right.$$
Therefore the orbits $SO(6)({\rm Exp}(\frac{\pi}{12}(e_1+e_2-2e_3))),\,
SO(6)({\rm Exp}(\frac{\pi}{12}(2e_1-e_2-e_3)))$ and 
$SO(6)({\rm Exp}(\frac{\pi}{4}(e_1-e_3)))$ are the only minimal 
singular orbits through one dimensional stratums of $\partial\widetilde C$ of 
the $SO(6)$-action.  

\vspace{1truecm}

\centerline{
\unitlength 0.1in
\begin{picture}( 44.1000, 21.4000)(  8.4000,-30.8000)
%
\special{pn 8}%
\special{pa 4000 1000}%
\special{pa 2840 2800}%
\special{fp}%
%
\special{pn 8}%
\special{pa 3600 1000}%
\special{pa 4770 2800}%
\special{fp}%
%
\special{pn 8}%
\special{pa 2590 2600}%
\special{pa 5000 2600}%
\special{fp}%
%
\special{pn 8}%
\special{pa 3360 1310}%
\special{pa 4300 1310}%
\special{fp}%
%
\special{pn 8}%
\special{pa 4820 2310}%
\special{pa 4480 2850}%
\special{fp}%
%
\special{pn 8}%
\special{pa 3800 1310}%
\special{pa 3800 2600}%
\special{dt 0.045}%
%
\special{pn 8}%
\special{pa 4640 2600}%
\special{pa 3390 1970}%
\special{dt 0.045}%
%
\special{pn 8}%
\special{pa 2980 2600}%
\special{pa 4220 1960}%
\special{dt 0.045}%
%
\special{pn 13}%
\special{sh 1}%
\special{ar 3800 2170 10 10 0  6.28318530717959E+0000}%
\special{sh 1}%
\special{ar 3800 2170 10 10 0  6.28318530717959E+0000}%
%
\special{pn 13}%
\special{sh 1}%
\special{ar 3800 2600 10 10 0  6.28318530717959E+0000}%
\special{sh 1}%
\special{ar 3800 2600 10 10 0  6.28318530717959E+0000}%
%
\special{pn 13}%
\special{sh 1}%
\special{ar 3380 1960 10 10 0  6.28318530717959E+0000}%
\special{sh 1}%
\special{ar 3380 1960 10 10 0  6.28318530717959E+0000}%
%
\special{pn 13}%
\special{sh 1}%
\special{ar 4220 1960 10 10 0  6.28318530717959E+0000}%
\special{sh 1}%
\special{ar 4220 1960 10 10 0  6.28318530717959E+0000}%
%
\special{pn 8}%
\special{pa 3690 2860}%
\special{pa 3800 2610}%
\special{dt 0.045}%
\special{sh 1}%
\special{pa 3800 2610}%
\special{pa 3756 2664}%
\special{pa 3780 2660}%
\special{pa 3792 2680}%
\special{pa 3800 2610}%
\special{fp}%
%
\special{pn 8}%
\special{pa 3160 3060}%
\special{pa 3290 2600}%
\special{dt 0.045}%
\special{sh 1}%
\special{pa 3290 2600}%
\special{pa 3254 2660}%
\special{pa 3276 2652}%
\special{pa 3292 2670}%
\special{pa 3290 2600}%
\special{fp}%
%
\special{pn 8}%
\special{pa 3230 1710}%
\special{pa 3370 1940}%
\special{dt 0.045}%
\special{sh 1}%
\special{pa 3370 1940}%
\special{pa 3352 1874}%
\special{pa 3342 1894}%
\special{pa 3318 1894}%
\special{pa 3370 1940}%
\special{fp}%
%
\special{pn 8}%
\special{pa 2890 2090}%
\special{pa 3180 2240}%
\special{dt 0.045}%
\special{sh 1}%
\special{pa 3180 2240}%
\special{pa 3130 2192}%
\special{pa 3134 2216}%
\special{pa 3112 2228}%
\special{pa 3180 2240}%
\special{fp}%
%
\special{pn 8}%
\special{pa 4500 1960}%
\special{pa 4230 1950}%
\special{dt 0.045}%
\special{sh 1}%
\special{pa 4230 1950}%
\special{pa 4296 1972}%
\special{pa 4284 1952}%
\special{pa 4298 1932}%
\special{pa 4230 1950}%
\special{fp}%
%
\special{pn 8}%
\special{pa 4260 1660}%
\special{pa 3810 2170}%
\special{dt 0.045}%
\special{sh 1}%
\special{pa 3810 2170}%
\special{pa 3870 2134}%
\special{pa 3846 2130}%
\special{pa 3840 2108}%
\special{pa 3810 2170}%
\special{fp}%
%
\special{pn 8}%
\special{pa 5120 2440}%
\special{pa 4770 2400}%
\special{dt 0.045}%
\special{sh 1}%
\special{pa 4770 2400}%
\special{pa 4834 2428}%
\special{pa 4824 2406}%
\special{pa 4840 2388}%
\special{pa 4770 2400}%
\special{fp}%
%
\special{pn 8}%
\special{pa 4290 1130}%
\special{pa 4150 1310}%
\special{dt 0.045}%
\special{sh 1}%
\special{pa 4150 1310}%
\special{pa 4208 1270}%
\special{pa 4184 1268}%
\special{pa 4176 1246}%
\special{pa 4150 1310}%
\special{fp}%
\put(44.7000,-28.9000){\makebox(0,0)[rt]{$\frac{\pi}{12}(e_1+e_2-2e_3)$}}%
\put(33.4000,-30.8000){\makebox(0,0)[rt]{$\beta_{12}^{-1}(0)$}}%
\put(51.5000,-18.7000){\makebox(0,0)[rt]{$\frac{\pi}{4}(e_1-e_3)$}}%
%
\put(46.6000,-20.3000){\makebox(0,0)[lb]{}}%
\put(47.4000,-14.6000){\makebox(0,0)[rt]{$\frac{\pi}{6}(e_1-e_3)$}}%
\put(35.4000,-14.9000){\makebox(0,0)[rt]{$\frac{\pi}{12}(2e_1-e_2-e_3)$}}%
\put(28.4000,-19.9000){\makebox(0,0)[rt]{$\beta_{23}^{-1}(0)$}}%
\put(56.9000,-23.9000){\makebox(0,0)[rt]{$\beta_{23}^{-1}(\frac{\pi}{2})$}}%
\put(45.4000,-9.4000){\makebox(0,0)[rt]{$\beta_{12}^{-1}(\frac{\pi}{2})$}}%
%
\special{pn 8}%
\special{pa 5250 2110}%
\special{pa 4410 2220}%
\special{dt 0.045}%
\special{sh 1}%
\special{pa 4410 2220}%
\special{pa 4480 2232}%
\special{pa 4464 2214}%
\special{pa 4474 2192}%
\special{pa 4410 2220}%
\special{fp}%
\put(57.9000,-20.1000){\makebox(0,0)[rt]{$\beta_{13}^{-1}(\frac{\pi}{2})$}}%
%
\special{pn 13}%
\special{sh 1}%
\special{ar 2970 2600 10 10 0  6.28318530717959E+0000}%
\special{sh 1}%
\special{ar 2970 2600 10 10 0  6.28318530717959E+0000}%
\put(30.5000,-26.3000){\makebox(0,0)[rt]{${\bf 0}$}}%
\end{picture}%
\hspace{5truecm}}

\vspace{0.5truecm}

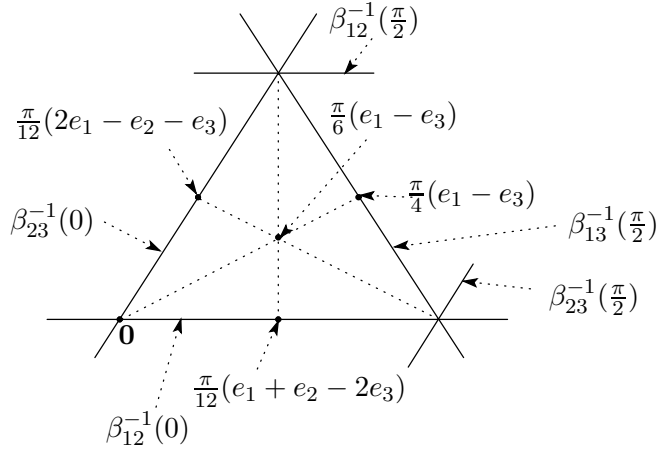
\centerline{{\bf Fig. 6}}

\vspace{0.5truecm}

\noindent
Next we shall investigate the divergences ${\rm div}\,X$ and 
${\rm div}\,X^{\widetilde{\sigma}_{k{\it l}}}$ ($(k,{\it l})\in\widehat I$).  
Take an orthonormal base $\{v_1,\,v_2\}$ of $\mathfrak b$ as in Example 1.  
Let ${\bf x}=\sum\limits_{i=1}^3x_ie_i=\sum\limits_{i=1}^2y_iv_i$.  Then, 
from $(4.12)$, we have 
$$\begin{array}{l}
\displaystyle{X_{\bf x}=-4\left(
\sqrt2\cot2\sqrt2y_1+(\sqrt2-\frac{1}{\sqrt2})
\cot(\sqrt2y_1+\sqrt6y_2)\right.}\\
\hspace{2truecm}\displaystyle{\left.-\frac{1}{\sqrt2}\cot
(-\sqrt2y_1+\sqrt6y_2)\right)v_1}\\
\hspace{1.4truecm}\displaystyle{-2\sqrt6\left(
\cot(\sqrt2y_1+\sqrt6y_2)+
\cot(-\sqrt2y_1+\sqrt6y_2)\right)v_2}
\end{array}$$
and hence 
$$({\rm div}\,X)_{\bf x}=16\left(\frac{1}{\sin^22\sqrt2y_1}+
\frac{1}{\sin^2(\sqrt2y_1+\sqrt6y_2)}+
\frac{1}{\sin^2(-\sqrt2y_1+\sqrt6y_2)}\right)\,>\,0.$$
Therefore, since $X$ has the only zero point and $X$ is as in Fig.1 over 
a collar neighborhood of $\widetilde{\sigma}_{k{\it l}}\,
((k,{\it l})\in\widehat I)$ by the proof of Theorem A, all the integral curves 
of $X$ through points other than its zero point $\frac{\pi}{6}(e_1-e_3)$ 
converge to points of $\partial\widetilde C$ in finite time.  
Similarly we can show ${\rm div}\,X^{\widetilde{\sigma}_{k{\it l}}}\,>\,0$ on 
$\widetilde{\sigma}_{k{\it l}}$ ($(k,{\it l})\in\widehat I$).  Hence all 
integral curves of $X^{\widetilde{\sigma}_{k{\it l}}}$ through points other 
than its zero point converge to points of 
$\partial\widetilde{\sigma}_{k{\it l}}$ in finite time.  
Therefore we obtain the following fact.  

\vspace{0.5truecm}

\noindent
{\bf Proposition 4.2.} {\sl The mean curvature flow having any non-minimal 
principal orbit of the $SO(6)$-action (on $SU(6)/Sp(3)$) as initial data 
converges to some singular orbit in finite time.  Also, the mean curvature 
flow having any non-minimal singular orbit through 
$\exp^{\perp}(\widetilde{\sigma}_{k{\it l}})$ of the 
$SO(6)$-action as initial data converges to one of two singular orbits through 
$\exp^{\perp}(\partial\widetilde{\sigma}_{k{\it l}})$ in finite time.}

\vspace{0.5truecm}

\noindent
{\it Example 3.} 
We consider the isotropy action $\displaystyle{
S(U(p)\times U(q))\curvearrowright SU(p+q)/S(U(p)\times U(q))}$ ($p<q$).  
Then, since $\triangle'$ is the root system of 
$SU(p+q)/S(U(p)\times U(q))$, it is of $({\rm B}_p)$-type and 
$\mathfrak b,\triangle'_+,\Pi$ and $\delta$ are described as 
$$\begin{array}{l}
\displaystyle{\mathfrak b={\rm Span}\{e_1,\cdots,e_p\},}\\
\displaystyle{\triangle'_+=\{\beta_i-\beta_j\,\vert\,
1\leq i<j\leq p\}\cup\{\beta_i\,\vert\,1\leq i\leq p\}}\\
\hspace{1.5truecm}
\displaystyle{\cup\{\beta_i+\beta_j\,\vert\,1\leq i<j\leq p\}\cup
\{2\beta_i\,\vert\,1\leq i\leq p\},}\\
\displaystyle{\Pi=\{\beta_i-\beta_{i+1}\,\vert\,1\leq i\leq p-1\}
\cup\{\beta_p\},\,\,\delta=\beta_1+\beta_2,}
\end{array}
\leqno{(4.16)}$$
where $(e_1,\cdots,e_p)$ is the orthonormal base of $\mathfrak b$ and 
$(\beta_1,\cdots,\beta_p)$ is the dual base of $(e_1,\cdots,e_p)$.  
For simplicity, we set $\beta_{ij}:=\beta_i-\beta_j\,\,(1\leq i<j\leq p),\,
I:=\{1,\cdots,p\}$ and 
$\widehat I:=\{(i,j)\in I^2\,\vert\,1\leq i<j\leq p\}$.  
It is shown that $\beta_{ij}$'s ($1\leq i<j\leq p$) and 
$\beta_i+\beta_j$'s ($1\leq i<j\leq p$) are of multiplicity $2$, 
$\beta_i$'s ($1\leq i\leq p$) are of multiplicity $2(q-p)$ and 
$2\beta_i$'s ($1\leq i\leq p$) are of multiplicity $1$.  
Since we consider the isotropy action, we have 
${\triangle'}^V_+={\triangle'}_+$ and ${\triangle'}^H_+=\emptyset$, that is, 
$$\begin{array}{l}
\displaystyle{\widetilde C=\{{\bf x}\in\mathfrak b\,\vert\,0<\beta_{ij}
({\bf x})<\pi\,\,((i,j)\in\widehat I),\,\,
0<(\beta_i+\beta_j)({\bf x})<\pi\,\,((i,j)\in\widehat I),}\\
\hspace{8.2truecm}\displaystyle{0<2\beta_i({\bf x})<\pi\,\,(i\in I)\}}\\
\hspace{0.4truecm}\displaystyle{=\{{\bf x}\in\mathfrak b\,\vert\,
\beta_{i,i+1}({\bf x})>0\,(1\leq i\leq p-1),\,\beta_p({\bf x})>0,\,\,
(\beta_1+\beta_2)({\bf x})<\pi\}.}
\end{array}$$
From $(4.5)$, we can describe $X$ explicitly as 
$$\begin{array}{l}
\displaystyle{X_{\bf x}=-\sum_{(i,j)\in\widehat I}
2\cot(x_i-x_j)(e_i-e_j)-\sum_{(i,j)\in\widehat I}
2\cot(x_i+x_j)(e_i+e_j)}\\
\hspace{1.2truecm}\displaystyle{-\sum_{i\in I}2(q-p)\cot x_i\cdot e_i
-\sum_{i\in I}\cot2x_i\cdot2e_i}\\
\hspace{0.7truecm}\displaystyle{
=-2\sum_{i\in I}\left(
\sum_{j\in I\setminus\{i\}}(\cot(x_i-x_j))+\cot(x_i+x_j))
+(q-p)\cot x_i+\cot2x_i\right)e_i}\\
\hspace{6truecm}\displaystyle{({\bf x}(=\sum_{i=1}^nx_ie_i)\in\widetilde C).}
\end{array}$$
According to this relation, a principal orbit 
$S(U(p)\times U(q))({\rm Exp}\,{\bf x})$ (${\bf x}\in\widetilde C$) is minimal 
if and only if the following relations hold:
$$\sum_{j\in I\setminus\{i\}}(\cot(x_i-x_j)+\cot(x_i+x_j))+(q-p)\cot x_i
+\cot2x_i=0\,\,\,\,(i=1,\cdots,p).
\leqno{(4.17)}$$
In the sequel, we consider the case of $p=2$.  Then $(4.17)$ is as follows:
$$\left\{
\begin{array}{l}
\displaystyle{\cot(x_1-x_2)+\cot(x_1+x_2)+(q-2)\cot x_1+\cot2x_1=0}\\
\displaystyle{\cot(x_2-x_1)+\cot(x_1+x_2)+(q-2)\cot x_2+\cot2x_2=0.}
\end{array}\right.\leqno{(4.18)}$$
Under $0<x_i-x_j<\pi$ ($(i,j)\in\widehat I$), $0<x_i+x_j<\pi$ 
($(i,j)\in\widehat I$) and $0<x_i<\frac{\pi}{2}$ ($i\in I$), 
the equation $(4.18)$ has the only solution, which we denote by 
$(\alpha_1,\alpha_2)$.  Therefore the orbit $S(U(2)\times U(q))
({\rm Exp}(\alpha_1e_1+\alpha_2e_2))$ is the only minimal principal orbit 
of the $S(U(2)\times U(q))$-action.  
Here we note that $\lim\limits_{q\to\infty}(\alpha_1,\alpha_2)
=(\frac{\pi}{2},\frac{\pi}{2})$.  
Denote by $\widetilde{\sigma}_1,\,\widetilde{\sigma}_2$ and 
$\widetilde{\sigma}_3$ one dimensional stratums of 
$\partial\widetilde C$ which are contained in $\beta_{12}^{-1}(0),\,
\beta_2^{-1}(0)$ and $(\beta_1+\beta_2)^{-1}(\pi)$, respectively.  
Then we have 
$$\begin{array}{l}
\displaystyle{X^{\widetilde{\sigma}_1}_{\bf x}=-2(2\cot2x_1+(q-2)\cot x_1)
(e_1+e_2),}\\
\displaystyle{X^{\widetilde{\sigma}_2}_{\bf x}=-2(q\cot x_1+\cot2x_1)e_1,}\\
\displaystyle{X^{\widetilde{\sigma}_3}_{\bf x}=-2(2\cot2x_1+(q-2)\cot x_1)
(e_1-e_2).}
\end{array}$$
Hence the orbit $S(U(2)\times U(q))({\rm Exp}\,{\bf x})$ 
(${\bf x}\in \widetilde{\sigma}_1$) is minimal if and only if 
$2\cot2x_1+(q-2)\cot x_1=0$ ($x_2=x_1$) holds, that is, 
$(x_1,x_2)=(\arctan\sqrt{q-1},\arctan\sqrt{q-1})$.  
Hence the orbit $S(U(2)\times U(q))
({\rm Exp}(\arctan\sqrt{q-1}(e_1+e_2)))$ is the only minimal singular orbit 
through $\widetilde{\sigma}_1$.  
Also, the orbit $S(U(2)\times U(q))({\rm Exp}\,{\bf x})$ 
(${\bf x}\in \widetilde{\sigma}_2$) is minimal if and only if 
$q\cot x_1+\cot2x_1=0$ ($x_2=0$) holds, that is, 
$(x_1,x_2)=(\arctan\sqrt{2q+1},0)$.  Hence the orbit $S(U(2)\times U(q))
({\rm Exp}(\arctan\sqrt{2q+1}e_1))$ is the only minimal singular orbit 
through $\widetilde{\sigma}_2$.  
Also, the orbit $S(U(2)\times U(q))({\rm Exp}\,{\bf x})$ 
(${\bf x}\in \widetilde{\sigma}_3$) is minimal if and only if 
$2\cot2x_1+(q-2)\cot x_1=0$ ($x_2=\pi-x_1$) holds, that is, 
$(x_1,x_2)=(\arctan\sqrt{q-1},\pi-\arctan\sqrt{q-1})$.  
Hence the orbit $S(U(2)\times U(q))
({\rm Exp}(\arctan\sqrt{q-1}e_1+(\pi-\arctan\sqrt{q-1})e_2))$ 
is the only minimal singular orbit through $\widetilde{\sigma}_3$.  

\vspace{1truecm}

\centerline{
\unitlength 0.1in
\begin{picture}( 44.0000, 22.3000)(  6.0000,-30.0000)
%
\special{pn 8}%
\special{pa 2590 2600}%
\special{pa 5000 2600}%
\special{fp}%
%
\special{pn 8}%
\special{pa 3910 1170}%
\special{pa 4850 1170}%
\special{fp}%
%
\special{pn 13}%
\special{sh 1}%
\special{ar 4150 2330 10 10 0  6.28318530717959E+0000}%
\special{sh 1}%
\special{ar 4150 2330 10 10 0  6.28318530717959E+0000}%
%
\special{pn 13}%
\special{sh 1}%
\special{ar 3820 2600 10 10 0  6.28318530717959E+0000}%
\special{sh 1}%
\special{ar 3820 2600 10 10 0  6.28318530717959E+0000}%
%
\special{pn 13}%
\special{sh 1}%
\special{ar 3600 1950 10 10 0  6.28318530717959E+0000}%
\special{sh 1}%
\special{ar 3600 1950 10 10 0  6.28318530717959E+0000}%
%
\special{pn 8}%
\special{pa 3650 2830}%
\special{pa 3812 2610}%
\special{dt 0.045}%
\special{sh 1}%
\special{pa 3812 2610}%
\special{pa 3756 2652}%
\special{pa 3780 2652}%
\special{pa 3788 2676}%
\special{pa 3812 2610}%
\special{fp}%
%
\special{pn 8}%
\special{pa 2860 2930}%
\special{pa 3212 2606}%
\special{dt 0.045}%
\special{sh 1}%
\special{pa 3212 2606}%
\special{pa 3148 2638}%
\special{pa 3172 2642}%
\special{pa 3176 2666}%
\special{pa 3212 2606}%
\special{fp}%
%
\special{pn 8}%
\special{pa 3330 1750}%
\special{pa 3586 1954}%
\special{dt 0.045}%
\special{sh 1}%
\special{pa 3586 1954}%
\special{pa 3546 1898}%
\special{pa 3544 1922}%
\special{pa 3520 1928}%
\special{pa 3586 1954}%
\special{fp}%
%
\special{pn 8}%
\special{pa 4640 1460}%
\special{pa 4380 1580}%
\special{dt 0.045}%
\special{sh 1}%
\special{pa 4380 1580}%
\special{pa 4450 1570}%
\special{pa 4428 1558}%
\special{pa 4432 1534}%
\special{pa 4380 1580}%
\special{fp}%
%
\special{pn 8}%
\special{pa 4520 2960}%
\special{pa 4228 2764}%
\special{dt 0.045}%
\special{sh 1}%
\special{pa 4228 2764}%
\special{pa 4272 2818}%
\special{pa 4272 2794}%
\special{pa 4294 2784}%
\special{pa 4228 2764}%
\special{fp}%
%
\special{pn 8}%
\special{pa 4890 990}%
\special{pa 4750 1170}%
\special{dt 0.045}%
\special{sh 1}%
\special{pa 4750 1170}%
\special{pa 4808 1130}%
\special{pa 4784 1128}%
\special{pa 4776 1106}%
\special{pa 4750 1170}%
\special{fp}%
\put(42.9000,-28.9000){\makebox(0,0)[rt]{\scriptsize{$\arctan\sqrt{q-1}(e_1+e_2)$}}}%
\put(31.0000,-29.6000){\makebox(0,0)[rt]{$\beta_{12}^{-1}(0)$}}%
\put(65.8000,-18.7000){\makebox(0,0)[rt]{\scriptsize{$\arctan\sqrt{q-1}e_1+(\pi-\arctan\sqrt{q-1})e_2$}}}%
%
\put(46.6000,-20.3000){\makebox(0,0)[lb]{}}%
\put(29.1000,-21.0000){\makebox(0,0)[rt]{$\beta_2^{-1}(0)$}}%
\put(52.1000,-7.7000){\makebox(0,0)[rt]{$\beta_{12}^{-1}(\pi)$}}%
%
\special{pn 13}%
\special{sh 1}%
\special{ar 2970 2600 10 10 0  6.28318530717959E+0000}%
\special{sh 1}%
\special{ar 2970 2600 10 10 0  6.28318530717959E+0000}%
\put(30.5000,-26.3000){\makebox(0,0)[rt]{${\bf 0}$}}%
%
\special{pn 8}%
\special{pa 4380 2810}%
\special{pa 4380 910}%
\special{fp}%
%
\special{pn 8}%
\special{pa 2790 2780}%
\special{pa 4650 900}%
\special{fp}%
%
\special{pn 13}%
\special{sh 1}%
\special{ar 3370 2190 10 10 0  6.28318530717959E+0000}%
\special{sh 1}%
\special{ar 3370 2190 10 10 0  6.28318530717959E+0000}%
%
\special{pn 8}%
\special{pa 3120 2010}%
\special{pa 3370 2190}%
\special{dt 0.045}%
\special{sh 1}%
\special{pa 3370 2190}%
\special{pa 3328 2136}%
\special{pa 3328 2160}%
\special{pa 3304 2168}%
\special{pa 3370 2190}%
\special{fp}%
\put(33.1000,-16.0000){\makebox(0,0)[rt]{$\frac{\pi}{2}e_1$}}%
%
\special{pn 8}%
\special{pa 2950 2200}%
\special{pa 3200 2380}%
\special{dt 0.045}%
\special{sh 1}%
\special{pa 3200 2380}%
\special{pa 3158 2326}%
\special{pa 3158 2350}%
\special{pa 3134 2358}%
\special{pa 3200 2380}%
\special{fp}%
\put(31.0000,-18.7000){\makebox(0,0)[rt]{$\arctan\sqrt{2q+1}e_1$}}%
\put(48.6000,-30.0000){\makebox(0,0)[rt]{$\beta_2^{-1}(\frac{\pi}{2})$}}%
%
\special{pn 13}%
\special{sh 1}%
\special{ar 4380 2210 10 10 0  6.28318530717959E+0000}%
\special{sh 1}%
\special{ar 4380 2210 10 10 0  6.28318530717959E+0000}%
%
\special{pn 8}%
\special{pa 4660 2040}%
\special{pa 4390 2200}%
\special{dt 0.045}%
\special{sh 1}%
\special{pa 4390 2200}%
\special{pa 4458 2184}%
\special{pa 4436 2174}%
\special{pa 4438 2150}%
\special{pa 4390 2200}%
\special{fp}%
\put(57.1000,-13.5000){\makebox(0,0)[rt]{$(\beta_1+\beta_2)^{-1}(\pi)$}}%
%
\special{pn 13}%
\special{pa 3970 2650}%
\special{pa 4210 2640}%
\special{fp}%
\special{sh 1}%
\special{pa 4210 2640}%
\special{pa 4144 2624}%
\special{pa 4158 2642}%
\special{pa 4144 2664}%
\special{pa 4210 2640}%
\special{fp}%
%
\special{pn 13}%
\special{pa 4430 2280}%
\special{pa 4430 2480}%
\special{fp}%
\special{sh 1}%
\special{pa 4430 2480}%
\special{pa 4450 2414}%
\special{pa 4430 2428}%
\special{pa 4410 2414}%
\special{pa 4430 2480}%
\special{fp}%
%
\special{pn 13}%
\special{pa 4220 2390}%
\special{pa 4340 2560}%
\special{fp}%
\special{sh 1}%
\special{pa 4340 2560}%
\special{pa 4318 2494}%
\special{pa 4310 2516}%
\special{pa 4286 2518}%
\special{pa 4340 2560}%
\special{fp}%
\put(49.9000,-22.9000){\makebox(0,0)[rt]{\scriptsize{$(q\to\infty)$}}}%
%
\special{pn 8}%
\special{pa 4150 2850}%
\special{pa 4470 2500}%
\special{fp}%
\put(42.1000,-24.2000){\makebox(0,0)[rt]{\scriptsize{$(q\to\infty)$}}}%
\put(42.8000,-26.8000){\makebox(0,0)[rt]{\scriptsize{$(q\to\infty)$}}}%
%
\special{pn 13}%
\special{pa 3450 2180}%
\special{pa 3610 2040}%
\special{fp}%
\special{sh 1}%
\special{pa 3610 2040}%
\special{pa 3548 2070}%
\special{pa 3570 2076}%
\special{pa 3574 2100}%
\special{pa 3610 2040}%
\special{fp}%
\put(41.2000,-20.9000){\makebox(0,0)[rt]{\scriptsize{$(q\to\infty)$}}}%
%
\special{pn 8}%
\special{pa 4470 1710}%
\special{pa 4160 2310}%
\special{dt 0.045}%
\special{sh 1}%
\special{pa 4160 2310}%
\special{pa 4208 2260}%
\special{pa 4184 2264}%
\special{pa 4174 2242}%
\special{pa 4160 2310}%
\special{fp}%
\put(52.6000,-16.2000){\makebox(0,0)[rt]{$\alpha_1e_1+\alpha_2e_2$}}%
\end{picture}%
\hspace{6truecm}}

\vspace{1truecm}

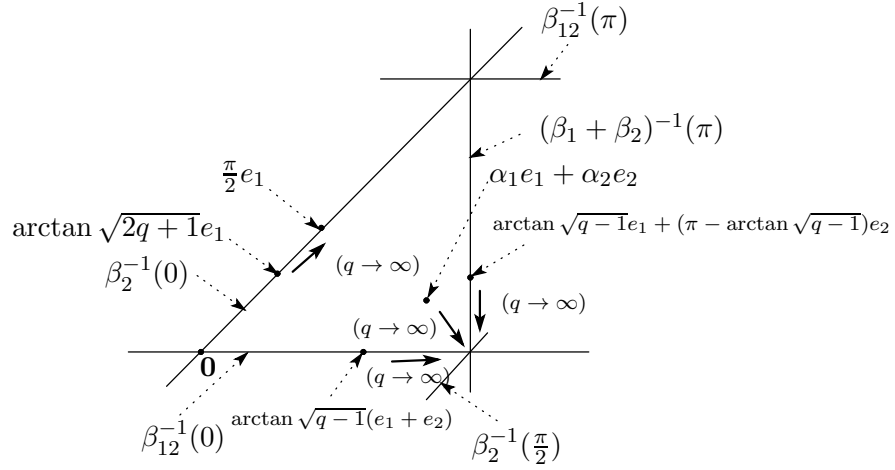
\centerline{{\bf Fig. 7.}}

\vspace{0.5truecm}

\noindent
Easily we can show ${\rm div}\,X\,>\,0$ (on $\widetilde C$).  
Therefore, since $X$ has the only zero point and $X$ is as in Fig.1 over 
a collar neighborhood of $\widetilde{\sigma}_i\,(i=1,2,3)$ 
by the proof of Theorem A, all the integral curves 
of $X$ through points other than its zero point converge to points of 
$\partial\widetilde C$ in finite time.  
Also, we can show ${\rm div}\,X^{\widetilde{\sigma}_i}\,>\,0$ on 
$\widetilde{\sigma}_i$ ($i=1,2,3$).  Hence all integral curves of 
$X^{\widetilde{\sigma}_i}$ through points other than its zero point converge 
to points of $\partial\widetilde{\sigma}_i$ in finite time.  
Therefore we obtain the following fact.  

\vspace{0.5truecm}

\noindent
{\bf Proposition 4.3.} {\sl The mean curvature flow having any non-minimal 
principal orbit of the $S(U(2)\times U(q))$-action 
(on $SU(2+q)/S(U(2)\times U(q))$) as initial data converges to some singular 
orbit in finite time.  Also, the mean curvature flow having any non-minimal 
singular orbit through $\exp^{\perp}(\widetilde{\sigma}_i)$ of the 
$S(U(2)\times U(q))$-action as initial data converges to one of two singular 
orbits through $\exp^{\perp}(\partial\widetilde{\sigma}_i)$ in finite time.}

\vspace{0.5truecm}

\noindent
{\it Example 4.} 
We consider the isotropy action 
$\displaystyle{S(U(p)\times U(p))\curvearrowright 
SU(2p)/S(U(p)\times U(p))}$.  
Then, since $\triangle'$ is the root system of 
$SU(2p)/S(U(p)\times U(p))$, it is of $({\rm C}_p)$-type 
and $\mathfrak b,\triangle'_+,\Pi$ and $\delta$ are described as 
$$\begin{array}{l}
\displaystyle{\mathfrak b={\rm Span}\{e_1,\cdots,e_p\},}\\
\displaystyle{\triangle'_+=\{\beta_i-\beta_j\,\vert\,
1\leq i<j\leq p\}\cup\{\beta_i+\beta_j\,\vert\,1\leq i<j\leq p\}}\\
\hspace{1.2truecm}\displaystyle{\cup\{2\beta_i\,\vert\,1\leq i\leq p\},}\\
\displaystyle{\Pi=\{\beta_i-\beta_{i+1}\,\vert\,1\leq i\leq p-1\}
\cup\{2\beta_p\},\,\,\delta=2\beta_1,}
\end{array}
\leqno{(4.19)}$$
where $(e_1,\cdots,e_p)$ is the orthonormal base of $\mathfrak b$ and 
$(\beta_1,\cdots,\beta_p)$ is the dual base of $(e_1,\cdots,e_p)$.  
For simplicity, we set $\beta_{ij}:=\beta_i-\beta_j\,\,(1\leq i<j\leq p),\,
I:=\{1,\cdots,p\}$ and 
$\widehat I:=\{(i,j)\in I^2\,\vert\,1\leq i<j\leq p\}$.  
It is shown that $\beta_{ij}$'s ($(i,j)\in\hat I$) and 
$\beta_i+\beta_j$'s ($(i,j)\in\hat I$) are of multiplicity $2$ and 
$2\beta_i$'s ($i\in I$) are of multiplicity $1$.  
Since we consider the isotropy action, we have 
${\triangle'}^V_+={\triangle'}_+$ and ${\triangle'}^H_+=\emptyset$, that is, 
$$\begin{array}{l}
\displaystyle{\widetilde C=\{{\bf x}\in\mathfrak b\,\vert\,0<\beta_{ij}
({\bf x})<\pi\,\,((i,j)\in\widehat I),\,\,
0<(\beta_i+\beta_j)({\bf x})<\pi\,\,((i,j)\in\widehat I),}\\
\hspace{5.5truecm}\displaystyle{0<2\beta_i({\bf x})<\pi\,\,(i\in I)\}}\\
\hspace{0.4truecm}\displaystyle{=\{{\bf x}\in\mathfrak b\,\vert\,
\beta_{i,i+1}({\bf x})>0\,(1\leq i\leq p-1),\,2\beta_p({\bf x})>0,\,\,
2\beta_1({\bf x})<\pi\}.}
\end{array}$$
From $(4.5)$, we can describe $X$ explicitly as 
$$\begin{array}{l}
\displaystyle{X_{\bf x}=-\sum_{(i,j)\in\widehat I}
2\cot(x_i-x_j)(e_i-e_j)-\sum_{(i,j)\in\widehat I}
2\cot(x_i+x_j)(e_i+e_j)}\\
\hspace{1.2truecm}\displaystyle{-\sum_{i\in I}\cot2x_i\cdot2e_i}\\
\hspace{0.7truecm}\displaystyle{
=-2\sum_{i\in I}\left(\sum_{j\in I\setminus\{i\}}(\cot(x_i-x_j)+\cot(x_i+x_j))
+\cot2x_i\right)e_i.}
\end{array}$$
According to this relation, a principal orbit 
$S(U(p)\times U(p))({\rm Exp}\,{\bf x})$ (${\bf x}\in\widetilde C$) is minimal 
if and only if the following relations hold:
$$\sum_{j\in I\setminus\{i\}}(\cot(x_i-x_j)+\cot(x_i+x_j))+\cot2x_i=0
\,\,\,\,(i\in I).
\leqno{(4.20)}$$
In the sequel, we consider the case of $p=2$.  Then $(4.20)$ is as follows:
$$\left\{
\begin{array}{l}
\displaystyle{\cot(x_1-x_2)+\cot(x_1+x_2)+\cot2x_1=0}\\
\displaystyle{\cot(x_2-x_1)+\cot(x_1+x_2)+\cot2x_2=0.}
\end{array}\right.\leqno{(4.21)}$$
Under $0<x_1-x_2<\pi,\,0<x_1+x_2<\pi,\,0<2x_1<\pi$ and $0<2x_2<\pi$, 
the equation $(4.21)$ has the only solution, which we denote by 
$(\alpha_1,\alpha_2)$.  
Therefore the orbit $S(U(2)\times U(2))
({\rm Exp}(\alpha_1e_1+\alpha_2e_2))$ is the only minimal principal orbit 
of the $S(U(2)\times U(2))$-action.  
Denote by $\widetilde{\sigma}_1,\,\widetilde{\sigma}_2$ and 
$\widetilde{\sigma}_3$ the one dimensional stratums of 
$\partial\widetilde C$ which are contained in $\beta_{12}^{-1}(0),\,
(2\beta_2)^{-1}(0)$ and $(2\beta_1)^{-1}(\pi)$.  Then we have 
$$\begin{array}{l}
\displaystyle{X^{\widetilde{\sigma}_1}_{\bf x}=-4\cot2x_1(e_1+e_2),}\\
\displaystyle{X^{\widetilde{\sigma}_2}_{\bf x}=-2(2\cot x_1+\cot2x_1)e_1,}\\
\displaystyle{X^{\widetilde{\sigma}_3}_{\bf x}=-2\cot2x_2e_2.}
\end{array}$$
Hence the orbit $S(U(2)\times U(2))({\rm Exp}\,{\bf x})$ 
(${\bf x}\in \widetilde{\sigma}_1$) is minimal if and only if 
$\cot2x_1=0$ ($x_2=x_1$) holds, that is, 
$(x_1,x_2)=(\frac{\pi}{4},\frac{\pi}{4})$.  
Therefore the orbit $S(U(2)\times U(2))$\newline
$({\rm Exp}(\frac{\pi}{4}(e_1+e_2)))$ is the only minimal singular orbit 
through $\widetilde{\sigma}_1$.  
Also, the orbit $S(U(2)\times U(2))({\rm Exp}\,{\bf x})$ 
(${\bf x}\in \widetilde{\sigma}_2$) is minimal if and only if 
$2\cot x_1+\cot2x_1=0$ ($x_2=0$) holds, that is, 
$(x_1,x_2)=(\arctan\sqrt{5},0)$.  Therefore the orbit $S(U(2)\times U(2))
({\rm Exp}(\arctan\sqrt{5}e_1))$ is the only minimal singular orbit 
through $\widetilde{\sigma}_2$.  
Also, the orbit $S(U(2)\times U(2))({\rm Exp}\,{\bf x})$ 
(${\bf x}\in \widetilde{\sigma}_3$) is minimal if and only if 
$\cot2x_2=0$ ($x_1=\frac{\pi}{2}$) holds, that is, 
$(x_1,x_2)=(\frac{\pi}{2},\frac{\pi}{4})$.  
Therefore the orbit $S(U(2)\times U(2))
({\rm Exp}(\frac{\pi}{4}(2e_1+e_2)))$ 
is the only minimal singular orbit through $\widetilde{\sigma}_3$.  

\vspace{1truecm}

\centerline{
\unitlength 0.1in
\begin{picture}( 39.4000, 22.1000)( 10.6000,-30.0000)
%
\special{pn 8}%
\special{pa 2590 2600}%
\special{pa 5000 2600}%
\special{fp}%
%
\special{pn 8}%
\special{pa 3910 1170}%
\special{pa 4850 1170}%
\special{fp}%
%
\special{pn 13}%
\special{sh 1}%
\special{ar 3940 2230 10 10 0  6.28318530717959E+0000}%
\special{sh 1}%
\special{ar 3940 2230 10 10 0  6.28318530717959E+0000}%
%
\special{pn 13}%
\special{sh 1}%
\special{ar 3820 2600 10 10 0  6.28318530717959E+0000}%
\special{sh 1}%
\special{ar 3820 2600 10 10 0  6.28318530717959E+0000}%
%
\special{pn 13}%
\special{sh 1}%
\special{ar 3600 1950 10 10 0  6.28318530717959E+0000}%
\special{sh 1}%
\special{ar 3600 1950 10 10 0  6.28318530717959E+0000}%
%
\special{pn 8}%
\special{pa 3650 2830}%
\special{pa 3812 2610}%
\special{dt 0.045}%
\special{sh 1}%
\special{pa 3812 2610}%
\special{pa 3756 2652}%
\special{pa 3780 2652}%
\special{pa 3788 2676}%
\special{pa 3812 2610}%
\special{fp}%
%
\special{pn 8}%
\special{pa 2860 2930}%
\special{pa 3212 2606}%
\special{dt 0.045}%
\special{sh 1}%
\special{pa 3212 2606}%
\special{pa 3148 2638}%
\special{pa 3172 2642}%
\special{pa 3176 2666}%
\special{pa 3212 2606}%
\special{fp}%
%
\special{pn 8}%
\special{pa 3330 1750}%
\special{pa 3586 1954}%
\special{dt 0.045}%
\special{sh 1}%
\special{pa 3586 1954}%
\special{pa 3546 1898}%
\special{pa 3544 1922}%
\special{pa 3520 1928}%
\special{pa 3586 1954}%
\special{fp}%
%
\special{pn 8}%
\special{pa 4640 1460}%
\special{pa 4380 1580}%
\special{dt 0.045}%
\special{sh 1}%
\special{pa 4380 1580}%
\special{pa 4450 1570}%
\special{pa 4428 1558}%
\special{pa 4432 1534}%
\special{pa 4380 1580}%
\special{fp}%
%
\special{pn 8}%
\special{pa 4520 2960}%
\special{pa 4228 2764}%
\special{dt 0.045}%
\special{sh 1}%
\special{pa 4228 2764}%
\special{pa 4272 2818}%
\special{pa 4272 2794}%
\special{pa 4294 2784}%
\special{pa 4228 2764}%
\special{fp}%
%
\special{pn 8}%
\special{pa 4890 990}%
\special{pa 4750 1170}%
\special{dt 0.045}%
\special{sh 1}%
\special{pa 4750 1170}%
\special{pa 4808 1130}%
\special{pa 4784 1128}%
\special{pa 4776 1106}%
\special{pa 4750 1170}%
\special{fp}%
\put(40.3000,-28.6000){\makebox(0,0)[rt]{\scriptsize{$\arctan\sqrt{5}e_1$}}}%
\put(31.0000,-29.6000){\makebox(0,0)[rt]{$(2\beta_2)^{-1}(0)$}}%
\put(54.3000,-16.7000){\makebox(0,0)[rt]{$\frac{\pi}{4}(2e_1+e_2)$}}%
%
\put(46.6000,-20.3000){\makebox(0,0)[lb]{}}%
\put(29.1000,-21.0000){\makebox(0,0)[rt]{$\beta_{12}^{-1}(0)$}}%
\put(54.1000,-7.9000){\makebox(0,0)[rt]{$(2\beta_2)^{-1}(\pi)$}}%
%
\special{pn 13}%
\special{sh 1}%
\special{ar 2970 2600 10 10 0  6.28318530717959E+0000}%
\special{sh 1}%
\special{ar 2970 2600 10 10 0  6.28318530717959E+0000}%
\put(30.5000,-26.3000){\makebox(0,0)[rt]{${\bf 0}$}}%
%
\special{pn 8}%
\special{pa 4380 2810}%
\special{pa 4380 910}%
\special{fp}%
%
\special{pn 8}%
\special{pa 2790 2780}%
\special{pa 4650 900}%
\special{fp}%
\put(33.1000,-16.0000){\makebox(0,0)[rt]{$\frac{\pi}{4}(e_1+e_2)$}}%
%
\special{pn 8}%
\special{pa 2950 2180}%
\special{pa 3200 2360}%
\special{dt 0.045}%
\special{sh 1}%
\special{pa 3200 2360}%
\special{pa 3158 2306}%
\special{pa 3158 2330}%
\special{pa 3134 2338}%
\special{pa 3200 2360}%
\special{fp}%
\put(47.1000,-30.0000){\makebox(0,0)[rt]{$(\beta_{12})^{-1}(\pi)$}}%
%
\special{pn 13}%
\special{sh 1}%
\special{ar 4380 1940 10 10 0  6.28318530717959E+0000}%
\special{sh 1}%
\special{ar 4380 1940 10 10 0  6.28318530717959E+0000}%
%
\special{pn 8}%
\special{pa 4660 1780}%
\special{pa 4390 1940}%
\special{dt 0.045}%
\special{sh 1}%
\special{pa 4390 1940}%
\special{pa 4458 1924}%
\special{pa 4436 1914}%
\special{pa 4438 1890}%
\special{pa 4390 1940}%
\special{fp}%
\put(54.1000,-13.6000){\makebox(0,0)[rt]{$(2\beta_1)^{-1}(\pi)$}}%
%
\special{pn 8}%
\special{pa 4150 2850}%
\special{pa 4470 2500}%
\special{fp}%
%
\special{pn 8}%
\special{pa 4650 2100}%
\special{pa 3940 2220}%
\special{dt 0.045}%
\special{sh 1}%
\special{pa 3940 2220}%
\special{pa 4010 2230}%
\special{pa 3994 2212}%
\special{pa 4002 2190}%
\special{pa 3940 2220}%
\special{fp}%
\put(55.5000,-20.0000){\makebox(0,0)[rt]{$\alpha_1e_1+\alpha_2e_2$}}%
\end{picture}%
\hspace{6.5truecm}}

\vspace{1truecm}

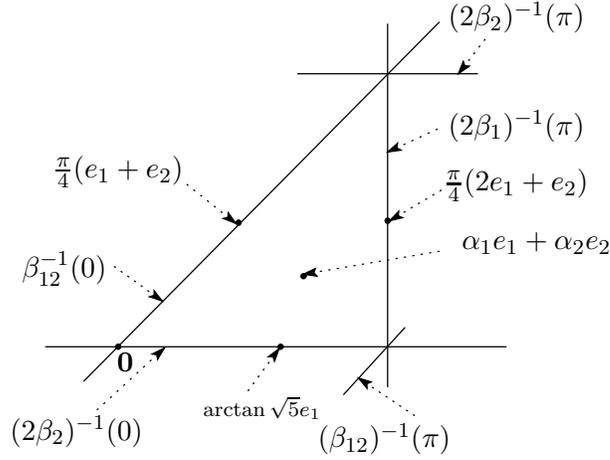
\centerline{{\bf Fig. 8.}}

\vspace{0.5truecm}

\noindent
Easily we can show ${\rm div}\,X\,>\,0$ (on $\widetilde C$).  
Therefore, since $X$ has the only zero point and $X$ is as in Fig.1 over 
a collar neighborhood of $\widetilde{\sigma}_i\,(i=1,2,3)$ 
by the proof of Theorem A, all the integral curves 
of $X$ through points other than its zero point converge to points of 
$\partial\widetilde C$ in finite time.  
Also, we can show ${\rm div}\,X^{\widetilde{\sigma}_i}\,>\,0$ on 
$\widetilde{\sigma}_i$ ($i=1,2,3$).  Hence all integral curves of 
$X^{\widetilde{\sigma}_i}$ through points other than its zero point converge 
to points of $\partial\widetilde{\sigma}_i$ in finite time.  
Therefore we obtain the following fact.  

\vspace{0.5truecm}

\noindent
{\bf Proposition 4.4.} {\sl The mean curvature flow having any non-minimal 
principal orbit of the $S(U(2)\times U(2))$-action 
(on $SU(4)/S(U(2)\times U(2))$) as initial data converges to some singular 
orbit in finite time.  Also, the mean curvature flow having any non-minimal 
singular orbit through $\exp^{\perp}(\widetilde{\sigma}_i)$ of the 
$S(U(2)\times U(2))$-action as initial data converges to one of two singular 
orbits through $\exp^{\perp}(\partial\widetilde{\sigma}_i)$ in finite time.}

\vspace{0.5truecm}

\noindent
{\it Example 5.} 
We consider the Hermann action 
$\displaystyle{SO(p+q)\curvearrowright SU(p+q)/S(U(p)\times U(q))}$ ($p<q$).  
Since $L/H\cap K$ is then equal to $SO(p+q)/SO(p)\times SO(q)$, the 
cohomogeneity of this action is equal to $p$, that is, it is the rank of 
$SU(p+q)/S(U(p)\times U(q))$.  Hence a maximal abelian subspace $\mathfrak b$ 
of $\mathfrak p\cap\mathfrak q$ is also a maximal abelian subspace of 
$\mathfrak p$ and hence 
$\triangle'$ is the root system of $SU(p+q)/S(U(p)\times U(q))$.  
Hence $\triangle'$ is of $(B_p)$-type and the quantities 
$\mathfrak b,\,\triangle'_+,\,\Pi$ and $\delta$ are as in $(4.16)$.  Let 
$\beta_i,\,\beta_{ij},\,I$ and $\widehat I$ be as in Example 3.  
We have ${\rm dim}\,\mathfrak p_{\beta_{ij}}={\rm dim}\,
\mathfrak p_{\beta_i+\beta_j}=2\,((i,j)\in\hat I),\,{\rm dim}\,
\mathfrak p_{\beta_i}=2(q-p)\,(i\in I)$ and 
${\rm dim}\,\mathfrak p_{2\beta_i}=1\,(i\in I)$.  
Also, since ${\triangle'}^V$ is the root system 
of $SO(p+q)/SO(p)\times SO(q)$, it is of $(B_p)$-type.  
The Satake diagram of $SO(p+q)/SO(p)\times SO(q)$ is as follows:
$$\left\{
\begin{array}{ll}
\displaystyle{\begin{xy}
\tiny
\ar@{-} (0,0)*\cir<3pt>{}="A";(3,0)
\ar@{.} (4,0);(6,0)
\ar@{-} (7,0);(9,0)*\cir<3pt>{}="B"
\ar@{-} "B";(12,0)
\ar@{-}  (12,0)*\txt{\Large{$\bullet$}};(15,0)
\ar@{.} (16,0);(18,0)
\ar@{-} (19,0);(21,0)
\ar@{=>} (21,0)*\txt{\Large{$\bullet$}};(25,0)
\ar@{} (25,0)*\txt{\Large{$\bullet$}}
\end{xy}
}&
\displaystyle{(p+q\,:\,{\rm odd})}\\
\displaystyle{
\begin{xy}
\tiny
\ar@{-} (0,0)*\cir<3pt>{}="A";(3,0)
\ar@{.} (4,0);(6,0)
\ar@{-} (7,0);(9,0)*\cir<3pt>{}="B"
\ar@{-} "B";(12,0)
\ar@{}  (12,0)*\txt{\Large{$\bullet$}}
\ar@{-} (12,0);(14,0)
\ar@{.} (15,0);(17,0)
\ar@{-} (18,0);(20,0)
\ar@{}  (20,0)*\txt{\Large{$\bullet$}}
\ar@{-} (20,0);(22.536,2.536)
\ar@{-} (20,0);(22.536,-2.536)
\ar@{} (22.536,2.536)*\txt{\Large{$\bullet$}}
\ar@{} (22.536,-2.536)*\txt{\Large{$\bullet$}}
\end{xy}
}&
\displaystyle{(p+q\,:\,{\rm even},\,\,q-p\geq 4)}\\
\displaystyle{
\begin{xy}
\tiny
\ar@{-} (0,0)*\cir<3pt>{}="A";(3,0)
\ar@{.} (4,0);(7,0)
\ar@{-} (8,0);(10,0)*\cir<3pt>{}="B"
\ar@{-} "B";(13.536,2.536)*\cir<3pt>{}="C"
\ar@{-} "B";(13.536,-2.5)*\cir<3pt>{}="D"
\ar @/^/@{<->}"C";"D"
\end{xy}}&
\displaystyle{(p+q\,:\,{\rm even},\,\,q-p=2).}
\end{array}\right.$$
Hence we have ${\triangle'}^V_+={\triangle'}_+\setminus
\{2\beta_i\,\vert\,i\in I\},\,{\rm dim}(\mathfrak p_{\beta_{ij}}\cap
\mathfrak q)={\rm dim}(\mathfrak p_{\beta_i+\beta_j}\cap\mathfrak q)=1$ 
($(i,j)\in\hat I$) and ${\rm dim}(\mathfrak p_{\beta_i}\cap\mathfrak q)=q-p$ 
($i\in I$).  
Furthermore we have ${\triangle'}^H_+=\triangle'_+,\,{\rm dim}
(\mathfrak p_{\beta_{ij}}\cap\mathfrak h)={\rm dim}
(\mathfrak p_{\beta_i+\beta_j}\cap\mathfrak h)=1$ ($(i,j)\in\hat I$), 
${\rm dim}(\mathfrak p_{\beta_i}\cap\mathfrak h)=q-p$ ($i\in I$) and 
${\rm dim}(\mathfrak p_{2\beta_i}\cap\mathfrak h)=1$ ($i\in I$).  
Therefore we have 
$$\widetilde C=\{{\bf x}\in\mathfrak b\,\vert\,
\beta_{i,i+1}({\bf x})>0\,\,(i=1,\cdots,p-1),\,\,\beta_p({\bf x})>0,\,\,
(\beta_1+\beta_2)({\bf x})<\frac{\pi}{2}\}.$$
From $(4.5)$, we can describe $X$ explicitly as 
$$\begin{array}{l}
\displaystyle{X_{\bf x}=\sum_{i\in I}2\tan2x_i\cdot e_i
-\sum_{(i,j)\in\widehat I}
2\cot2(x_i-x_j)(e_i-e_j)}\\
\hspace{1.2truecm}\displaystyle{-\sum_{i\in I}2(q-p)\cot2x_i\cdot e_i
-\sum_{(i,j)\in\widehat I}2\cot2(x_i+x_j)(e_i+e_j)}\\
\hspace{0.7truecm}\displaystyle{=2\sum_{i\in I}\left(
\tan2x_i-(q-p)\cot2x_i-\sum_{i\in I\setminus\{i\}}
(\cot2(x_i-x_j)+\cot2(x_i+x_j))\right)e_i.}
\hspace{6truecm}\displaystyle{({\bf x}(=\sum_{i\in I}x_ie_i)\in\widetilde C).}
\end{array}$$
According to this relation, a principal orbit 
$SO(p+q)({\rm Exp}\,{\bf x})$ (${\bf x}\in\widetilde C$) is minimal 
if and only if the following relations hold:
$$\tan2x_i-(q-p)\cot2x_i-\sum_{j\in I\setminus\{i\}}
(\cot2(x_i-x_j)+\cot2(x_i+x_j))=0\,\,\,\,(i\in I).
\leqno{(4.22)}$$
In the sequel, we consider the case of $p=2$.  Then $(4.22)$ is as follows:
$$\left\{
\begin{array}{l}
\displaystyle{\tan2x_1-(q-2)\cot2x_1-\cot2(x_1-x_2)-\cot2(x_1+x_2)=0}\\
\displaystyle{\tan2x_2-(q-2)\cot2x_2-\cot2(x_2-x_1)-\cot2(x_1+x_2)=0.}
\end{array}
\right.\leqno{(4.23)}$$
Under $0<x_1-x_2<\frac{\pi}{2},\,0<x_1+x_2<\frac{\pi}{2},\,
0<x_i<\frac{\pi}{2}\,(i=1,2)$ and $-\frac{\pi}{2}<2x_i<\frac{\pi}{2}\,
(i=1,2)$, the equation $(4.23)$ has the only solution, which we denote by 
$(\alpha_1,\alpha_2)$.  
Here we note that $\lim\limits_{q\to\infty}(\alpha_1,\alpha_2)
=(\frac{\pi}{4},\frac{\pi}{4})$.  
Denote by $\widetilde{\sigma}_1,\,\widetilde{\sigma}_2$ and 
$\widetilde{\sigma}_3$ the one dimensional stratums of 
$\partial\widetilde C$ which are contained in $\beta_{12}^{-1}(0),\,
\beta_2^{-1}(0)$ and $(\beta_1+\beta_2)^{-1}(\frac{\pi}{2})$, respectively.  
Then we have 
$$\begin{array}{l}
\displaystyle{X^{\widetilde{\sigma}_1}_{\bf x}
=2(\tan2x_1-(q-2)\cot2x_1-\cot4x_1)(e_1+e_2),}\\
\displaystyle{X^{\widetilde{\sigma}_2}_{\bf x}=2(\tan2x_1-q\cot2x_1)e_1,}\\
\displaystyle{X^{\widetilde{\sigma}_3}_{\bf x}=2(\tan2x_1-(q-2)\cot2x_1
-\cot4x_1)(e_1-e_2).}
\end{array}$$
Hence the orbit $SO(q+2)({\rm Exp}\,{\bf x})$ 
(${\bf x}\in \widetilde{\sigma}_1$) is minimal if and only if 
$\tan2x_1-(q-2)\cot2x_1-\cot4x_1=0$ ($x_2=x_1$) holds, that is, 
$(x_1,x_2)=(\frac12\arctan\sqrt{\frac23q-1},$\newline
$\frac12\arctan\sqrt{\frac23q-1})$.  
Hence the orbit $SO(q+2)({\rm Exp}(\frac12\arctan\sqrt{\frac23q-1}(e_1+e_2)))$ 
is the only minimal singular orbit through $\widetilde{\sigma}_1$.  
Also, the orbit $SO(q+2)({\rm Exp}\,{\bf x})$ 
(${\bf x}\in \widetilde{\sigma}_2$) is minimal if and only if 
$\tan2x_1-q\cot2x_1=0$ ($x_2=0$) holds, that is, 
$(x_1,x_2)=(\frac12\arctan\sqrt q,0)$.  
Hence the orbit $SO(q+2)({\rm Exp}(\frac12\arctan\sqrt q e_1))$ 
is the only minimal singular orbit through $\widetilde{\sigma}_2$.  
Also, the orbit $SO(q+2)({\rm Exp}\,{\bf x})$ 
(${\bf x}\in \widetilde{\sigma}_3$) is minimal if and only if 
$\tan2x_1-(q-2)\cot2x_1-\cot4x_1=0$ ($x_2=\frac{\pi}{2}-x_1$) holds, that is, 
$(x_1,x_2)=(\frac{\pi}{2}-\frac12\arctan\sqrt{\frac23q-1},
\frac12\arctan\sqrt{\frac23q-1})$.  
Hence the orbit $SO(q+2)({\rm Exp}
((\frac{\pi}{2}-\frac12\arctan\sqrt{\frac23q-1})e_1
+\frac12\arctan\sqrt{\frac23q-1}e_2)$ 
is the only minimal singular orbit through $\widetilde{\sigma}_3$.  

\vspace{0.2truecm}

\centerline{
\unitlength 0.1in
\begin{picture}( 66.9000, 24.4000)(-16.9000,-29.5000)
%
\special{pn 8}%
\special{pa 2590 2600}%
\special{pa 5000 2600}%
\special{fp}%
%
\special{pn 8}%
\special{pa 3910 1170}%
\special{pa 4850 1170}%
\special{fp}%
%
\special{pn 13}%
\special{sh 1}%
\special{ar 4150 2330 10 10 0  6.28318530717959E+0000}%
\special{sh 1}%
\special{ar 4150 2330 10 10 0  6.28318530717959E+0000}%
%
\special{pn 13}%
\special{sh 1}%
\special{ar 3820 2600 10 10 0  6.28318530717959E+0000}%
\special{sh 1}%
\special{ar 3820 2600 10 10 0  6.28318530717959E+0000}%
%
\special{pn 13}%
\special{sh 1}%
\special{ar 3600 1950 10 10 0  6.28318530717959E+0000}%
\special{sh 1}%
\special{ar 3600 1950 10 10 0  6.28318530717959E+0000}%
%
\special{pn 8}%
\special{pa 3650 2830}%
\special{pa 3812 2610}%
\special{dt 0.045}%
\special{sh 1}%
\special{pa 3812 2610}%
\special{pa 3756 2652}%
\special{pa 3780 2652}%
\special{pa 3788 2676}%
\special{pa 3812 2610}%
\special{fp}%
%
\special{pn 8}%
\special{pa 2860 2930}%
\special{pa 3212 2606}%
\special{dt 0.045}%
\special{sh 1}%
\special{pa 3212 2606}%
\special{pa 3148 2638}%
\special{pa 3172 2642}%
\special{pa 3176 2666}%
\special{pa 3212 2606}%
\special{fp}%
%
\special{pn 8}%
\special{pa 3330 1750}%
\special{pa 3586 1954}%
\special{dt 0.045}%
\special{sh 1}%
\special{pa 3586 1954}%
\special{pa 3546 1898}%
\special{pa 3544 1922}%
\special{pa 3520 1928}%
\special{pa 3586 1954}%
\special{fp}%
%
\special{pn 8}%
\special{pa 4640 1460}%
\special{pa 4380 1580}%
\special{dt 0.045}%
\special{sh 1}%
\special{pa 4380 1580}%
\special{pa 4450 1570}%
\special{pa 4428 1558}%
\special{pa 4432 1534}%
\special{pa 4380 1580}%
\special{fp}%
%
\special{pn 8}%
\special{pa 4500 2950}%
\special{pa 4174 2818}%
\special{dt 0.045}%
\special{sh 1}%
\special{pa 4174 2818}%
\special{pa 4228 2862}%
\special{pa 4224 2838}%
\special{pa 4244 2824}%
\special{pa 4174 2818}%
\special{fp}%
%
\special{pn 8}%
\special{pa 4890 990}%
\special{pa 4750 1170}%
\special{dt 0.045}%
\special{sh 1}%
\special{pa 4750 1170}%
\special{pa 4808 1130}%
\special{pa 4784 1128}%
\special{pa 4776 1106}%
\special{pa 4750 1170}%
\special{fp}%
\put(43.5000,-28.7000){\makebox(0,0)[rt]{\scriptsize{$\frac12\arctan\sqrt{\frac23q-1}(e_1+e_2)$}}}%
\put(29.4000,-29.5000){\makebox(0,0)[rt]{$\beta_{12}^{-1}(0)$}}%
\put(73.1000,-19.4000){\makebox(0,0)[rt]{\scriptsize{$(\frac{\pi}{2}-\frac12\arctan\sqrt{\frac23q-1})e_1+\frac12\arctan\sqrt{\frac23q-1}e_2$}}}%
%
\put(46.6000,-20.3000){\makebox(0,0)[lb]{}}%
\put(29.1000,-21.0000){\makebox(0,0)[rt]{$\beta_2^{-1}(0)$}}%
\put(52.8000,-7.9000){\makebox(0,0)[rt]{$\beta_{12}^{-1}(\frac{\pi}{2})$}}%
%
\special{pn 13}%
\special{sh 1}%
\special{ar 2970 2600 10 10 0  6.28318530717959E+0000}%
\special{sh 1}%
\special{ar 2970 2600 10 10 0  6.28318530717959E+0000}%
\put(30.3000,-26.5000){\makebox(0,0)[rt]{${\bf 0}$}}%
%
\special{pn 8}%
\special{pa 4380 2810}%
\special{pa 4380 910}%
\special{fp}%
%
\special{pn 8}%
\special{pa 2790 2780}%
\special{pa 4650 900}%
\special{fp}%
%
\special{pn 13}%
\special{sh 1}%
\special{ar 3370 2190 10 10 0  6.28318530717959E+0000}%
\special{sh 1}%
\special{ar 3370 2190 10 10 0  6.28318530717959E+0000}%
%
\special{pn 8}%
\special{pa 3120 2010}%
\special{pa 3370 2190}%
\special{dt 0.045}%
\special{sh 1}%
\special{pa 3370 2190}%
\special{pa 3328 2136}%
\special{pa 3328 2160}%
\special{pa 3304 2168}%
\special{pa 3370 2190}%
\special{fp}%
\put(33.1000,-16.0000){\makebox(0,0)[rt]{$\frac{\pi}{4}e_1$}}%
%
\special{pn 8}%
\special{pa 2950 2190}%
\special{pa 3200 2370}%
\special{dt 0.045}%
\special{sh 1}%
\special{pa 3200 2370}%
\special{pa 3158 2316}%
\special{pa 3158 2340}%
\special{pa 3134 2348}%
\special{pa 3200 2370}%
\special{fp}%
\put(31.0000,-18.7000){\makebox(0,0)[rt]{$\frac12\arctan\sqrt{q}e_1$}}%
\put(49.6000,-29.3000){\makebox(0,0)[rt]{$\beta_2^{-1}(\frac{\pi}{4})$}}%
%
\special{pn 13}%
\special{sh 1}%
\special{ar 4380 2210 10 10 0  6.28318530717959E+0000}%
\special{sh 1}%
\special{ar 4380 2210 10 10 0  6.28318530717959E+0000}%
%
\special{pn 8}%
\special{pa 4660 2040}%
\special{pa 4390 2200}%
\special{dt 0.045}%
\special{sh 1}%
\special{pa 4390 2200}%
\special{pa 4458 2184}%
\special{pa 4436 2174}%
\special{pa 4438 2150}%
\special{pa 4390 2200}%
\special{fp}%
\put(56.8000,-13.4000){\makebox(0,0)[rt]{$(\beta_1+\beta_2)^{-1}(\frac{\pi}{2})$}}%
%
\special{pn 13}%
\special{pa 3970 2650}%
\special{pa 4210 2640}%
\special{fp}%
\special{sh 1}%
\special{pa 4210 2640}%
\special{pa 4144 2624}%
\special{pa 4158 2642}%
\special{pa 4144 2664}%
\special{pa 4210 2640}%
\special{fp}%
%
\special{pn 13}%
\special{pa 4430 2280}%
\special{pa 4430 2480}%
\special{fp}%
\special{sh 1}%
\special{pa 4430 2480}%
\special{pa 4450 2414}%
\special{pa 4430 2428}%
\special{pa 4410 2414}%
\special{pa 4430 2480}%
\special{fp}%
%
\special{pn 13}%
\special{pa 4220 2390}%
\special{pa 4340 2560}%
\special{fp}%
\special{sh 1}%
\special{pa 4340 2560}%
\special{pa 4318 2494}%
\special{pa 4310 2516}%
\special{pa 4286 2518}%
\special{pa 4340 2560}%
\special{fp}%
\put(49.9000,-22.9000){\makebox(0,0)[rt]{\scriptsize{$(q\to\infty)$}}}%
%
\special{pn 8}%
\special{pa 4150 2850}%
\special{pa 4470 2500}%
\special{fp}%
\put(42.1000,-24.2000){\makebox(0,0)[rt]{\scriptsize{$(q\to\infty)$}}}%
\put(42.8000,-26.8000){\makebox(0,0)[rt]{\scriptsize{$(q\to\infty)$}}}%
%
\special{pn 13}%
\special{pa 3450 2180}%
\special{pa 3610 2040}%
\special{fp}%
\special{sh 1}%
\special{pa 3610 2040}%
\special{pa 3548 2070}%
\special{pa 3570 2076}%
\special{pa 3574 2100}%
\special{pa 3610 2040}%
\special{fp}%
\put(41.2000,-20.9000){\makebox(0,0)[rt]{\scriptsize{$(q\to\infty)$}}}%
%
\special{pn 8}%
\special{pa 3040 2674}%
\special{pa 2752 2372}%
\special{fp}%
%
\special{pn 8}%
\special{pa 2470 2470}%
\special{pa 2830 2460}%
\special{dt 0.045}%
\special{sh 1}%
\special{pa 2830 2460}%
\special{pa 2764 2442}%
\special{pa 2778 2462}%
\special{pa 2764 2482}%
\special{pa 2830 2460}%
\special{fp}%
\put(24.3000,-24.0000){\makebox(0,0)[rt]{$\beta_1^{-1}(0)$}}%
%
\special{pn 8}%
\special{pa 4160 960}%
\special{pa 4510 1290}%
\special{fp}%
%
\special{pn 8}%
\special{pa 4310 690}%
\special{pa 4220 1010}%
\special{dt 0.045}%
\special{sh 1}%
\special{pa 4220 1010}%
\special{pa 4258 952}%
\special{pa 4234 960}%
\special{pa 4220 940}%
\special{pa 4220 1010}%
\special{fp}%
\put(45.6000,-5.1000){\makebox(0,0)[rt]{$\beta_1^{-1}(\frac{\pi}{2})$}}%
%
\special{pn 8}%
\special{pa 4550 1790}%
\special{pa 4160 2320}%
\special{dt 0.045}%
\special{sh 1}%
\special{pa 4160 2320}%
\special{pa 4216 2278}%
\special{pa 4192 2278}%
\special{pa 4184 2254}%
\special{pa 4160 2320}%
\special{fp}%
\put(52.7000,-16.2000){\makebox(0,0)[rt]{$\alpha_1e_1+\alpha_2e_2$}}%
\end{picture}%
\hspace{15truecm}}

\vspace{1truecm}

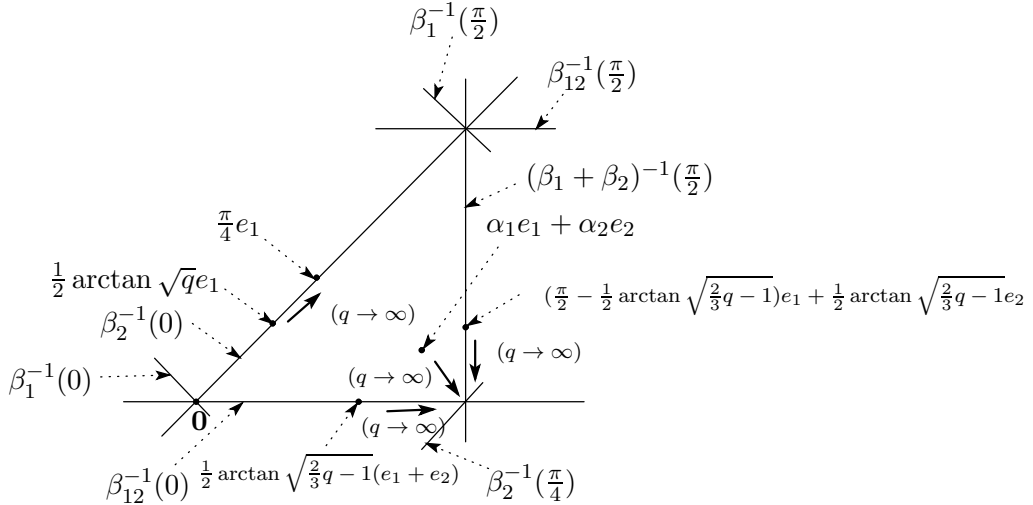
\centerline{{\bf Fig. 9.}}

\vspace{0.5truecm}

\noindent
Easily we can show ${\rm div}\,X\,>\,0$ (on $\widetilde C$).  
Therefore, since $X$ has the only zero point and $X$ is as in Fig.1 over a 
collar neighborhood of $\widetilde{\sigma}_i$ ($i=1,2,3$) 
by the proof of Theorem A, all the integral curves 
of $X$ through points other than its zero point converge to points of 
$\partial\widetilde C$ in finite time.  
Also, we can show ${\rm div}\,X^{\widetilde{\sigma}_i}\,>\,0$ on 
$\widetilde{\sigma}_i$ ($i=1,2,3$).  Hence all integral curves of 
$X^{\widetilde{\sigma}_i}$ through points other than its zero point converge 
to points of $\partial\widetilde{\sigma}_i$ in finite time.  
Therefore we obtain the following fact.  

\vspace{0.2truecm}

\noindent
{\bf Proposition 4.5.} {\sl The mean curvature flow having any non-minimal 
principal orbit of the $SO(q+2)$-action 
(on $SU(q+2)/S(U(2)\times U(q))$) as initial data converges to some singular 
orbit in finite time, where $q\,>\,2$.  
Also, the mean curvature flow having any non-minimal 
singular orbit through $\exp^{\perp}(\widetilde{\sigma}_i)$ of the 
$SO(q+2)$-action as initial data converges to one of two singular 
orbits through $\exp^{\perp}(\partial\widetilde{\sigma}_i)$ in finite time.}

\vspace{0.2truecm}

\noindent
{\it Example 6.} 
We consider the Hermann action 
$\displaystyle{SO(2p)\curvearrowright SU(2p)/S(U(p)\times U(p))}$.  
Since $L/H\cap K$ is then equal to $SO(2p)/SO(p)\times SO(p)$, 
the cohomogeneity of this action is equal to 
$p\,(={\rm rank}(SU(2p)/S(U(p)\times U(p))))$.  
Hence a maximal abelian subspace $\mathfrak b$ of $\mathfrak p\cap\mathfrak q$ 
is also a maximal abelian subspace of $\mathfrak p$ 
and hence $\triangle'$ is the root system of $SU(2p)/S(U(p)\times U(p))$.  
Hence it is of $(C_p)$-type and the quantities $\mathfrak b,\,\triangle'_+,\,
\Pi$ and $\delta$ are as in $(4.19)$.  
Let $\beta_i,\,\beta_{ij},\,I$ and $\widehat I$ be as in 
Example 4.  We have ${\rm dim}\,\mathfrak p_{\beta_{ij}}
={\rm dim}\,\mathfrak p_{\beta_i+\beta_j}=2$ ($(i,j)\in\hat I$) and 
${\rm dim}\,\mathfrak p_{2\beta_i}=1$ ($i\in I$).  
Also, ${\triangle'}^V$ is the root system of $SO(2p)/SO(p)\times SO(p)$.  
Hence it is of $(D_p)$-type and we have 
${\triangle'}^V_+=\{\beta_{ij}\,\vert\,(i,j)\in\hat I\}\cup\{\beta_i+\beta_j\,
\vert\,(i,j)\in\hat I\}$ and 
${\rm dim}(\mathfrak p_{\beta_{ij}}\cap\mathfrak q)
={\rm dim}(\mathfrak p_{\beta_i+\beta_j}\cap\mathfrak q)=1$ 
($(i,j)\in\hat I$).  Hence we have ${\triangle'}^H_+=\triangle'_+$ and each 
root of ${\triangle'}^H_+$ has multiplicity $1$.  Therefore we have 
$$\widetilde C=\{{\bf x}\in\mathfrak b\,\vert\,
\beta_{i,i+1}({\bf x})>0\,\,(i=1,\cdots,p-1),\,\,
(\beta_{p-1}+\beta_p)({\bf x})>0,\,\,2\beta_1({\bf x})<\frac{\pi}{2}\}$$
and 
$$\begin{array}{l}
\displaystyle{X_{\bf x}=-\sum_{(i,j)\in\widehat I}
2\left(\cot2(x_i-x_j)(e_i-e_j)+\cot2(x_i+x_j)(e_i+e_j)\right)}\\
\hspace{1.2truecm}\displaystyle{+\sum_{i\in I}\tan2x_i\cdot2e_i}\\
\hspace{0.7truecm}\displaystyle{
=2\sum_{i\in I}\left(-\sum_{j\in I\setminus\{i\}}
(\cot2(x_i-x_j)+\cot2(x_i+x_j))+\tan2x_i\right)e_i.}
\end{array}$$
According to this relation, a principal orbit 
$SO(2p)({\rm Exp}\,{\bf x})$ (${\bf x}\in\widetilde C$) is minimal 
if and only if the following relation holds:
$$\sum_{j\in I\setminus\{i\}}(\cot2(x_i-x_j)+\cot2(x_i+x_j))-\tan2x_i=0
\,\,\,\,(i\in I).
\leqno{(4.24)}$$
In the sequel, we consider the case of $p=2$.  Then $(4.24)$ is as follows:
$$\left\{
\begin{array}{l}
\displaystyle{\cot2(x_1-x_2)+\cot2(x_1+x_2)-\tan2x_1=0}\\
\displaystyle{\cot2(x_2-x_1)+\cot2(x_1+x_2)-\tan2x_2=0.}
\end{array}\right.\leqno{(4.25)}$$
Under $x_1-x_2>0,\,x_1+x_2>0$ and $x_1<\frac{\pi}{4}$, 
the equation $(4.25)$ has the only solution 
$(x_1,x_2)=(\frac12\arctan\sqrt2,0)$.  Therefore the orbit 
$SO(4){\rm Exp}(\frac12\arctan\sqrt2e_1)$ is the only minmal principal orbit 
of the $SO(4)$-action.  
Denote by $\widetilde{\sigma}_1,\,\widetilde{\sigma}_2$ and 
$\widetilde{\sigma}_3$ the one dimensional stratums of 
$\partial\widetilde C$ which are contained in $\beta_{12}^{-1}(0),\,
(\beta_1+\beta_2)^{-1}(0)$ and $(2\beta_1)^{-1}(\frac{\pi}{2})$, 
respectively.   Then we have 
$$\begin{array}{l}
\displaystyle{X^{\widetilde{\sigma}_1}_{\bf x}=
-2(\cot2(x_1+x_2)-\tan2x_1)e_1-2(\cot2(x_1+x_2)-\tan2x_2)e_2}\\
\displaystyle{X^{\widetilde{\sigma}_2}_{\bf x}=
-2(\cot2(x_1-x_2)-\tan2x_1)e_1-2(\cot2(x_2-x_1)-\tan2x_2)e_2}\\
\displaystyle{X^{\widetilde{\sigma}_3}_{\bf x}=
-2(\cot2(x_1-x_2)+\cot2(x_1+x_2))e_1}\\
\hspace{1.5truecm}\displaystyle{-2(\cot2(x_2-x_1)+\cot2(x_1+x_2)-\tan2x_2)e_2.}
\end{array}$$
Hence the orbit $SO(4)({\rm Exp}\,{\bf x})$ 
(${\bf x}\in \widetilde{\sigma}_1$) is minimal if and only if 
$\cot4x_1=\tan2x_1$ ($x_2=x_1$) holds, that is, 
$(x_1,x_2)=(\frac{\pi}{12},\frac{\pi}{12})$.  
Therefore the orbit $SO(4)({\rm Exp}(\frac{\pi}{12}(e_1+e_2)))$ is 
the only minimal singular orbit through $\widetilde{\sigma}_1$.  
Also, the orbit $SO(4)({\rm Exp}\,{\bf x})$ 
(${\bf x}\in \widetilde{\sigma}_2$) is minimal if and only if 
$\cot4x_1=\tan2x_1$ ($x_2=-x_1$) holds, that is, 
$(x_1,x_2)=(\frac{\pi}{12},-\frac{\pi}{12})$.  Therefore the orbit 
$SO(4)({\rm Exp}(\frac{\pi}{12}(e_1-e_2)))$ is the only minimal singular 
orbit through $\widetilde{\sigma}_2$.  
Also, the orbit $SO(4)({\rm Exp}\,{\bf x})$ 
(${\bf x}\in \widetilde{\sigma}_3$) is minimal if and only if 
$\cot(\frac{\pi}{2}-2x_2)+\cot(\frac{\pi}{2}+2x_2)=0$ and 
$-\cot(\frac{\pi}{2}-2x_2)+\cot(\frac{\pi}{2}+2x_2)-\tan2x_2=0$ 
($x_1=\frac{\pi}{4}$) hold, 
that is, $(x_1,x_2)=(\frac{\pi}{4},0)$.  
Therefore the orbit $SO(4)({\rm Exp}(\frac{\pi}{4}e_1))$ 
is the only minimal singular orbit through $\widetilde{\sigma}_3$.  

\vspace{0.2truecm}

\centerline{
\unitlength 0.1in
\begin{picture}( 52.0000, 20.6000)( -2.0000,-28.5000)
%
\special{pn 8}%
\special{pa 2590 2600}%
\special{pa 5000 2600}%
\special{fp}%
%
\special{pn 8}%
\special{pa 3910 1170}%
\special{pa 4850 1170}%
\special{fp}%
%
\special{pn 13}%
\special{sh 1}%
\special{ar 3940 2240 10 10 0  6.28318530717959E+0000}%
\special{sh 1}%
\special{ar 3940 2240 10 10 0  6.28318530717959E+0000}%
%
\special{pn 13}%
\special{sh 1}%
\special{ar 4160 2600 10 10 0  6.28318530717959E+0000}%
\special{sh 1}%
\special{ar 4160 2600 10 10 0  6.28318530717959E+0000}%
%
\special{pn 13}%
\special{sh 1}%
\special{ar 3600 1950 10 10 0  6.28318530717959E+0000}%
\special{sh 1}%
\special{ar 3600 1950 10 10 0  6.28318530717959E+0000}%
%
\special{pn 8}%
\special{pa 3990 2830}%
\special{pa 4152 2610}%
\special{dt 0.045}%
\special{sh 1}%
\special{pa 4152 2610}%
\special{pa 4096 2652}%
\special{pa 4120 2652}%
\special{pa 4128 2676}%
\special{pa 4152 2610}%
\special{fp}%
%
\special{pn 8}%
\special{pa 3330 1750}%
\special{pa 3586 1954}%
\special{dt 0.045}%
\special{sh 1}%
\special{pa 3586 1954}%
\special{pa 3546 1898}%
\special{pa 3544 1922}%
\special{pa 3520 1928}%
\special{pa 3586 1954}%
\special{fp}%
%
\special{pn 8}%
\special{pa 4640 1460}%
\special{pa 4380 1580}%
\special{dt 0.045}%
\special{sh 1}%
\special{pa 4380 1580}%
\special{pa 4450 1570}%
\special{pa 4428 1558}%
\special{pa 4432 1534}%
\special{pa 4380 1580}%
\special{fp}%
%
\special{pn 8}%
\special{pa 4890 990}%
\special{pa 4750 1170}%
\special{dt 0.045}%
\special{sh 1}%
\special{pa 4750 1170}%
\special{pa 4808 1130}%
\special{pa 4784 1128}%
\special{pa 4776 1106}%
\special{pa 4750 1170}%
\special{fp}%
\put(43.2000,-28.4000){\makebox(0,0)[rt]{$\frac{\pi}{12}(e_1-e_2)$}}%
\put(35.4000,-28.5000){\makebox(0,0)[rt]{$(\beta_1+\beta_2)^{-1}(0)$}}%
\put(54.9000,-21.1000){\makebox(0,0)[rt]{$\frac{\pi}{12}(e_1+e_2)$}}%
%
\put(46.6000,-20.3000){\makebox(0,0)[lb]{}}%
\put(29.1000,-21.0000){\makebox(0,0)[rt]{$(2\beta_1)^{-1}(\frac{\pi}{2})$}}%
\put(57.0000,-7.9000){\makebox(0,0)[rt]{$(\beta_1+\beta_2)^{-1}(\frac{\pi}{2})$}}%
%
\special{pn 13}%
\special{sh 1}%
\special{ar 4380 2600 10 10 0  6.28318530717959E+0000}%
\special{sh 1}%
\special{ar 4380 2600 10 10 0  6.28318530717959E+0000}%
\put(45.5000,-26.5000){\makebox(0,0)[rt]{${\bf 0}$}}%
%
\special{pn 8}%
\special{pa 4380 2810}%
\special{pa 4380 910}%
\special{fp}%
%
\special{pn 8}%
\special{pa 2790 2780}%
\special{pa 4650 900}%
\special{fp}%
\put(33.1000,-16.0000){\makebox(0,0)[rt]{$\frac{\pi}{4}e_1$}}%
%
\special{pn 8}%
\special{pa 2940 2190}%
\special{pa 3190 2370}%
\special{dt 0.045}%
\special{sh 1}%
\special{pa 3190 2370}%
\special{pa 3148 2316}%
\special{pa 3148 2340}%
\special{pa 3124 2348}%
\special{pa 3190 2370}%
\special{fp}%
%
\special{pn 13}%
\special{sh 1}%
\special{ar 4380 2410 10 10 0  6.28318530717959E+0000}%
\special{sh 1}%
\special{ar 4380 2410 10 10 0  6.28318530717959E+0000}%
%
\special{pn 8}%
\special{pa 4660 2240}%
\special{pa 4390 2400}%
\special{dt 0.045}%
\special{sh 1}%
\special{pa 4390 2400}%
\special{pa 4458 2384}%
\special{pa 4436 2374}%
\special{pa 4438 2350}%
\special{pa 4390 2400}%
\special{fp}%
\put(51.6000,-13.6000){\makebox(0,0)[rt]{$\beta_{12}^{-1}(0)$}}%
%
\special{pn 8}%
\special{pa 4380 2600}%
\special{pa 3600 1960}%
\special{dt 0.045}%
%
\special{pn 8}%
\special{pa 4800 1960}%
\special{pa 3930 2230}%
\special{dt 0.045}%
\special{sh 1}%
\special{pa 3930 2230}%
\special{pa 4000 2230}%
\special{pa 3982 2214}%
\special{pa 3988 2192}%
\special{pa 3930 2230}%
\special{fp}%
\put(58.0000,-18.5000){\makebox(0,0)[rt]{$\frac12\arctan\sqrt2e_1$}}%
%
\special{pn 8}%
\special{pa 2970 2290}%
\special{pa 2970 2750}%
\special{fp}%
%
\special{pn 8}%
\special{pa 2730 2390}%
\special{pa 2970 2480}%
\special{dt 0.045}%
\special{sh 1}%
\special{pa 2970 2480}%
\special{pa 2916 2438}%
\special{pa 2920 2462}%
\special{pa 2902 2476}%
\special{pa 2970 2480}%
\special{fp}%
\put(26.8000,-23.2000){\makebox(0,0)[rt]{$\beta_{12}^{-1}(\frac{\pi}{2})$}}%
%
\special{pn 8}%
\special{pa 3320 2810}%
\special{pa 3500 2600}%
\special{dt 0.045}%
\special{sh 1}%
\special{pa 3500 2600}%
\special{pa 3442 2638}%
\special{pa 3466 2640}%
\special{pa 3472 2664}%
\special{pa 3500 2600}%
\special{fp}%
\end{picture}%
\hspace{9truecm}}

\vspace{1truecm}

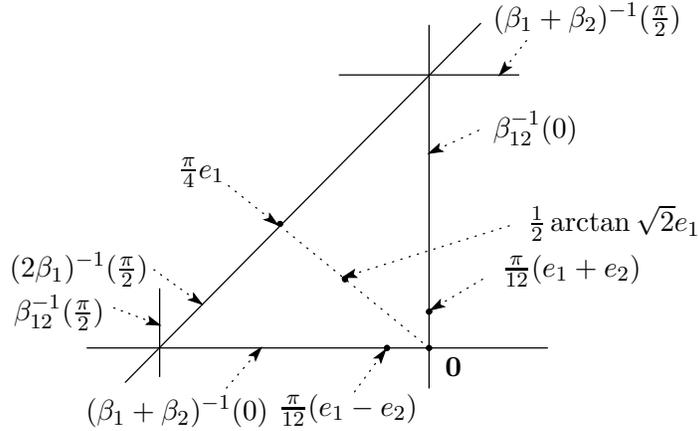
\centerline{{\bf Fig. 10.}}

\vspace{0.5truecm}

\noindent
Easily we can show ${\rm div}\,X\,>\,0$ (on $\widetilde C$).  
Therefore, since $X$ has the only zero point and $X$ is as in Fig.1 over a 
collar neighborhood of $\widetilde{\sigma}_i\,(i=1,2,3)$ 
by the proof of Theorem A, all the integral curves 
of $X$ through points other than its zero point converge to points of 
$\partial\widetilde C$ in finite time.  
Also, we can show ${\rm div}\,X^{\widetilde{\sigma}_i}\,>\,0$ on 
$\widetilde{\sigma}_i$ ($i=1,2,3$).  Hence all integral curves of 
$X^{\widetilde{\sigma}_i}$ through points other than its zero point converge 
to points of $\partial\widetilde{\sigma}_i$ in finite time.  
Therefore we obtain the following fact.  

\vspace{0.5truecm}

\noindent
{\bf Proposition 4.6.} {\sl The mean curvature flow having any non-minimal 
principal orbit of the $SO(4)$-action 
(on $SU(4)/S(U(2)\times U(2))$) as initial data converges to some singular 
orbit in finite time.  Also, the mean curvature flow having any non-minimal 
singular orbit through $\exp^{\perp}(\widetilde{\sigma}_i)$ of the 
$SO(4)$-action as initial data converges to one of two singular 
orbits through $\exp^{\perp}(\partial\widetilde{\sigma}_i)$ in finite time.}


\vspace{1truecm}

\centerline{{\bf References}}

\vspace{0.5truecm}

{\small 
\noindent
[BV] J. Berndt and L. Vanhecke, 
Curvature-adapted submanifolds, 
Nihonkai Math. J. 

{\bf 3} (1992) 177-185.

\noindent
[BCO] J. Berndt, S. Console and C. Olmos, Submanifolds and holonomy, 
Research Notes 

in Mathematics 434, CHAPMAN $\&$ HALL/CRC Press, Boca Raton, London, New York 

Washington, 2003.

\noindent
[Ch] U. Christ, 
Homogeneity of equifocal submanifolds, J. Differential Geometry {\bf 62} 

(2002) 1-15.

\noindent
[Co] L. Conlon, 
Remarks on commuting involutions, 
Proc. Amer. Math. Soc. {\bf 22} (1969) 

255-257.

\noindent
[GT] O. Goertsches and G. Thorbergsson, 
On the Geometry of the orbits of Hermann actions, 

Geom. Dedicata {\bf 129} (2007) 101-118.

\noindent
[H] S. Helgason, 
Differential geometry, Lie groups and symmetric spaces, Academic Press, 

New York, 1978.

\noindent
[HLO] E, Heintze, X. Liu and C. Olmos, Isoparametric submanifolds and a 
Chevalley-

type restriction theorem, Integrable systems, geometry, and topology, 151-190, 

AMS/IP Stud. Adv. Math. 36, Amer. Math. Soc., Providence, RI, 2006.

\noindent
[HPTT] E. Heintze, R.S. Palais, C.L. Terng and G. Thorbergsson, 
Hyperpolar actions 

on symmetric spaces, Geometry, topology and physics for Raoul Bott 
(ed. S. T. Yau), 

Conf. Proc. Lecture Notes Geom. Topology {\bf 4}, Internat. Press, Cambridge, 
MA, 1995 

pp214-245.

\noindent
[Koi1] N. Koike, On proper Fredholm submanifolds in a Hilbert space arising 
from subma-

nifolds in a symmetric space, Japan. J. Math. {\bf 28} (2002) 61-80.

\noindent
[Koi2] N. Koike, Actions of Hermann type and proper complex equifocal 
submanifolds, 

Osaka J. Math. {\bf 42} (2005) 599-611.

\noindent
[Kol] A. Kollross, A Classification of hyperpolar and cohomogeneity one 
actions, Trans. 

Amer. Math. Soc. {\bf 354} (2001) 571-612.

\noindent
[LT] X. Liu and C. L. Terng, The mean curvature flow for isoparametric 
submanifolds, 

Duke Math. J. {\bf 147} (2009) 157-179.

\noindent
[P] R.S. Palais, 
Morse theory on Hilbert manifolds, Topology {\bf 2} (1963) 299-340.

\noindent
[PT] R.S. Palais and C.L. Terng, Critical point theory and submanifold 
geometry, Lecture 

Notes in Math. {\bf 1353}, Springer, Berlin, 1988.

\noindent
[T1] C.L. Terng, 
Isoparametric submanifolds and their Coxeter groups, 
J. Differential Ge-

ometry {\bf 21} (1985) 79-107.

\noindent
[T2] C.L. Terng, 
Proper Fredholm submanifolds of Hilbert space, 
J. Differential Geometry 

{\bf 29} (1989) 9-47.

\noindent
[TT] C.L. Terng and G. Thorbergsson, 
Submanifold geometry in symmetric spaces, J. Diff-

erential Geometry {\bf 42} (1995) 665-718.

\noindent
[Z] X. P. Zhu, Lectures on mean curvature flows, Studies in Advanced Math., 
AMS/IP, 2002.

\vspace{0.5truecm}

{\small 
\rightline{Department of Mathematics, Faculty of Science}
\rightline{Tokyo University of Science, 26 Wakamiya-cho}
\rightline{Shinjuku-ku, Tokyo 162-8601 Japan}
\rightline{(koike@ma.kagu.tus.ac.jp)}
}

\end{document}